\documentclass[11pt]{article}
\usepackage[english]{babel}
\usepackage[cp850]{inputenc}
\usepackage[dvips]{graphicx}
\usepackage[left=2.5cm,right=2cm,top=3cm,bottom=3cm]{geometry}
\usepackage{amsmath,amsfonts,amsthm,amssymb}
\usepackage[usenames, dvipsnames]{color}
\usepackage{fancyhdr}
\usepackage{stmaryrd}
\usepackage[colorlinks=true,citecolor=red,linkcolor=blue,urlcolor=blue,pdfstartview=FitH]{hyperref}
\usepackage{dsfont}
\usepackage{xcolor}
\usepackage{epsfig}

\font\tenmath=msbm10 scaled 1200

\font\sevenmath=msbm7 scaled 1200
\font\fivemath=msbm5 scaled 1200 

\newfam\mathfam \textfont\mathfam=\tenmath
\scriptfont\mathfam=\sevenmath \scriptscriptfont\mathfam=\fivemath

\def\R{{\mathbb{R}}}
\def\N{{\mathbb{N}}}
\def\E{{\mathbb{E}}}
\def\L{{\cal L}}
\def\A{{\cal A}}
\def\P{{\mathbb{P}}}

\def\F{{\cal F}}
\newcommand {\lessim} {\ {\raise-.5ex\hbox{$\buildrel<\over\sim$}}\ }

\newtheorem{Definition}{Definition}[section]

\newtheorem{theo}{Theorem}[section]
\newtheorem{lem}{Lemma}[section]
\newtheorem{prop}{Proposition}[section]

\newtheorem{pro}{Proof}[section]

\newfam\mathfam \textfont\mathfam=\tenmath
\scriptfont\mathfam=\sevenmath \scriptscriptfont\mathfam=\fivemath

\def \^#1{\if#1i{\accent"5E\i}\else{\accent"5E#1}\fi}

\def \g{\gamma}

\def \cadre{(\Omega, {\cal A}, \P)}
\def \ds{\displaystyle}
\def \cqfd{\quad\Box}
\def \ms{\medskip}

\def \bs{\bigskip}

\def \Tr{\hbox{\bf w}}

\title{\bf Nonlinear Randomized Urn Models: a Stochastic Approximation Viewpoint}
 
\author{\textsc{Sophie Laruelle} \thanks{Universit\'e Paris-Est, Laboratoire d'Analyse et de Math\'ematiques Appliqu\'ees, UMR 8050, UPEMLV, UPEC, CNRS, F-94010, Cr\'eteil, France E-mail: \texttt{sophie.laruelle@u-pec.fr}. Part of this work was achieved   at  Laboratoire de Math\'ematiques Appliqu\'ees aux Syst\`emes, \'Ecole Centrale de Paris, Grande voie des vignes, 92295 Ch\^atenay-Malabry Cedex, France.} \and \textsc{Gilles Pag\`es} \thanks{Laboratoire de Probabilit\'es, Statistique et Mod\'elisation, UMR~8002, campus Pierre et Marie Curie, Sorbonne-Universit\'e, F-75252 Paris, France. E-mail: \texttt{gilles.pages@sorbonne-universite.fr}}}

\begin{document}

\maketitle 

\begin{abstract}
This paper extends the link between stochastic approximation ($SA$) theory and randomized urn models developed in~\cite{LarPag}, and their applications to clinical trials introduced in~\cite{BaiHu,BaiHu2,BaiHuShe}. We no longer assume that  the drawing rule is  uniform among the balls of  the urn (which contains $d$ colors), but can be reinforced by a function $f$. This  is a way to  model  risk aversion. Firstly, by considering that $f$ is concave or convex and by reformulating the dynamics of the urn composition as an $SA$ algorithm with remainder, we derive the $a.s.$ convergence and the asymptotic normality (Central Limit Theorem, $CLT$) of the normalized procedure by calling upon the so-called $ODE$ and $SDE$ methods. An in-depth analysis of  the case $d=2$ exhibits two different behaviors: A single equilibrium point when $f$ is concave, and when $f$ is convex, a transition phase from  a single attracting equilibrium to a system with two attracting and one repulsive equilibrium points. The last setting is solved using results on non-convergence toward noisy and noiseless ``traps"   in order to  deduce the $a.s.$ convergence  toward one of the attracting points. Secondly, the special case of a P\'olya urn (when the addition rule is the $I_d$ matrix) is analyzed, still using result from $SA$ theory about ``traps''. Finally, these results are applied to a function with regular variation and to an optimal asset allocation in Finance.
\end{abstract}


\bs 
\noindent {\em 2010 AMS classification:} 62L20, 62E20, 62L05
secondary: 62F12, 62P10.

\section{Introduction}


In this paper, we introduce and study in-depth a class of generalized P\'olya urns (with $d$ colors) characterized by their nonlinear drawing rules. These models appear as a generalization of  randomized urn models originally devised for clinical trials which takes into account the risk aversion attitude of the agent. Randomized urn models have been extensively investigated by various authors (see~\cite{BaiHu,BaiHu2,BaiHuShe}) during the last twenty years based on {\em ad hoc} martingale arguments to solve $a.s.$ convergence as well as its rate of convergence. In a recent paper~\cite{LarPag} (see also~\cite{LarPag2, Zhang}), we revisited, unified and often extended these results by showing that they can be established by relying on the main results of  Stochastic Approximation ($SA$) theory, especially $a.s.$ convergence and weak convergence rate (Central Limit Theorem ($CLT$)). Although this analysis is more  demanding than for urn models with linear drawing rules,  this $SA$ ``toolbox''  turns out to be very efficient  (see also~\cite{Ben2} and the references therein). $SA$ deals with the asymptotic behavior of zero search stochastic recursive procedures and goes back to the seminal paper by Robbins \& Monro in the 1950's. Since then, this theory has been developed extensively by many authors (see~\cite{KusCla, KusYin, BMP, Duf, Duf2} and the references therein for an overview and historical notes) and has been applied in various directions (Automatic Control, Mathematical Psychology, Artificial Neural Networks, Statistics, Stochastic Control, Numerical Probability, etc.). 

Considering nonlinear drawing rules leads, once written as a recursive stochastic algorithm, to dynamics for the (normalized) urn composition having several {\em equilibrium points}. By equilibrium point, which belongs to the $SA$ terminology adopted in the sequel, we mean a zero of the mean field of  the stochastic algorithm which in practice corresponds to a potential asymptotic composition of the urn. These equilibrium  points include not only local attractors, but also ``parasitic'' ones (repeller, saddle points, etc). This is a major difference with the linear case investigated in~\cite{LarPag} since, this time,  we will need to call upon the whole machinery of $SA$, in particular ``second order'' results about these parasitic noisy and noiseless  equilibrium points, sometimes called ``traps'' in the $SA$ literature (see~\cite{BraDuf, Duf,LamPagTar}, see also~\cite{Laz, Pem,Ben}). Taking advantage of these results, we will establish the $a.s.$ convergence, or strong consistency, of the (normalized) urn composition even in presence of multiple attractors. Then, we will analyze its weak convergence rate.

\smallskip
Let us be more precise on the urn model under consideration in this paper. We consider an urn containing balls of (at most) $d$ different types (or colors). All random variables involved in the model are supposed to be defined on the same probability space $\cadre$. Denote by $Y_0=(Y_0^i)_{i=1,\ldots,d}\in\R_+^d\setminus\{0\}$ the initial composition of the urn, where $Y_0^i$ is the number of balls of type $i\in\{1,\ldots,d\}$ (of course a more natural, though not mandatory, assumption would be $Y_0\!\in \N^d\setminus\{0\}$). The urn composition at draw $n$ is denoted by $Y_n=(Y_n^i)_{i=1,\ldots,d}$. At the $n^{th}$ stage, one draws randomly (according to a law defined further on) a ball from the urn with instant replacement. If the drawn ball is of type $j$, then the urn composition is updated by adding $D^{ij}_n$ balls of type $i$, for every  $i\in\{1,\ldots,d\}$. The procedure is then iterated. The urn composition at stage $n$, modeled by an $\R^d$-valued vector $Y_n$, satisfies the following recursive updating rule  between times $n$ and $n+1$:
\begin{equation}\label{dynamic}
	Y_{n+1}=Y_{n}+D_{n+1}X_{n+1},\quad n\geq0, \quad Y_0\!\in \R_+^d\setminus\{0\}, 
\end{equation}
where $D_n=(D^{ij}_n)_{1\leq i,j\leq d}$ is the {\em addition rule matrix}  and  $X_n:(\Omega,\A,\P)\rightarrow\{e^1,\cdots,e^d\}$ models the  type of the drawn ball  at time $n$ ($\{e^1,\cdots,e^d\}$ denotes  the canonical basis of $\R^d$ with $e^j$ standing for type $j$). We assume that there is no extinction {\it i.e.} $Y_n\!\in \R_+^d\setminus\{0\}$ $a.s.$ for every $n\ge 1$: This  is always  the case if all the entries $D_n^{ij}$ are $a.s.$ non-negative (see~\cite{LarPag}). The filtration of the model is defined by $\F_n=\sigma(Y_0,X_k,D_k,1\leq k\leq n)$, $n\geq0$. 

The {\em generating matrices} are defined as the ${\cal F}_n$-compensator of the additions rule sequence {\it i.e.} 
\begin{equation}\label{GenMatrix}
H_n=\left[\E\left(D^{ij}_n\,|\,\F_{n-1}\right)\right]_{1\leq i,j\leq d}, \; n\geq1.
\end{equation}
We will also assume that the sequence of generating matrices $a.s.$ converges toward a limiting generating matrix denoted by $H$.

\medskip
We will  mostly investigate  the asymptotic behavior of skewed drawing rules of the  form
%
	\begin{equation}\label{ConstructX}
		\forall i\in\{1,\ldots,d\}, \quad \P\big(X_{n+1}=e^i\big | \F_n\big)=\frac{f\big(Y_n^i/(n+\Tr(Y_0))\big)}{\sum_{j=1}^df(Y_n^j/(n+\Tr\big(Y_0))\big)}, \quad n\ge 0,
	\end{equation}
	where  $\Tr(y)=y^1+\cdots+y^d$ denotes the {\em weight} of a vector $y=(y^1,\ldots,y^d)^t\!\in \R_+^d$ and  the function
	\begin{equation}\label{eq:fX}
	\mbox{$f:\R_+\to\R_+$ is  non-decreasing    with $f(0)=0$, $f(1)=1$.}
	\end{equation}
The function $f$  usually satisfies an {\em additional convexity or concavity property}. Such a rule will be called  {\em normalized $f$-skewed empirical frequency based drawing rule}, $f$ being the {\em skewing function}, non-trivial  when $f\neq\operatorname{Id}_{\R_+}$. When $f =\operatorname{Id}_{\R_+}$, we retrieve a more standard drawing rule based on the regular empirical frequency of the types in the urn (see~\cite{BaiHu, LarPag} among others). 
	
\smallskip 
In fact, we will see that the result obtained for this family of frequency-based drawing rules allows us to elucidate a second way to skew the drawing, called {\em normalized $f$-skewed distribution}, based this time on the {\em number} of  balls of each type in the urn, namely 
\begin{equation}\label{ConstructX2}
		\forall i\in\{1,\ldots,d\}, \quad \P\big(X_{n+1}=e^i\,\big| \,\F_n\big)=\frac{f(Y_n^i)}{\sum_{j=1}^df(Y_n^j)}, \quad n\ge 0,
\end{equation}
where the function
\begin{equation}\label{eq:fX2}
	\mbox{$f:\R_+\to\R_+$ is non-decreasing,  with $f(0)=0$ and is  regularly varying with index  $\alpha>0$}
\end{equation}
({\it i.e.} for every $ t>0$, $\frac{f(t x)}{f(x)}\underset{x\to+\infty}{\longrightarrow}t^{\alpha}$).

Moreover, we will make the assumption that $D_n$ and $X_n$ are conditionally independent given ${\cal F}_{n-1}$ (see~${\bf (A2)}$ further on). Such a drawing procedure can be  performed by using an exogenous  i.i.d. sequence $(U_n)_{n\ge 1}$ of random variables with uniform distribution on the unit interval, independent of the sequence $(D_n)_{n \ge 1}$, to simulate the above conditional probabilities. 

\ms
Let us remark that, when $f={\rm Id}_{\R_+}$, both updating rules~\eqref{ConstructX} and~\eqref{ConstructX2} coincide. In fact, the normalized $f$-skewed distribution drawing rule  will appear as a by-product of the first one (see Section~\ref{sec:regvar}) by noting that, if $f$ is bounded on every interval $(0,M]$ and regularly varying with index $\alpha>0$, then $\frac{f(tx)}{f(x)}\underset{x\to\infty}{\longrightarrow}t^{\alpha}$ uniformly in $t$ on every interval $(0,T]$, $0<T<+\infty$ (see Theorem 1.5.2 p.22 in~\cite{BinGolTeu}). Thus, if $\frac{Y_n}{n+\sum_{i=1}^dY_0^i}$ lies in a compact set, then
$$
\max_{1\leq i\leq d}\left|\frac{f(Y_n^i)}{f(n+\sum_{i=1}^dY_0^i)}-\left(\frac{Y_n}{n+\sum_{i=1}^dY_0^i}\right)^{\alpha}\right|\underset{n\to+\infty}{\longrightarrow}0.
$$
Then, we will conclude by applying the result related to the $f$-skewed empirical frequency based drawing rule  to the   functions $x\mapsto x^{\alpha}$, $\alpha>0$. 

\medskip

In this paper we both randomize a single urn in the sense that the addition rule matrix (defined by~(\ref{GenMatrix})) itself can be random (as introduced in~\cite{Jan,BaiHu2,BaiHuShe} and studied with $SA$ theory in~\cite{LarPag}) and investigate wide non-parametric classes of convex and concave drawing rules (see~\eqref{ConstructX} and~\eqref{ConstructX2}). $A.s.$ convergence of the urn composition and that of the drawing rule is proved as well as their convergence  rate, either in a weak or strong sense, depending on the structure of the updating rules. In the particular case of  P\'olya's urns ({\it i.e.} when the addition rule matrix $D_n$ is equal to identity), but implemented here with a convex skewed drawing rule (for generalized P\'olya urn, see for example~\cite{HilLanSud,Zhu,ColCotLic}), ``noiseless traps'' may appear ({\it i.e.} unstable equilibria, noiseless since they lie  at the boundary of the state space). To determine whether they are parasitic, we  develop a dedicated approach, close in spirit to that introduced in~\cite{LamPagTar} and~\cite{LamPag2} the analysis of adaptive bandit algorithms where the authors establish  a kind of ``oracle'' inequality relying on a specified martingale (see Lemma~\ref{LemMgle} further on). These results highlight the efficiency of $SA$ theory, even in presence of noiseless repulsive equilibrium points.
 
Recently, a system of P\'olya's urns with graph based interactions and a ``power drawing'' rule has also been investigated using $SA$ techniques in~\cite{BenBenCheLim,CheLuc} and the $a.s.$ convergence of the normalized urn composition is established.

Generalized P\'olya Urn models ($GPU$) have been widely studied in the literature with different points of view: Martingale method (see $e.g.$~\cite{Fre}), algebraic approach (see $e.g.$~\cite{Pou}), reinforcement process (see $e.g.$~\cite{Pou}), branching process (see $e.g.$~\cite{Jan}), stochastic approximation (see for example~\cite{Ben2,LasMaiSel}), contraction method (see~\cite{KnaNei}). These models also have applications to many areas: Biology, random walks and clinical trials, statistics and learning, computer science,  psychology, economics or finance for instance (see~\cite{Zhu}). 

In these adaptive models, the key point is the  updating rules  of the urn composition  after each drawing given here by~\eqref{ConstructX} and~\eqref{ConstructX2}. Basically, we will show that (a normalized version of) this urn composition can be formulated  as a classical recursive stochastic algorithm with step $\gamma_n=\frac{1}{n+\Tr(Y_0)}$ where $\Tr(Y_0)$ denotes the number of balls in the urn at time $0$. Doing so, we will be in position to first establish the $a.s.$ convergence of the procedure by calling upon the so-called Ordinary Differential Equation method ($ODE$ method) toward a finite set of equilibrium points (but usually not reduced to a single point). As a second step, we will rely on $a.s.$ non-convergence results toward traps (see~\cite{BraDuf,Duf}) and on $a.s.$ convergence in presence of multiple targets (see~\cite{BenHir,Duf, ForPag}). As a third step, we entirely elucidate the rate of convergence (namely a weak rate through a $CLT$ or an $a.s.$ rate) by using the Stochastic Differential Equation method ($SDE$ method, see $e.g.$~\cite{Duf2,BMP}). The three main theoretical results from $SA$ are recalled in a self-contained form in the  Appendix. Proofs of such results can be found in classical textbooks on $SA$ (\cite{BMP,Duf,Duf2,KusYin}). As for the $CLT$, they go back to~\cite{KusCla} and~\cite{Bou}, see also~\cite{Zhang} more recent results. We will verify once again  how powerful these general theorems are  to solve such  questions, sparing tedious computations and  repetitive proofs.

Among many  fields of application, we present in Section~\ref{cinq} an adaptive asset allocation procedure based on a reinforcement principle relying on non-linear randomized urns, illustrated by a first numerical test.  One may also consider a similar procedure  as a strategy to update the composition of a portfolio or even a whole fund, based on the (recent) past performances of the assets.

\medskip
The paper is organized as follows. Section~\ref{deux} presents the framework of skewed randomized urn models with the required assumptions on both the addition rule matrices and the generating matrices. After rewriting the dynamics of the urn composition as an $SA$ procedure in Section~\ref{sec:sa}, we analyze in Section~\ref{sec:exist} the equilibrium points and their stability for the associated $ODE$ when the $f$-skewed drawing rule is convex/concave. An in-depth analysis of the $2$-color urn is carried out  in Section~\ref{trois}. We exhibit several kinds of behaviors: When $f$ is concave, there is always a unique stable equilibrium point and when $f$ is convex three generic situations may occur with one, tow or three equilibrium points (one being parasitic in the last two settings).
By calling upon  $SA$ result on traps, we prove the $a.s.$ convergence towards one of the attracting equilibrium points; then we derive from the $SDE$ method all the possible rates of convergence. In Section~\ref{quatre}, we study the case of  P\'olya urns~--~urn dynamics whose  addition rule matrix equals to identity~--~updated by a skewed drawing rule. We rely on  methods  borrowed from the analysis of  adaptive bandit algorithms to prove the convergence towards the ``targeted urn composition" and the non-convergence towards traps. Finally, in Section~\ref{cinq}, we first transfer our results to the second type of drawing rule and we conclude  by an   application   to portfolio allocation.

\medskip
\noindent {\sc Notations.} For $u=(u^i)_{i=1,\ldots,d}\in\R^d$, $\left(\cdot\left|\right.\cdot\right)$ denote the Euclidean inner product and $\left\| u\right\|$ its related norm of the column vector $u\!\in\R^d$, $\Tr(u)=\sum_{k=1}^du^k$ denotes  its ``weight'',  $u^t$ denotes its transpose, $u\otimes v = [u^iv^j]_{i,j=1,\ldots,d}$; $|\!|\!|A|\!|\!|$ denotes the operator norm of the matrix $A\in{\cal M}_{d,q}(\R)$ with $d$ rows and $q$ columns with respect to the two canonical Euclidean norms. When $d\!=\!q$, ${\rm Sp}(A)$ denotes the set of eigenvalues of $A$. $\mathbf{1}\!=\! (1\cdots1)^t$ denotes the unit column vector in $\R^d$, $I_d$ denotes the $d\times d$ identity matrix, ${\rm diag}(u)=[\delta_{ij}u_i]_{1\le i,j\le d}$, where $\delta_{ij}$ stands for the Kronecker symbol and ${\cal S}_d=\left\{u\in\R_+^d:\sum_{i=1}^du^i=1\right\}$ denotes the canonical simplex. ${\cal E}_d=\{y\in{\cal S}_d: h(y)=0\}$ denotes the set of zeros of $h$ called  equilibrium points in reference to the $ODE$ $\dot y = -h(y)$.
 
\section{Skewed randomized urn models}\label{deux}

\subsection{Main assumptions and definitions}

With the notations and definitions described in the introduction, we are in position to formulate the main assumptions needed to establish the $a.s.$ convergence of the urn composition. 

\medskip
\noindent ${\bf (A1)}$ 
$\equiv\left\{\begin{array}{ll}
	(i) &  \hskip-0.2cm\mbox{{\em Addition rule matrix:}  For every $n\ge 1$, the  matrix $D_n$ $a.s.$ has non-negative entries.} \\
	
(ii) &  \hskip-0.2cm\mbox{{\em Generating matrix:}  For every $n\ge 1$,  the generating matrix  $H_n=(H^{ij}_n)_{1\leq i,j \leq d}$} \\
	    & \hskip-0.2cm\mbox{$a.s.$ satisfies:} \hskip2cm\forall\,j \in\{1,\ldots,d\}, \quad \displaystyle\sum_{i=1}^d H^{ij}_n=c>0. \\
	   
	(iii) & \hskip-0.2cm\mbox{{\em Starting value:} The starting urn composition vector $Y_0\!\in \R_+^d\setminus\{0\}$.}
\end{array}\right.$

\bigskip
The constant $c$ is known as the {\em balance} of the urn. In fact, we may assume without loss of generality, up to a renormalization of $Y_n$, that $c=1$. As a matter of fact, we set for every $n\geq0$, $\widehat{Y}_n=\frac{Y_n}{c}$ and $\widehat{D}_{n+1}=\frac{D_{n+1}}{c}$ where $(Y_n, X_n,D_n)$ satisfies~\eqref{dynamic}, then the couple $(\widehat{Y}_n,X_n,\widehat{D}_{n})_{n\geq 1}$, still satisfies the dynamics~(\ref{dynamic}), namely
$$
\widehat{Y}_{n+1}=\widehat{Y}_{n}+\widehat{D}_{n+1}X_{n+1},\quad n\geq0, \quad \widehat{Y}_0\!\in \R_+^d\setminus\{0\},
$$
whereas $\widehat{H}_n=\E[\widehat{D}_n|\F_{n-1}]$ satisfies now ${\bf (A1)}$-$(iii)$ with $c=1$ {\it i.e.} $\widehat{H}_n$ is {\em co-stochastic} (in the sense that its transpose is a stochastic matrix). Assumptions ${\bf (A1)}$-$(i)$\&$(iii)$ combined with the drawing rule~\eqref{dynamic} ensure that $Y_n\in\R^d_+\setminus\{0\}$, for every $n\geq0$.

From now on, throughout the paper, we will consider this normalized balanced version, still denoted by $Y_n$ and $D_n$ for convenience.

\medskip
\noindent ${\bf (A2)}$ $\equiv\left\{\begin{array}{ll}
(i) &\hskip-0.25cm   \mbox{The addition rule $D_n$ and the drawing procedure $X_n$ are conditionally independent} \\
    & \mbox{given $\F_{n-1}$.} \\
(ii) &
	\forall j\in\{1,\ldots,d\}, \quad\sup_{n\geq1}\E\left[\left\|D^{\cdot j}_n\right\|^{2}\,|\,\F_{n-1}\right]<+\infty\quad a.s. \label{A2} \\
	&\hskip 3 cm \Longleftrightarrow \forall\, i,j\in\{1,\ldots,d\}, \quad \sup_{n\geq1}\E\left[(D_n^{ij})^{2}\,|\,\F_{n-1}\right]<+\infty  \quad a.s. \nonumber
	\end{array}\right.$ 
	
\medskip
\noindent  where $D^{\cdot\,j}_n=(D^{ij}_n)_{i=1,\ldots,d}$ (column vector). 

\bigskip
\noindent ${\bf (A3)}$ There exists an irreducible $d\times d$ matrix $H$ (with non-negative entries) such that 
\begin{equation}\label{A3}
	H_n\overset{a.s.}{\underset{n\rightarrow+\infty}{\longrightarrow}} H \quad \mbox{and} \quad \sum_{n\geq1}|\!|\!|H_n-H|\!|\!|^2<+\infty \quad a.s.
\end{equation}
$H$ is called the {\em limiting generating matrix}. 

\bigskip 
The central object of interest of this paper will be the {\em  (quasi-)normalized composition} of the urn at time $n$ defined by
\begin{equation}\label{eq:tildeYn}
\widetilde Y_n = \frac{Y_n}{n+\Tr(Y_0)}
\end{equation}
(where the weight function is defined in the introduction). The reason for introducing  such a  renormalization factor   is that, as established further on in the next subsection,
$$
\forall n\geq0, \quad \E\left[\Tr(Y_n)\right]=n+\Tr(Y_0).
$$
Then, one may guess that  $\displaystyle\widetilde{Y}_n$ is close to the simplex ${\cal S}_d=\left\{u\!\in\R_+^d: \sum_{i=1}^du^i=1\right\}$ and will possibly asymptotically lie in it. Even note that when the matrices $D_n$ are themselves co-stochastic, $\widetilde Y_n$ is ${\cal S}_d$-valued. Therefore, this is a natural deterministic way to normalize the urn composition vector.


\begin{Definition}[Skewing  functions and skewed drawing rules] $(a)$ A function $f:\R_+\to \R_+$ satisfying
\begin{equation}\label{eq:def-f}
\mbox{$f$ non-decreasing,  {\em convex or concave}, $f(0)=0$ and $f(1)=1$ and $\{f>0\}= (0,+\infty)$}
\end{equation}
is called a {\em skewing function}. 

\medskip
\noindent $(b)$ Assuming $({\bf A1})$-$(i)$ \& $(iii)$, the {\em $f$-skewed drawing rule} $(X_n)_{n\ge 1}$ induced by a skewing function $f$ is defined by 
\begin{equation}\label{LoiXFConvex}
\forall\, i\in\{1,\ldots,d\}, \quad \P\big(X_{n+1}=e^i\big | \F_n\big)=\frac{f(\widetilde{Y}_n^i)}{\sum_{j=1}^df(\widetilde{Y}_n^j)}, \quad n\ge 0,
\end{equation}
where $\displaystyle\widetilde{Y}_n$ is defined by~\eqref{eq:tildeYn}, and $(Y_n,X_n, D_n)_{n\ge 1}$ satisfies~\eqref{dynamic}. 
\end{Definition}

Of course, such a drawing rule is really {\em skewed} only if $f\not\equiv\operatorname{Id}_{\R_+}$. Note that this definition is consistent under $({\bf A1})$-$(i)$ \& $(iii)$ since $Y_n \!\in \R_+^d\setminus\{0\}$, for every $n\ge 0$ so that $\widetilde Y_n\!\in \R_+^d\setminus \{0\}$. 

\smallskip 
We will extensively use that a skewing function $f$ always satisfies that the function $\xi\in\R_+\mapsto \frac{f(\xi)}{\xi}$ is monotonic (non-decreasing if $f$ is convex, non-increasing if $f$ is concave) and strictly monotonic if the concavity/convexity of $f$ is itself strict.

\subsection{Representation as a stochastic algorithm}\label{sec:sa}

The starting point,  like in~\cite{LarPag},  is to reformulate the dynamics (\ref{dynamic})-(\ref{ConstructX}) as a recursive stochastic algorithm in order to take advantage of classical results from $SA$ Theory to elucidate the asymptotic properties ($a.s.$ convergence) of both the urn composition $Y_n$ and the ball drawing rule $X_n$. To do so, we start from~(\ref{dynamic}) with $Y_0\!\in\R_+^d\setminus\{0\}$. For every $n\geq0$, we note that
\begin{equation}\label{dynq2}
Y_{n+1}=Y_n+D_{n+1}X_{n+1}=Y_n+\E\left[D_{n+1}X_{n+1}\,|\,\F_{n}\right]+\Delta M_{n+1},
\end{equation}
where 
\begin{equation}\label{eq:DeltaMn}
\Delta M_{n+1}:=D_{n+1}X_{n+1}-\E\left[D_{n+1}X_{n+1}\,|\,\F_{n}\right]
\end{equation} 
is an $\F_{n}$-local martingale increment (integrability follows from ${\bf (A2)}$-$(ii)$). By the definition~\eqref{GenMatrix} of the generating matrix $H_n$, we have, owing to the conditional independence assumption ${\bf (A2)}$-$(i)$,
\begin{eqnarray*}
\E\left[D_{n+1}X_{n+1}\,|\,\F_{n}\right]&=&\sum_{i=1}^d\E\left[D_{n+1}\mathds{1}_{\left\{X_{n+1}=e^i\right\}}\,|\,\F_{n}\right]e^i=\sum_{i=1}^d\E\left[D_{n+1}\,|\,\F_{n}\right]\P\left(X_{n+1}=e^i\,|\,\F_{n}\right)e^i\\
 &=&H_{n+1}\sum_{i=1}^d\frac{f(\widetilde{Y}_n^i)}{\Tr(\widetilde{f}(\widetilde{Y}_n))}e^i=H_{n+1}\frac{\widetilde{f}(\widetilde{Y}_n)}{\Tr(\widetilde{f}(\widetilde{Y}_n))}
\end{eqnarray*}
where, for every $y=(y^1,\ldots,y^d)^t\!\in \R_+^d\setminus \{0\}$,  
\begin{equation}\label{eq:ftilde}
\widetilde{f}\big((y^1,\ldots,y^d)^t\big)=\left(f(y^i)\right)_{1\leq i\leq d}\!\in\R_+^d\setminus \{0\}
\end{equation} 
is a column vector, so that 
\begin{equation}\label{ASC1}
Y_{n+1}=Y_n+H_{n+1}\frac{\widetilde{f}(\widetilde{Y}_n)}{\Tr(\widetilde{f}(\widetilde{Y}_n))}+\Delta M_{n+1}.
\end{equation}
Now we can derive a stochastic approximation for the normalized urn composition $\displaystyle\widetilde Y_n=\frac{Y_n}{n+\Tr(Y_0)}$, $n\geq0$. First, we have, for every $n\ge 0$, 
\begin{equation}\label{ASC2}
\frac{Y_{n+1}}{n+1+\Tr(Y_0)}=\frac{Y_n}{n+\Tr(Y_0)}+\frac{1}{n+1+\Tr(Y_0)}\left(H_{n+1}\frac{\widetilde{f}(\widetilde{Y}_n)}{\Tr(\widetilde{f}(\widetilde{Y}_n))}-\frac{Y_n}{n+\Tr(Y_0)}\right)+\frac{\Delta M_{n+1}}{n+1+\Tr(Y_0)}.
\end{equation}
Consequently, the sequence  $(\widetilde{Y}_n)_{n\ge 0}$, satisfies the canonical recursive stochastic approximation procedure starting from $\widetilde Y_0\!\in \R_+^d\setminus \{0\}$, 
\begin{equation}
\widetilde{Y}_{n+1}= \widetilde{Y}_n+\frac{1}{n+1+\Tr(Y_0)}\left(H_{n+1}\frac{\widetilde{f}(\widetilde{Y}_n)}{\Tr(\widetilde{f}(\widetilde{Y}_n))}-\widetilde{Y}_n\right)+\frac{1}{n+1+\Tr(Y_0)}\Delta M_{n+1}
\end{equation}	 
or, equivalently, 
\begin{equation}
\widetilde{Y}_{n+1}=\widetilde{Y}_n-\gamma_{n+1}h(\widetilde Y_n)+\gamma_{n+1}\left(\Delta M_{n+1}+r_{n+1}\right) 
\label{ASNLConvex}
\end{equation}	   
where the {\em mean field} $h:\R^d\setminus\{0\}\to\R^d$ of the procedure is defined by
\begin{equation}\label{defhcvx}
	h(y):=\E\left[\left.\widetilde{Y}_n-H\frac{\widetilde{f}(\widetilde{Y}_n)}{\Tr(\widetilde{f}(\widetilde{Y}_n))}\,\right|\,\widetilde{Y}_n=y\right]=y-H\frac{\widetilde{f}(y)}{\Tr\big(\widetilde{f}(y)\big)},
\end{equation}
$\gamma_n :=\frac{1}{n+\Tr(Y_0)}$ is its step parameter and
\begin{equation}\label{resteConvex}
r_{n+1}:=(H_{n+1}-H)\frac{\widetilde{f}(\widetilde{Y}_n)}{\Tr(\widetilde{f}(\widetilde{Y}_n))}
\quad \mbox{ is an $\F_n$-adapted remainder term.}
\end{equation}

\subsection{Boundedness of the normalized urn composition}\label{sec:bound}

Our first task is to establish the $a.s.$ boundedness of the sequence $(\widetilde{Y}_n)_{n\geq0}$. 
By summing up the components of $\widetilde{Y}_n$ in~(\ref{ASC1}), we obtain
$$
\Tr(Y_{n+1})=\Tr(Y_n)+\frac{\Tr(H_{n+1}\widetilde{f}(\widetilde{Y}_n))}{\Tr(\widetilde{f}(\widetilde{Y}_n))}+\Tr(\Delta M_{n+1}).
$$
Using that the transpose of the generating matrix $H_{n+1}$ is a stochastic matrix by ${\bf (A1)}$-$(ii)$ (with $c=1$), we obtain
$$
\Tr(H_{n+1}\widetilde{f}(\widetilde{Y}_n))=\sum_{i=1}^d(H_{n+1}\widetilde{f}(\widetilde{Y}_n))_i=\sum_{i=1}^d\sum_{j=1}^dH^{ij}_{n+1}f(\widetilde{Y}_n^j)=\sum_{j=1}^d\left(\sum_{i=1}^dH^{ij}_{n+1}\right)f(\widetilde{Y}_n^j)=\Tr(\widetilde{f}(\widetilde{Y}_n)).
$$
Consequently, for every $n\geq0$,
\begin{equation}\label{ASTrace}
\Tr(Y_{n+1})=\Tr(Y_n)+1+\Tr(\Delta M_{n+1}).
\end{equation}

\noindent Let  $N_0=0$ and $N_n:=\sum_{k=1}^nX_k$, $n\geq1$, denote  the number of times each type of ball was drawn between draws 1 and $n$. For every $n\geq0$,
\begin{eqnarray*}N_{n+1}&=&N_n+X_{n+1}=N_n+\frac{\widetilde{f}(\widetilde{Y}_n)}{\Tr(\widetilde{f}(\widetilde{Y}_n))}+\Delta \widetilde{M}_{n+1},\\
\mbox{where }\qquad\Delta \widetilde{M}_{n+1}&:=&X_{n+1}-\E\left[X_{n+1}\left|\right.\F_n\right]=X_{n+1}-\frac{\widetilde{f}(\widetilde{Y}_n)}{\Tr(\widetilde{f}(\widetilde{Y}_n))}
\end{eqnarray*} 
is an $\F_n$-martingale increment. Thus, $\widetilde{N}_n:=\displaystyle \frac{N_n}{n}$ satisfies, still for every $n\ge 0$, 
$$
\widetilde{N}_{n+1}=\widetilde{N}_n-\frac{1}{n+1}\left(\widetilde{N}_n-\frac{\widetilde{f}(\widetilde{Y}_n)}{\Tr(\widetilde{f}(\widetilde{Y}_n))}\right)+\frac{1}{n+1}\Delta\widetilde{M}_{n+1}.
$$

\begin{prop}\label{ThmCvxTrace}
Let $(Y_n)_{n\geq0}$ be the urn composition sequence defined by~(\ref{dynamic})-(\ref{ConstructX}). 

\smallskip
\noindent $(a)$ Under the assumptions {\bf (A1)} and {\bf (A2)},
$$
\frac{\Tr(Y_n)}{n+\Tr(Y_0)}\overset{a.s.}{\underset{n\rightarrow+\infty}{\longrightarrow}} 1.
$$
\noindent \smallskip \noindent $(b)$ If the addition rule matrices $D_n$ themselves are co-stochastic, then $\Tr(Y_n)=n+\Tr(Y_0)$, and the sequence $(\widetilde{Y}_n)_{n\geq0}$ lives in the simplex ${\cal S}_d$.
\end{prop}

\noindent {\bf Proof.} $(a)$ We derive from the identity
$$
D_{n+1}X_{n+1}=\sum_{j=1}^dD^{\cdot\,j}_{n+1}\mathds{1}_{\{X_{n+1}=e^j\}}, \quad n\geq0,
$$
that
$$
\left\|D_{n+1}X_{n+1}\right\|^2=\sum_{j=1}^d\left\|D^{\cdot\, j}_{n+1}\right\|^2\mathds{1}_{\{X_{n+1}=e^j\}}.
$$
Hence, owing to ${\bf (A2)}$-$(ii)$,
\begin{eqnarray*}
\E\left[\left\|D_{n+1}X_{n+1}\right\|^2\,|\,\F_n\right]&=&\sum_{j=1}^d\E\left[\left\|D^{\cdot j}_{n+1}\right\|^2\,|\,\F_n\right]\P\left(X_{n+1}=e^j\,|\,\F_n\right) \\
	&\leq&\,\sup_{n\geq0}\sup_{1\leq j\leq d}\E\left[\left\|D^{\cdot j}_{n+1}\right\|^2\,|\,\F_n\right]<+\infty\quad a.s.
\end{eqnarray*}
Consequently $\sup_{n\geq1}\E\left[\left\|\Delta M_{n+1}\right\|^2\,|\,\F_n\right]<+\infty$ $a.s.$ and  thanks to the strong law of large numbers for conditionally $L^2$-bounded local martingale increments, we have $\frac{M_n}{n}\underset{n\rightarrow+\infty}{\longrightarrow}0$ $a.s.$ Finally, it follows from (\ref{ASTrace}) that
$$
\frac{\Tr(Y_n)}{n+\Tr(Y_0)}=1+\frac{\Tr(M_n)}{n+\Tr(Y_0)}\overset{a.s.}{\underset{n\rightarrow+\infty}{\longrightarrow}} 1.
$$

\noindent \smallskip \noindent $(b)$ In this case $\Tr(M_n)=0$, consequently for every $n\geq0$, $\Tr(\widetilde{Y}_n)=1$. \hfill$\cqfd$

\subsection{Existence of equilibrium points}\label{sec:exist}

As written in~(\ref{ASNLConvex}), the urn dynamics  appears  as a recursive   zero search algorithm with  mean field  $h:[0,1]^d\setminus\{0\} \to\R^d$ whose potential limiting points lies   in the  the  canonical simplex ${\cal S}_d$ defined by 
$$
{\cal S}_d=\Tr^{-1}\{1\}=\left\{y\in\R^d_+\,|\,\Tr(y)=1\right\}.
$$
More generally, since the components of $\widetilde{Y}_n=\frac{Y_n}{n+\Tr(Y_0)}$ are non-negative by construction and~--~under $({\bf A1})$-$({\bf A2})$~--~$\Tr(\widetilde{Y}_n)=\frac{\Tr(Y_n)}{n+\Tr(Y_0)}\underset{n\rightarrow+\infty}{\longrightarrow}1$ $a.s.$, it is clear that $\P(d\omega)$-$a.s.$, the sequence $(\widetilde{Y}_n(\omega))_{n\geq0}$ is bounded and that the set ${\cal Y}_{\infty}(\omega)$ of its limiting values lies in the simplex ${\cal S}_d$. Consequently, we search ${\cal S}_d$-valued equilibrium points  {\it i.e.} points $y\in{\cal S}_d$ such that $h(y)=0$ where $h$ is given by~\eqref{defhcvx}.

\medskip
Throughout this section $f $ denotes a skewing function in the sense of~\eqref{eq:def-f} and $H$ is a {\em deterministic}  (co-stochastic) matrix, intended to be a  limiting generating matrix as soon as it is irreducible.

%

\begin{prop}\label{prop:minmax} 
$(a)$ Let $H$ be a deterministic  co-stochastic matrix.  The function $\varphi_{_H}:[0,1]^d\setminus\{0\}\to {\cal S}_d$ defined by $\varphi_{_H}(y) =  H\frac{\widetilde{f}(y)}{\Tr\big(\widetilde{f}(y)\big)}$ 
has at least one fixed point. As a consequence $h$ has at least one zero $y^*$ and $\{h=0\}\subset {\cal S}_d$. (When $H$ is bi-stochastic, $y(d):= \frac 1d \mbox{\bf 1}$ is such a zero of $h$.)

\smallskip
\noindent 
$(b)$ If, furthermore, for every $i,j\in\{1,\ldots,d\}$, $H^{ij}>0$, then for every zero $y^*$ of $h$ in ${\cal S}_d$, 
$$
\forall\, i\in\{1,\ldots,d\},\quad \min_{1\le j\le d}H^{ij}\leq y^{*,i}\leq\max_{1\le j\le d} H^{ij}.
$$
\end{prop}

\noindent {\bf Proof.} $(a)$ The function $\varphi_{_H}$ is well-defined on $[0,1]^d\setminus\{0\}$ since $\Tr\big(\widetilde{f}(y)\big)>0$ on  $[0,1]^d\setminus\{0\}$ owing to the fact that $f>0$ on $(0,1)$. The function $\varphi:y\mapsto\frac{\widetilde{f}(y)}{\Tr(\widetilde{f}(y))}$ clearly maps  $[0,1]^d\setminus\{0\}$ into ${\cal S}_d$. So does $\varphi_{_H}$ since $H$ is co-stochastic and subsequently maps ${\cal S}_d$ into ${\cal S}_d$. Since $h(y)=y-\varphi_{_H}(y)$ by~\eqref{defhcvx}, it follows that $\{h=0\}\subset {\cal S}_d$. 

Now, both  functions $y\mapsto\widetilde{f}(y)$ and $y\mapsto\Tr\big(\widetilde{f}(y)\big)$ are continuous on ${\cal S}_d$. Moreover $\Tr\big(\widetilde{f}(y)\big)\geq f\left(\frac{1}{d}\right)>0$ since, for  any  $y=(y^1,\ldots,y^d)^t\!\in{\cal S}_d$, there exists $i_0\in\{1,\ldots,d\}$ such that $y^{i_0}\geq\frac{1}{d}$ so that $\Tr\big(\widetilde{f}(y)\big)\geq f(y^{i_0})\geq f\left(\frac{1}{d}\right)>0$. Therefore, $\varphi_{_H}$ is continuous and  maps  ${\cal S}_d$ into ${\cal S}_d$.  Then, by Brouwer's Theorem, $\varphi_H$ has  at least one fixed point {\it i.e.} $\{h=0\}\neq\varnothing $. The last claim is obvious since $\varphi(y(d))=y(d)$ and $H\mbox{\bf 1}= \mbox{\bf 1}$.

\smallskip
\noindent $(b)$ Let $i\!\in\{1,\ldots,d\}$.  It follows from the identity $\displaystyle \sum_{j=1}^dH^{ij}\frac{f(y^{*,j})}{\Tr(\widetilde{f}(y^*))}=y^{*,i}$,  that  $\displaystyle\min_{1\le j\le d}H^{ij}\!\le\! y^{*i}\!\le\! \max_{1\le j\le d} H^{ij}$ since $f$ is non-negative.~$\qquad\cqfd$ 

\begin{prop}[Bi-stochastic case]\label{prop:fixpt} Let ${\cal E}_d:=\{h=0\} = \big\{y\!\in \R_+^d\setminus\{0\}: h(y)=0\big\}\subset {\cal S}_d$ denote the (non-empty) set of equilibrium points.

\smallskip 
\noindent$(a)$ If $H$ is bi-stochastic, then $y(d):= \frac{1}{d}{\bf 1}\!\in{\cal E}_d$.

\smallskip
\noindent  $(b)$ If $H\!=\!I_d$, then $\big\{\widetilde e_{_I},\, I\subset\{1,\ldots,d\},I\neq\varnothing\big\}\subset{\cal E}_d$, where $\widetilde e_{_I}=  \frac{1}{|I|}\sum_{i\in I}e^i$ with $(e^i)_{1\leq i\leq d}$  the canonical basis of $\R^d$.

\smallskip
\noindent 	$(c)$ If $H\!=\!I_d$ and $f$ is a strictly convex or strictly concave skewing function, then
$$
{\cal E}_d=\big\{\widetilde e_{_I}, I\subset\{1,\ldots,d\},I\neq\varnothing\big\}.
$$
\noindent $(d)$ If   $H$ is  bi-stochastic and  irreducible, then ${\cal E}_d\subset\overset{\circ}{{\cal S}_d}=\left\{y\in(0,1)^d:\sum_{i=1}^dy^i=1\right\}$.
	
\smallskip
\noindent $(e)$ If,  furthermore,  $f$ is strictly concave, then ${\cal E}_d=\left\{\frac{1}{d}{\bf 1}\right\}$.
\end{prop}

\noindent {\bf Proof.} 
$(a)$ Set $y(d)=\frac{1}{d}{\bf 1}\in{\cal S}_d$. Thus, $\frac{f\left(\frac{1}{d}\right)}{d\cdot f(\frac{1}{d})}=\frac{1}{d}$ since $f\left(\frac{1}{d}\right)>0$ and, consequently, $\displaystyle \sum_{j=1}^dH^{ij}\frac{1}{d}=1\times\frac{1}{d}$ for every $i\!\in\{1,\ldots,d\}$, 
	therefore $h\big(y(d)\big)=y(d)-y(d)=0$.

\smallskip
\noindent \smallskip \noindent $(b)$ Let $I\neq\varnothing$. Then
$$
f(\widetilde e_{_I}^i)=\left\{\begin{array}{ccc} 0 & \mbox{if} & i\notin I \\ f\left(\frac{1}{|I|}\right) >0& \mbox{if} & i\in I \end{array}\right..
$$
Hence,  $\Tr(\widetilde{f}(\widetilde e_{_I}))=|I|f\left(\frac{1}{|I|}\right)$ as  well. Consequently $h(\widetilde  e_{_I})=0$.
	
\smallskip
\noindent $(c)$ Let $y^*\in{\cal E}_d$. Assume that there exists $y^{*,i_0}$, $y^{*,i_1}$ such that $0<y^{*,i_0}<y^{*,i_1}\leq1$. Then $y^{*,i_0}=\frac{f(y^{*,i_0})}{\Tr(\widetilde{f}(y^*))}$ and $y^{*,i_1} = \frac{f(y^{*,i_1})}{\Tr(\widetilde{f}(y^*))}$ so that $\frac{f(y^{*,i_0})}{y^{*,i_0}}=\frac{f(y^{*,i_1})}{y^{*,i_1}}=\Tr(\widetilde{f}(y^*))$. Now, if $  f$ or $-f$  is strictly convex, then the function $\xi\in\R_+\mapsto\frac{f(\xi)}{\xi}$ is strictly monotonic since $f(0)=0$. This yields a contradiction. 

As a consequence, there exists $\xi_0\!\in (0,1]$ such that $y^{*,i}\in\{0,\xi_0\}$, $i=1,\ldots,d$. Consequently,  if $I_{\xi_0}=\{i:y^{*,i}=\xi_0\}$, $\Tr(\widetilde{f}(y^*))=|I_{\xi_0}|f(\xi_0)$ and, for every $i\in I_{\xi_0}$,  $\frac{f(y^{*,i})}{\Tr(\widetilde{f}(y^*))}=\frac{f(\xi_0)}{|I_{\xi_0}|f(\xi_0)}=\frac{1}{|I_{\xi_0}|}$, {\it i.e.} $h(y^*)=\frac{1}{|I_{\xi_0}|}\sum_{i\in I_{\xi_0}}e^i$.  Hence, $y^*=\frac{1}{|I_{\xi_0}|}\sum_{i\in I_{\xi_0}}e^i$ so that $\xi_0=1$ and $y^*\!\in \big\{\tilde e_I,\, I\subset\{1,\ldots,d\}, \, I\neq \varnothing\big\}$.

\smallskip 
\noindent $(d)$ Let $i_0$ be such that $y^{*,i_0}=\min_iy^{*,i}$. If $y^*\notin\overset{\circ}{{\cal S}_d}$, then  $y^{*,i_0}=0$ so that
$$
\sum_{j=1}^dH^{i_0j}\frac{f(y^{*,j})}{\Tr(\widetilde{f}(y^*))}=0
$$
and $i_0\!\in I_0^*=\{i:  y^{*,i}=0\}\neq\varnothing$. Let $I_{>0}^*=I\setminus I_0^*= \{j: y^{*,j}>0\}\neq\varnothing$ since $y^*\!\in {\cal S}_d$. Let us show by induction  that  $I_{>0}^*= \{j: y^{*,j}>0\}$ and $I_0^*=\{j= y^{*,i}_j=0\}$ are not connected by any power of $H$ which will contradict the irreducibility of $H$. 
Let $i_0\!\in I^*_0$ and  $j_0\!\in I_{>0}^*$, then the above equality implies $H^{i_0j_0}=0$ since $f(y^{*,j_0})>0$. Now assume that $H^k_{ij}=0$ for every $(i,j)\!\in I^*_0\times I^*_{>0}$. Then
\[
H^k_{i_0j_0} = \sum_{\ell=1}^d H_{i_0\ell}H^{k-1}_{\ell j_0} = \sum_{\ell\in I^*_{>0}}H_{i_0\ell} H^{k-1}_{\ell j_0} +\sum_{\ell\in I^*_0}H_{i_0\ell} H^{k-1}_{\ell j_0} = \sum_{\ell\in I^*_0}H_{i_0\ell} H^{k-1}_{\ell j_0} =0.
\]
This contradicts the irreducibility of $H$ since 
$i_0$ and $j_0$ are not connected through a power of $H$. 

\smallskip
\noindent $(e)$  We know that $\frac{1}{d}{\bf 1}\in{\cal E}_d$ and that, any $y^*\in{\cal E}_d$, $\min_iy^{*,i}>0$. Let $i_0$ be such that $y^{*,i_0}=\min_iy^{*,i}$. Then
$$
\sum_{j=1}^dH^{i_0j}f(y^{*,j})\geq\sum_{j=1}^dH^{i_0j}f(y^{*,i_0})=f(y^{*,i_0})
$$
since $H$ is stochastic. On the other hand, using the concavity of $f$, we derive that
$$
\forall \,y\!\in{\cal S}_d, \quad \Tr\big(\widetilde{f}(y)\big)=d\sum_{j=1}^d\frac{1}{d}f(y_j)\leq d\cdot f\Big(\frac{1}{d}\sum_{j=1}^d y_j\Big)=d\cdot f\Big(\frac{1}{d}\Big).
$$
As a consequence
$$y^{*,i_0}=\frac{\sum_{j=1}^dH^{i_0j}f(y^{*,j})}{\Tr(\widetilde{f}(y^*))}\geq
\frac{f(y^{*,i_0})}{d\cdot f\left(\frac{1}{d}\right)},
$$
which can be rewritten as $\frac{f(y^{*,i_0})}{y^{*,i_0}}\leq\frac{f(1/d)}{1/d}$. As $\xi\mapsto\frac{f(\xi)}{\xi}$ is (strictly) decreasing since $f$ is strictly concave and $f(0)=0$, it implies that $y^{*,i_0}\geq\frac{1}{d}$ which in turn implies  that $y^*=\frac{1}{d}{\bf 1}$ since $y^*\in{\cal S}_d$.~$\cqfd$
%
%
%
\paragraph{Remark.} When $H$ is not bi-stochastic but simply co-stochastic, we have no closed form for an $\overset{\circ}{{\cal S}_d}$-valued equilibrium and we could not manage to establish uniqueness even if $f$ is (strictly) concave.

\bigskip
Note that, as $f$ and $\Tr$ are defined on $[0,1]^d\setminus\{0\}$, $\widetilde f \big([0,1]^d\setminus\{0\}\big)\subset [0,1]^d\setminus\{0\}$ and $\Tr\big(\widetilde f(y)\big)>0$ on $[0,1]^d\setminus\{0\}$ since $f>0$ on $(0,1)$. Hence on may define the function 
\begin{equation}\label{eq:defvarphi}
\varphi:y\mapsto \frac{\widetilde f(y)}{\Tr(\widetilde f(y))}\quad \mbox{ on }\quad [0,1]^d\setminus\{0\}.
\end{equation}

\begin{lem}\label{lem:jacphi} 
If $f$ is differentiable on in the neighbourhood of $\frac 1d$,  then the functions $\widetilde{f}$, $\Tr\circ \widetilde f$ and $\varphi$   defined  on $[0,1]^d\setminus\{0\}$ are  bounded  on $(0,1]^d$ and the last two functions are differentiable  in the neighbourhood of $y(d)=\frac{1}{d}{\bf 1}$.   Moreover, $\varphi(y(d))= y(d)$ and the Jacobian of $\varphi$ at $y(d)$ is  given by 
$$
J_{\varphi}(y(d))=\begin{pmatrix} a & b & \ldots & \ldots & b \\ b & a & b & \ldots & b \\ \vdots & b & \ddots & \ddots & \vdots \\ \vdots & \vdots & \ddots & \ddots & b\\ b & b & \ldots & b & a\end{pmatrix} \; 
\mbox{ with } \; 
a=\frac{f'(1/d)(d-1)}{d^2f(1/d)} \quad \mbox{and} \quad b=-\frac{f'(1/d)}{d^2f(1/d)}=-\frac{a}{d-1}.
$$
As a symmetric matrix, $J_{\varphi}(y(d))$ is diagonalizable in the orthogonal group ${\cal O}(d,\R)$ with eigenvalues 0 (associated to the eigenvector ${\bf 1}$) and $\frac{f'(1/d)}{d\,f(1/d)}$ with eigenspace the hyperplane ${\bf 1}^{\bot}=\left\{u\in\R^d:\sum_{i=1}^du^i=0\right\}$.
\end{lem}

\noindent {\bf Proof.}   If $y\in(0,1)^d$ and $f$ is differentiable in the neighborhood of $y^1,\ldots,y^d$, then $\varphi$ is differentiable at $y$ and
\begin{equation}\label{eq:varphi}
\frac{\partial\varphi^i}{\partial y^j}(y)=-\frac{f(y^i)f'(y^j)}{\Tr\big(\widetilde{f}(y)\big)^2}+\delta_{ij}\frac{f'(y^i)}{\Tr\big(\widetilde{f}(y)\big)},\; 1\le i,j\le d.
\end{equation}
The form of $J_{\varphi}(y(d))$ follows.
Elementary and classical computations show that, for every real numbers $a$, $b$, 

$$
{\rm det}\begin{pmatrix}  a & b & \ldots & \ldots & b \\ b & a & b & \ldots & b \\ \vdots & b & \ddots & \ddots & \vdots \\ \vdots & \vdots & \ddots & \ddots & b\\ b & b & \ldots & b & a \end{pmatrix}=(a-b)^{d-1}(a+b(d-1)),
$$
which in turn implies that the characteristic polynomial of $J_{\varphi}(y(d))$ is given by $(a-b-\lambda)^{d-1}(a+b(d-1)-\lambda)$ so that the eigenvalues are
$$
\lambda_0=a+b(d-1)=0 \mbox{ (order 1) }\quad\mbox{and}\quad \lambda_1=a-b=\frac{f'(1/d)}{df(1/d)} \mbox{ (order $d-1$)}.
$$
The eigenspace associated to $\lambda_0$ is clearly $\R{\bf 1}$ and that associated to $\lambda_1$ is ${\bf 1}^{\bot}=\left\{u\in\R^d:\right.$ $\left.\sum_{i=1}^du^i=0\right\}$, so that $J_{\varphi}(y(d))_{\left|{\bf 1}^{\bot}\right.}=\lambda_1 {I}_{d\left|{\bf 1}^{\bot}\right.}$.   \hfill $\cqfd$

\bigskip
To go beyond and provide convergence results, especially in the bi-stochastic case, we will rely on an important tool in $SA$ theory, the $ODE$ method. This method makes the connection between  the asymptotic behaviour of the stochastic algorithm with mean field $h$ with that of  $ODE_h\equiv \dot\theta=-h(\theta)$ (for some key results on this theory, we refer to the Appendix). 

\begin{prop} $(a)$ From any $\xi\!\in {\cal S}_d$, there exists at least one ${\cal S}_d$-valued solution $(y(\xi,t))_{t\in \R_+}$ to $ODE_h$ starting from $\xi$.

\smallskip
\noindent $(b)$ From any $\xi\!\in \overset{\circ}{{\cal S}_d}= {\cal S}_d\cap (0,1]^d$,  there exists a unique solution  to $ODE_h$ starting from $\xi$, denoted $(\Phi(\xi,t))_{t\in \R_+}$ is $\overset{\circ}{{\cal S}_d}$-valued.
\end{prop}

\noindent {\bf Proof.} $(a)$ Let $\eta_0>0$. As $h$ is continuous ${\cal S}_d \cap[\eta_0,1]^d  \to {\bf 1}^{\bot}=\{\sum_{i=1}^d u ^i=0\}$, it is clear that, on the one hand, the Euler scheme of $ODE_h$ with step $\frac 1n$ starting from 
$\xi\!\in {\cal S}_d$ defined by $\bar y^n_{\frac{k+1}{n}}=\bar y_{\frac kn}^n-\frac 1n h(\bar y^n_{\frac kn})$ is ${\cal S}_d$-valued since $\Tr\big(h(\bar y^n_{\frac kn})\big)$ is constant equal to $\Tr(y)=1$, and, on the other hand, converges owing to Peano's theorem to a solution of $ODE_h$.

\smallskip
\noindent $(b)$ The function  $u\mapsto f(u)$ is Lipschitz continuous on $[\eta_0, 1]$, with Lipschitz coefficient $f_\ell'(\eta_0)\le 1$ since $f'$ is concave. Hence the (canonical extensions of) the  functions $\widetilde f$ and $\mathbf{w}(\widetilde f)$   on $[\eta_0, 1]^d$ are Lipschitz as well. Both are also bounded. Moreover $\Tr\big(\widetilde{f}(y)\big)\geq f\left(\frac{1}{d}\right)>0$ because for any  $y=(y^1,\ldots,y^d)^t\!\in{\cal S}_d$, there exists $i_0\in\{1,\ldots,d\}$ such that $y^{i_0}\geq\frac{1}{d}$ so that $\Tr\big(\widetilde{f}(y)\big)\geq f(y^{i_0})\geq f\left(\frac{1}{d}\right)>0$. As a consequence, $\mathbf{w}(\widetilde f) \ge \frac 12f\big(\frac 1d\big)$ on a neigbourhood of ${\cal S}_d$ in $[0,1]^d$. Therefore, the functions $\varphi_{_H}:y\mapsto H \frac{\widetilde{f}(y)}{\Tr\big(\widetilde{f}(y)\big)}$   from  ${\cal S}_d\cap [\eta_0,1]^d $ to ${\cal S}_d$ and $h= I_d-\varphi_{_H}$ from ${\cal S}_d\cap [\eta_0,1]^d  $ to $\{\sum_{i=1}^d u ^i=0\}$ are Lipschitz continuous too in that neighbourhood. This guarantees the existence of an ${\cal S}_d$-valued flow for $ODE_h$  on ${\cal S}_d\cap (0,1]^d= \overset{\circ}{{\cal S}_d}$. 

\begin{Definition}[Stable equilibrium point]\label{def:stabeq} An element $y^*$ of ${\cal S}_d$ is a {\em  stable equilibrium} for  $ODE_h\equiv \dot{y}=-h(y)$ on ${\cal S}_d$ if  there is  (compact) neighborhood ${\cal K}^*$ of $y^*$ in ${\cal S}_d$ such  that
\[
\lim_{t\to +\infty}\sup\Big\{\big|y(\xi,t)-y^*\big|, \, \xi\!\in {\cal K}^*,\; y(\xi,\cdot)\; {\cal S}_d\mbox{-valued solution to $ODE_h$}, \, y(\xi,0)=0\Big\} =0
\]
with a slight abuse of notation since possibly several solutions may start from $\xi \!\in \partial {\cal S}_d$. See also Theorem~\ref{ThmODE}$(b)$.
\end{Definition}

\paragraph{Remarks.}  $\bullet$  If the ${\cal S}_d$ -valued  flow of  $ODE_h$ is well-defined in the neighbourhhood of $y^* $ in ${\cal S}_d$, this  boils down to showing that  $y(\xi,t)$ is ${\cal S}_d$-valued and converges to $y^*$ as $t\to+\infty$, uniformly with respect to $\xi$   in  a (compact) neighbourhood of $y^*$ in ${\cal S}_d$.

\smallskip
\noindent $\bullet$ Since $h$ is differentiable, an equilibrium $y^*$  is attractive if all the eigenvalues of $J_{h}(y^*)_{|{\bf 1}^{\bot}} = \big(I_{d}-H J_{\varphi}(y^*)\big)_{|{\bf 1}^{\bot}}$  have (strictly) positive real parts, see~\cite{Ben}. If one of these eigenvalues has a negative real part then the equilibrium is {\em unstable} and if all eigenvalues have negative real parts the equilibrium is called a {\em repeller}. 

\begin{prop} \label{pro:stabinstab} Assume $H$ be a bi-stochastic matrix (keep in mind that  $y(d)$ is  a zero of $h$).
Let $ \lambda_1=\frac{f'(1/d)}{d\cdot f(1/d)} $ and let  $\mu_{\max}$ be the  eigenvalue of $H_{|{\bf 1}^{\bot}}$ with the highest real part.
 
\smallskip
\noindent $(a)$   If  $\Re e(\mu_{\max})<\frac{1}{\lambda_1}$, then $y(d)=\mbox{\bf 1}/d$ is always a stable equilibrium    of $ODE_h\equiv\dot{y}=-h(y)$. Thus, if 
\[
\left(\lambda_1 \le 1\,\mbox{ and  $1\!\notin {\rm Sp}(H_{|{\bf 1}^{\bot}})$}\right) \, \mbox{or}\, \left(\lambda_1<1 \right),
 \]
then $y(d)$ is always a stable equilibrium. 

\smallskip In particular the  above left  condition is satisfied if $H$ is irreducible and $f$ is concave, whereas the right one is always fulfilled if $f$ is strictly concave over $(0, 1/d)$. 
 
\smallskip
\noindent $(b)$ If $\Re e(\mu_{\max})>\frac{1}{\lambda_1}$, then $y(d)$ is unstable (it is even a repeller when $H=I_d$).  Note that if $f$ is convex (resp. strictly over $(0, \frac 1d)$), then $\lambda_1\ge 1$ (resp. $>1$).
\end{prop}
%
%
%
\noindent {\bf Remark.}
%
\noindent $\bullet$ If $\lambda_1>1$ ($e.g.$ since $f$ is strictly convex), then $\frac{1}{\lambda_1}<1$ and the two opposite situations $\Re e(\mu_{\max})<\frac{1}{\lambda_1}$  and $\Re e(\mu_{\max})>\frac{1}{\lambda_1}$ may occur a priori {\it i.e.} $y(d)$ may switch from uniform stability to instability (see Section~\ref{trois} for the two-type case: $d=2$).

\smallskip
\noindent $\bullet$ Note that if $f$ is strictly convex and $1\!\in {\rm Sp}(H_{|{\bf 1}^{\bot}})$, then $y(d)$ is always unstable.

\bigskip
\noindent {\bf Proof.} $(a)$ It follows from Lemma~\ref{lem:jacphi} that $J_{\varphi}(y(d))= \lambda_1I_{d |{\bf 1}^{\bot}}$. Hence
\begin{eqnarray*}
	J_h(y(d))_{\left|{\bf 1}^{\bot}\right.}&=&I_{d|{\bf 1}^{\bot}}-\big(HJ_{\varphi}(y(d))\big)_{ |{\bf 1}^{\bot} }={I}_{d |{\bf 1}^{\bot}}-H\lambda_1{I_d}_{\left|{\bf 1}^{\bot}\right.}=(I_d-\lambda_1H)_{|{\bf 1}^{\bot} }.
\end{eqnarray*}
Hence  $\mbox{Sp}(J_h(y(d))_{|{\bf 1}^{\bot}})=\big\{1-\lambda_1\mu, \mu\in\mbox{Sp}(H_{|{\bf 1}^{\bot}})\big\}\subset \mathbb{C}$.

Every $\mu\!\in {\rm Sp}(H)$ satisfies $|\mu|\leq1$ since $H$ is stochastic.  If  $\lambda_1<1$, then $|\lambda_1\mu_{\max}|<1$ and consequently $\Re e(1-\lambda_1\mu_{\max})>0$ which ensures that $y(d)$ is attracting. The other  case follows likewise. 

When $H$ is irreducible, the eigenvalue $1$ is simple owing to the Perron-Frobenius Theorem so that $1\!\notin {\rm Sp}(H_{|{\bf 1}^{\bot}})$.

\smallskip
\noindent $(b)$ If $f$ is convex (and $f>0$ on $(0,1)$), then $0< f(1/d) \le f'(1/d)/d$ since $f(0)=0$ so that $\lambda_1\ge 1$. This inequality is strict if  $f$ is strictly convex over $(0, \frac 1d)$.~$\qquad\cqfd$ 

\begin{prop} \label{prop:flot}$(a)$ If $H$ is bi-stochastic and irreducible and $f$ is  strictly concave, then $ODE_h$ has at least one solution starting from every $\xi\!\in{\cal S}_d$. The flow of $ODE_h\equiv \dot{y}=-h(y)$, denoted by $\big(y(\xi,t)\big)_{\xi\in \overset{\circ}{{\cal S}_d}}$, $t\geq 0$, exists on  $\overset{\circ}{{\cal S}_d}$ and is ${\cal S}_d$-valued. Furthermore (still with $y(d)=\frac{{\bf 1}}{d}$),  
\[
\lim_{t\to +\infty}\sup\Big\{\Big|y(\xi,t)-y(d)\Big|, \, \xi\!\in {\cal S}_d,\; y(\xi,.)\; \mbox{solution to $ODE_h$}, \, y(\xi,0)=0\Big\} =0.
\]

\smallskip
\noindent $(b)$ If  $H$ is simply bi-stochastic (and possibly not irreducible), then the  flow, denoted by $y(\xi,t)$, converges toward $y(d)$  for every  $\xi\!\in \overset{\circ}{{\cal S}_d}$.
\end{prop}

\noindent {\bf Proof.} $(a)$  Let us denote by $y(\xi,t)$ an  ${\cal S}_d$-valued solution starting from $\xi \!\in {\cal S}_d$.

Let $\xi\in{\cal S}_d\!\setminus\!\{y(d)\}$
and let $i(t)$ be the right continuous function  such that  $y^{i(t)}(\xi,t)=\min_jy^j(\xi,t)\!\in [0,\frac 1d]$. We know that $y^{i(0)}(\xi,0)=\xi^{i(0)}<\frac{1}{d}$ since $\xi\neq y(d)$. Note that the function $g_i:t\mapsto y^{i(t)}(\xi,t)$ is continuous and right  differentiable with a   right derivative given by 
$$
(g_{i})_r'(\xi,t)=\sum_{j=1}^dH^{i(t)j}\frac{f(y^j(\xi,t))}{\Tr(\widetilde{f}(y(\xi,t)))}-y^{i(t)}(\xi,t)\geq 1\times\frac{f(y^{i(t)}(\xi,t))}{\Tr(\widetilde{f}(y(\xi,t)))}-y^{i(t)}(\xi,t).
$$
Note that $\Tr(\widetilde{f}(y(\xi,t)))\leq d\cdot f\left(\frac{1}{d}\right)$ since $f$ is concave. It follows that
\begin{equation}\label{star}
(g_{i})_r'(\xi,t)\geq \frac{f(y^{i(t)}(\xi,t))}{d \cdot f\left(\frac{1}{d}\right)}-y^{i(t)}(\xi,t)\ge 0,
\end{equation}
since $u\mapsto\frac{f(u)}{u}$ is non-increasing  ($f$ is  concave and    $f(0)=0$) and $y^{i(t)}(\xi,t)\le \frac 1d$ for every $t\ge 0$.

\smallskip If $\xi^{i(0)}=0$, then $y^{i(t)}(\xi,t)\ge 0$, for every $t\ge 0$. Assume there exists $\varepsilon_0>0$ such that $y^{i(t)}(\xi,t)= 0$ for every $t\!\in (0, \varepsilon_0]$, then one derives from the integrated form of $ODE_h$ that  
\[
\int_0^{\varepsilon_0}\left(\sum_{j=1}^dH^{i(t)j}\frac{f(y^j(\xi,s))}{\Tr(\widetilde{f}(y(\xi,s)))}\right)ds = 0,
\]
so that, as the above integrand is non-negative and continuous, $s\mapsto \sum_{j=1}^dH^{i(s)j}\frac{f(y^j(\xi,s))}{\Tr(\widetilde{f}(y(\xi,s)))}\equiv 0$ on $[0, \varepsilon_0]$. Let $I_0=\{i\, \mbox{ s.t. }\, y^{i}(\xi,\varepsilon_0)=0,\, 1\le i\le d\}$ and $I_1=I_0^c$. Then, it follows that $i(\varepsilon_0)\!\in I_0$ and $H^{i(\varepsilon_0)j}=0$, for every $j\!\in I_1$, which is not empty since $\sum_j y^j(\xi,t)=1$. Hence, reasoning like in Proposition~\ref{prop:fixpt}$(d)$,  $i(\varepsilon_0)$ (and more generally any element of $I_0$) and $I_1$ are not $H$-connected which contradicts the irreducibility. Consequently, $y^{i(t)}(\xi,t)>0$, $t\!\in(0,\varepsilon_0]$, $i.e.$  $y(\xi,t)$ lives in $\overset{\circ}{{\cal S}_d}$ and, consequently,  lives in $\overset{\circ}{{\cal S}_d}$ for every $t\ge 0$ by monotonicity of $y^{i(t)}(\xi,t)$.

\smallskip If $\xi^{i(0)}>0$, the same conclusion follows simply from the monotonicity   of $y^{i(t)}(\xi,t)$. 

Finally, for every   $t>0$, $y^{i(t)}(\xi,t)$ is positive and increasing. We can integrate~\eqref{star} between a fixed $\varepsilon>0$ and $t>\varepsilon$ as follows
%
$$-\log \left(\frac{d}{\xi^{i(0)}(\varepsilon)}\right)\ge \log\left(\frac{y^{i(t)}(\xi,t)}{\xi^{i(0)}(\varepsilon)}\right)\ge \int_\varepsilon^t\left(\underbrace{\frac{f(y^{i(s)}(\xi,s))}{y^{i(s)}(\xi,s)}-\frac{f\left(1/d\right)}{1/d}}_{\ge 0}\right)ds \ge 0
$$
and the latter integral  is increasing in $t$ as long  as $y^{i(t)}(\xi,t)<\frac 1d$. As the left hand side of the above string of inequalities is finite, it follows that $\int_0^{+\infty}\left(\frac{f(y^{i(s)}(\xi,s))}{y^{i(s)}(\xi,s)}-\frac{f\left(1/d\right)}{1/d}\right)ds<+\infty$. Then, combining that $u\mapsto\frac{f(u)}{u}$ is decreasing (by strict concavity) on $(0,+\infty)$ with the fact that $y^{i(t)}(\xi,t)$ is increasing and positive, we derive that  the function $t\mapsto\frac{f(y^{i(t)}(\xi,t))}{y^{i(t)}(\xi,s)}- \frac{f\left(1/d\right)}{1/d}$ is decreasing. The finiteness of the above integral implies that $ \displaystyle \frac{f\left(1/d\right)}{1/d}-\frac{f(y^{i(t)}(\xi,t))}{y^{i(t)}(\xi,t)}\to 0$ as $t\to+\infty$ or,  equivalently, that
$$
y^{i(t)}(\xi,t)=\min_{1\le i\le d} y^i(\xi,t)\to\frac1d\quad\mbox{ as }\quad t\to+\infty.
$$  
The convergence of every component $y^i(\xi,t)$ follows since $y(\xi,t)\!\in {\cal S}_d$. It remains to prove that the convergence holds uniformly in the starting value (the existence of the flow $i.e$ uniqueness starting from the boundary of the simplex is not mandatory).  First note that, $h$ being bounded, the family of all possible solutions $((y(\xi,t)_{t\ge 0})_{\xi \!\in {\cal S}_d})$  of $ODE_h$ are $\|h\|_{\infty}$-Lipschitz continuous  since
\[
y(\xi,t)= \xi -\int_0^t h\big(y(\xi,s)\big)ds.
\]
Assume there exists $\xi_n \to \xi_{\infty}$ and $t_n\to +\infty$ such that $\big|y(\xi_n,t_n)-y(d)\big| \ge \varepsilon_0\ge 0$. By Proposition~\ref{lem:jacphi}, $\frac{\mbox{\bf 1}}{d}$  is an attractor of $ODE_h$,  hence there exists $\eta_0>0$ such that  $\sup_{|\xi-y(d)|\le  \eta_0} \big|y(\xi,t)-y(d)\big|\le \varepsilon_0$, $t\ge t_0$ (the flow does exist in the neighbourhood of $y(d)$ so that $y(\xi,\cdot)$ is unique). Consequently, for every $t\!\in [0,t_n-t_0]$ (at least for $n$ large enough), $|y(\xi_n,t)-y(d)|> \eta_0$. By Ascoli's Theorem, up to an extraction, one may assume that $y(\xi_n,\cdot)$ converges on compact sets of $\R_+$ toward a solution $y(\xi_\infty,\cdot)$ of $ODE_h$ since $h$ is continuous. Then, letting $n$ go to infinity shows that this  solution satisfies $|y(\xi_\infty,t)- y(d)| \ge \eta_0$ for every $t\ge 0$. This contradicts the fact that any solution starting from ${\cal S}_d$ converges toward $y(d)$.  We can conclude that all the solutions of $ODE_h$ starting from ${\cal S}_d$ converge uniformly toward $y(d)$.

\smallskip
\noindent  $(b)$ is obvious given the above proof. \hfill $\cqfd$

\bigskip
\noindent {\bf Application~(I):  \texorpdfstring{$A.s.$}{A.s.} convergence of the algorithm when $f$ is concave.}
When $f$ is   concave and $H$ is bi-stochastic, the mean point $y(d)= \frac{\mbox{\bf 1}}{d}$  of the simplex is the unique $a.s.$ target of the urn composition. 
\begin{prop}[When $f$ is concave $y(d)$ is the  target]\label{prop:convy(d)}
 If $H$ is (deterministic and) bi-stochastic, irreducible and $f$ is a strictly concave skewing function, then 
\[
\widetilde Y_n \overset{a.s.}{\underset{n\to+\infty}{\longrightarrow}} y(d)= \frac{\mbox{\bf 1}}{d}.
\]
\end{prop}

\noindent {\bf Proof.} It follows from Proposition~\ref{prop:fixpt}$(e)$ that the flow of $ODE_h$  uniformly converges to   $y(d)$. Hence, 
 Theorem~\ref{ThmODE}$(b)$ from the  Appendix  implies that  the set $\Theta^{\infty}$ of the limiting values of $ODE_h$ is reduced to $\{\frac 1d \mbox{\bf 1}\}$ which completes the proof. \hfill$\cqfd$

\paragraph{Application~(II): The mean point $y(d)$ may be a noisy trap when $f$ is convex.}
We now give a (partial) result in the convex setting (a more precise one is provided in  Section~\ref{quatre} devoted to randomized P\'olya's urns, that is when $H=I_d$).  We keep  the  notations introduced Propositions~\ref{prop:minmax} to~\ref{pro:stabinstab}: $\lambda_1= \frac{f'(1/d)}{df(1/d)}$ and $\mu_{\max}$ the eigenvalues of $H_{|\mbox{\bf 1}^{\perp}}$ with the highest real part. 

As $H$ is (bi-)stochastic $\Re e(\mu_{\max})\le 1$. If $H$ is also irreducible, by the Perron-Frobenius Theorem, we know that $\Re e(\mu_{\max}) <1$ since  the eigenspace of the  eigenvalue $1$ is one-dimensional, hence equal to $\R\mbox{\bf 1}$.  On the other hand, when $f$ is convex then $\lambda_1\ge 1$ and even $>1$ if $f$ is {\it e.g.} strictly convex over $(0, 1/d)$. Consequently there exist matrices $H$ for which the  inequality $\Re e(\mu_{\max})>\frac{1}{\lambda_1}$ is satisfied, implying by Proposition~\ref{pro:stabinstab} that $y(d)$  is unstable for $ODE_h$. We show that $y(d)$ may be a noisy trap. 

\begin{prop}[When  $y(d)$ becomes a trap] Assume that $f$ is convex and $H$ is  a bi-stochastic matrix satisfying $\Re e(\mu_{\max})>\frac{1}{\lambda_1}$ so that  $y(d)$ is unstable. Assume  that $\ds \sup_{n\geq1} \E\left[|\!|\!|D_{n+1}|\!|\!|\left|\right.{\cal F}_{n-1}\right]\le L\!\in \R_+$.

\smallskip
\noindent  $(a)$  If, for every  $v\!\in \R^d$ with $\Tr(v)=0$ and $\|v\|=1$, 
\[
 \liminf_n \E\big(\|D_{n+1}v\|^2\,|\, {\cal F}_{n}\big)>0,
\]
then $\displaystyle \P\big(\widetilde Y_n \to y(d)\big) =0$.

\smallskip 
\noindent $(b)$ Assume that $H$ is also symmetric. Let $v_{\mu}$ be a unitary eigenvector attached to an eigenvalue $\mu$ of $H_{|\mbox{\bf 1}^{\perp}}$ such that $\Re e(\mu)>\frac{1}{\lambda_1}$. If 
\[
 \liminf_n \E\left[\|D_{n+1}v_{\mu}\|^2\,|\, {\cal F}_{n}\right]>0,
\]
then $\displaystyle \P\big(\widetilde Y_n \to y(d)\big) =0$.
\end{prop}

\noindent {\bf Proof.} We will focus on claim $(b)$ (claim $(a)$ relies on the same formal proof applied to any nonzero vector $v$ with $0$ weight). 
We want to apply Theorem~\ref{ThmPiege} from the  Appendix.  The first point to be checked is that the Jacobian $J_h$ is locally Lipschitz continuous in the neighbourhood of $y(d)$. Note that $h= I_d-H\varphi$ where $\varphi $ is defined in~\eqref{eq:defvarphi}.  This follows from the expression~\eqref{eq:varphi} for $\frac{\partial \varphi_i}{\partial x_j}$  having in mind that $f'$ is non-decreasing by convexity of $f$. Then, we have to check Assumption~\eqref{BruitMgle}. Elementary computations show that 
\[
\E\left[(\Delta M_{n+1}|v_{\mu})^2\,|\, {\cal F}_{n}\right]= v_{\mu}^t\left[ \E\left[D_{n+1}{\rm diag}(\varphi_{_H}(\widetilde Y _n)) D_{n+1}^t\,|\, {\cal F}_{n}\right]-2\varphi_{_H}(\widetilde Y_n)\otimes \widetilde Y_n +\widetilde Y_n ^{\otimes 2}\right]v_{\mu}.
\]
On the event $\big\{\widetilde Y_n \to y(d)\big\}$, we derive, owing to Assumption  ({\bf A2}),  the continuity of $\varphi_{_H}$ and $\varphi_{_H}(y(d))=y(d)$, that  
\[
\E\left[(\Delta M_{n+1}|v_{\mu})^2\,|\, {\cal F}_{n}\right]= v_{\mu}^t \E\left[D_{n+1}{\rm diag}\big(y(d)\big)D_{n+1}^t\,|\, {\cal F}_{n}\right]v_{\mu}- (v_{\mu}|y(d))^2 + o(1).
\]
Now,  ${\rm diag}(y(d))=\frac 1d I_d$, and $v_{\mu}\!\in \mathbf{1}^{\perp}$ since  $(v_{\mu}|y(d))^2=0$. Consequently,  
\[
\liminf_n \E\left[(\Delta M_{n+1}|v_{\mu})^2\,|\, {\cal F}_{n}\right]= \frac 1d \liminf_n \E\left[\|D_{n+1}v_{\mu}\|^2\,|\, {\cal F}_n\right] >0. 
\]
Now, we  note that
\[
\|D_{n+1}v_{\mu}\|\ge \frac{\|D_{n+1}v_{\mu}\|^2}{|\!|\!|D_{n+1}|\!|\!|}\ge \frac{1}{L},
\]
so that $\displaystyle \liminf_n \E\big[|(\Delta M_{n+1}|v_{\mu})|\,|\, {\cal F}_{n}\big]>0$. Now, as  $H$ is  symmetric and leaves $\mbox{\bf 1}^{\perp}$ stable,  $H$ is self-adjoint on ${\bf 1}^{\perp}$. Hence, there exists an orthonormal basis of ${\bf 1}^{\perp}$ containing $v_{\mu}$ so that, 
\[
 \liminf_n \E\left[\|\Delta M_{n+1}\|^2\,|\, {\cal F}_{n}\right]\ge \liminf_n \E\left[(\Delta M_{n+1}|v_{\mu})^2\,|\, {\cal F}_{n}\right]>0.\qquad \cqfd
 \]

%

\section{Bi-dimensional non linear randomized urn model}\label{trois}

When the skewing function $f$ is convex the situation becomes much more involved: Thus, if $H$ is bi-stochastic and irreducible, then $y(d)$ is  an equilibrium of $h$ and ${\cal E}_d\subset \overset{\circ}{{\cal S}_d}$ (cf. Proposition~\ref{prop:fixpt}$(d)$), but ${\cal E}_d$ is not reduced to $y(d)$ (which is a unstable for $ODE_h$ when $H=I_d$). 

To start elucidating this case, we limit ourselves in this paper to a two-type urn ($d=2)$ and an irreducible matrix $H$. We will see, as expected, that the asymptotic behavior of the urn is much more  involved since a phase transition appears.

In such a simplified setting, the irreducible  co-stochastic generating matrix $H$ can be written as follows
$$
H=\begin{pmatrix}p_1 & 1-p_2 \cr 1-p_1 & p_2 \end{pmatrix}, \qquad 0<p_i<1, \quad i=1,2,
$$
and the  mean function $h$ of  the model is still given by~(\ref{defhcvx}).

In this section, we assume for  simplicity that {\em $f$ is differentiable on $(0,1]$}, so that the mean field $h$ is differentiable too on $\overset{\circ}{{\cal S}_2}$. We saw that analyzing the $a.s.$ convergence properties essentially boils down to elucidate the behavior of $ODE_h$ on ${\cal S}_2$ which in turn can be reduced to a one dimension differential system since the simplex ${\cal S}_2$ can be parametrized by $(u,1-u)$, $u\!\in [0,1]$, and $h$ is  $\mbox{\bf 1}^{\perp}=\{z:z^1+z^2=0\}$-valued on $\overset{\circ}{{\cal S}_2}$ so that  $h^2=1-h^1$. Thus, the asymptotic analysis of $ODE_h$ is equivalent to that of  
\[
ODE_{h_0} \equiv \dot{u} = -h_0(u)\quad \mbox{ with }\quad h_0(u)= h^1(u,1-u),\quad u\!\in [0,1].
\]
Elementary computations show that, for every $u\!\in [0,1]$, 
\begin{equation}\label{eq:h0}
h_0(u)=u-\frac{p_1f(u)+(1-p_2)f(1-u)}{f(u)+f(1-u)}
\end{equation}
and
\begin{equation}\label{eq:h'0}
h'_0(u)=1-(p_1+p_2-1)\frac{f'(u)f(1-u)+f(u)f'(1-u)}{\left(f(u)+f(1-u)\right)^2}.
\end{equation}
We will now determine the combinatorics and the nature of the equilibrium points depending on the parameters $p_1$ and $p_2$.

Note that $H$ is bi-stochastic if and only if $p_1= p_2$. Note also that $h_0\big(\frac12\big)=\displaystyle\frac{p_2-p_1}{2}$ whatever the skewing function $f$ is.

\subsection{Counting equilibrium points}
%

It follows from~\eqref{eq:h0} that the equation $h_0(u)=0$ reads
\begin{equation}\label{EConvex}
	(p_1-u)f(u)+(1-p_2-u)f(1-u)=0.
\end{equation}
As preliminary  remarks, note  that:

\smallskip
--  if $p_1=1-p_2$, then $u^*=p_1$ is the unique solution of~\eqref{EConvex} because $f>0$ on $(0,1]$ and $h'_0(u^*)= 1>0$ so that $u^*$ is a stable equilibrium. 

\smallskip
-- $h'_0$ being symmetric w.r.t. $\frac12$, $h''_0$ is antisymmetric w.r.t. $\frac12$. Hence $h''_0\left(\frac12\right)=0$.

\smallskip
-- When the generating matrix $H$ is bi-stochastic (if and only if $p_1=p_2$), then $u^*= \frac 12$ is always solution to~\eqref{EConvex}. Thus, if $h'_0\left(\frac12\right)=0$, as $h''_0\left(\frac12\right)=0$, the status of the equilibrium $u^*$ requires  to investigate higher order (see Example  in Section~\ref{subsec:StabODEh2d} devoted to the bi-stochastic case).

\bigskip
Let $u^*$ be a solution of~(\ref{EConvex}). From Proposition~\ref{prop:minmax}$(b)$, we have that any zero $u^*$ of $h_0$ lies in $I^*:=[p_1\wedge(1-p_2),p_1\vee(1-p_2)]$ (or in  its interior $\overset{\circ}{I^*}$).  

\begin{prop}\label{PropEq} Let $p_1$, $p_2\!\in (0,1)$ and $f$ a skew-function.

\smallskip
\noindent $(a)$ If $p_1+p_2\leq1$, then~(\ref{EConvex}) has a unique solution $u^*$ lying in $I^*$ (and in $\overset{\circ}{I^*}$ if $p_1+p_2<1$). 

\smallskip
\noindent $(b)$  If $p_1+p_2>1$ and $f$ is concave, then~(\ref{EConvex}) has a unique solution $u^*$ lying in~$\overset{\circ}{I^*}$. 

\smallskip
\noindent $(c)$ If $p_1+p_2>1$ and $f$ is strictly convex, then~(\ref{EConvex}) may have one, two or three solutions lying in $\overset{\circ}{I^*}$. However, when  $f'(1)\le  \frac{1}{p_1+p_2-1}$,  then $h_0$ is increasing   and~(\ref{EConvex}) subsequently has a unique zero $u^*$  lying in $\overset{\circ}{I^*}$. 

\end{prop}

\noindent{\bf Proof.} $(a)$ 
As $f$ is increasing and non-negative, $h_0'>0$ when $p_1+p_2\leq 1$, therefore $h_0$ is increasing with $h_0(0)=p_2-1<0$ and $h_0(1)=1-p_1>0$, so $h_0$ has a unique zero lying in $I^*$.

\smallskip
\noindent $(b)$ Assume that $p_1+p_2>1$. Then it is obvious that $p_1$ and $1-p_2$ are not solutions of~(\ref{EConvex}) which can be rewritten as
$$
g_1(u)=g_2(u)  \quad\mbox{on}\quad {\cal J}=[0,1]\setminus\{p_1,1-p_2\},
$$
where
\begin{equation}\label{defg1g2}
g_1(u)=\frac{f(1-u)}{p_1-u}, \quad u\!\in[0,1]\setminus\{p_1\} \quad\mbox{and}\quad g_2(u)=\frac{f(u)}{u-1+p_2}, \quad u\!\in[0,1]\setminus\{1-p_2\}.
\end{equation}
Let us compute the first derivative of these functions: We obtain, for $u\!\in {\cal J}$,
$$
g'_1(u)=\frac{f(1-u)-f'(1-u)(p_1-u)}{(p_1-u)^2} \quad\mbox{and}\quad g'_2(u)=\frac{f'(u)(u-1+p_2)-f(u)}{(u-1+p_2)^2}.
$$
As $f$ is concave non-negative and  $f(0)=0$, $f(u)-uf'(u)\geq0$, $u\!\in[0,1]$. Let us show that, as  $0<p_1,p_2<1$, $g'_1>0$ and $g'_2<0$ on $I^*$. In fact, 
\[
g'_2(u)= \frac{u f'(u)-f(u)-(1-p_2)f'(u)}{(u-1+p_2)^2}\le 0\quad \mbox{since}\quad f'(u)\ge 0,
\]
and if $u f'(u)-f(u)-(1-p_2)f'(u)=0$, then  $(1-p_2)f'(u)\le 0$. Hence $f'(u)=0$ since $p_2<1$ so that $f=0$ on $[0,u]$ which contradicts the positivity of $f$ on $(0,1)$ as a skew-function. This shows that $g'_2<0$. One shows likewise that $g'_1>0$. Hence~(\ref{EConvex}) has a unique solution $u^*\in \overset{\circ}{I^*}$. 

\smallskip
\noindent $(c)$ If $p_1+p_2>1$ and $f$ is strictly convex, we set
$$
\varphi_1(u)=(p_1-u)f(u) \quad\mbox{and}\quad \varphi_2(u)=(u-1+p_2)f(1-u).
$$
 Then we have three possibilities: When $f$ is strictly convex, the functions $\varphi_1$ and $\varphi_2$ have a maximum on $I^*$. Let us set $u_1$, $u_2\!\in I^*$ such that $\varphi'_1(u_1)=0$ and $\varphi'_2(u_2)=0$. As $p_1+p_2-1>0$ and $f$ is strictly convex, $u_2<u_1$ and, depending on the relative position of $\varphi_1(u_2)$ and $\varphi_2(u_2)$, we obtain the three announced results. \hfill$\cqfd$

\clearpage
\noindent {\bf Example of phase transition~(I) : Convex power skewing function.} Let us illustrate the strictly convex setting  by considering the family of strictly convex functions  $f(u)= f_{\alpha}(u)=u^{\alpha}$, $\alpha>1$. 

Then Equation $h_0(u)=0$ (see~(\ref{EConvex})) can be rewritten as
\[
\varphi_1(u)=\varphi_2(u) \quad\mbox{on}\quad I^*,
\]
where
\[
\varphi_1(u)=(p_1-u)f(u)=(p_1-u)u^{\alpha}\quad\mbox{and}\quad \varphi_2(u)=(u-1+p_2)f(1-u)=(u-1+p_2)(1-u)^{\alpha}.
\]
Then, the first derivatives of these functions read
\[
\varphi'_1(u)=\big(\alpha p_1-(\alpha+1)u\big)u^{\alpha-1}\quad\mbox{and}\quad \varphi'_2(u)=\big((\alpha+1)(1-u)-\alpha p_2\big)(1-u)^{\alpha-1}.
\]
Each of these two derivatives $\varphi'_i$ has only one zero lying in $I^*$ which are the maximum of each function $\varphi_i$, $i=1,2$: $u_1=\frac{\alpha p_1}{\alpha+1}$ for $\varphi'_1$ and $u_2=1-\frac{\alpha p_2}{\alpha +1}$ for $\varphi'_2$. As $p_1+p_2-1>0$, we have that $1-\frac{\alpha p_2}{\alpha +1}<\frac{\alpha p_1}{\alpha+1}$. Besides, we have $\varphi_1(1-p_2)=(p_1+p_2-1)(1-p_2)^{\alpha}>0$, $\varphi_1(p_1)=0$, $\varphi_2(1-p_2)=0$ and $\varphi_2(p_1)=(p_1+p_2-1)(1-p_1)^{\alpha}>0$. Therefore we have three possibilities (of course consistent with claim~$(c)$ of the above proposition):
\begin{itemize}
	\item[$\bullet$] If $\varphi_1\left(1-\frac{\alpha p_2}{\alpha +1}\right)>\varphi_2\left(1-\frac{\alpha p_2}{\alpha +1}\right)$, then~(\ref{EConvex}) has a unique solution lying in $I^*$.
	\item[$\bullet$] If $\varphi_1\left(1-\frac{\alpha p_2}{\alpha +1}\right)=\varphi_2\left(1-\frac{\alpha p_2}{\alpha +1}\right)$, then~(\ref{EConvex}) has two solutions lying in $I^*$.
	\item[$\bullet$] If $\varphi_1\left(1-\frac{\alpha p_2}{\alpha +1}\right)<\varphi_2\left(1-\frac{\alpha p_2}{\alpha +1}\right)$, then~(\ref{EConvex}) has three solutions lying in $I^*$.
\end{itemize}

\smallskip
In this parametric setting,  as $f'(1)= \alpha$, it follows from Proposition~\ref{PropEq} that the simpler  criterion for the existence of a unique attractive equilibrium reads 
\[
1<p_1+p_2 \le 1+\frac{1}{\alpha}.
\]
To exhibit a transition  phase with three equilibrium points, this condition need to  be violated. However it is not a sufficient condition. This is illustrated by Figure~\ref{fig:phasetransition} below: Choose $p_1$ and $p_2$ such that $p_1+p_2>1$. The critical $\alpha_0$ is clearly strictly larger than the barrier $1/(p_1+p_2-1)$ suggested by the criterion. For example, for $p_1=0.7$ and $p_2=0.75$, $\alpha_0\approx3.09>1/(p_1+p_2-1)=\frac{20}{9}$.

\begin{figure}[!ht]\label{fig:phasetransition}
\centering
\includegraphics[width=10cm]{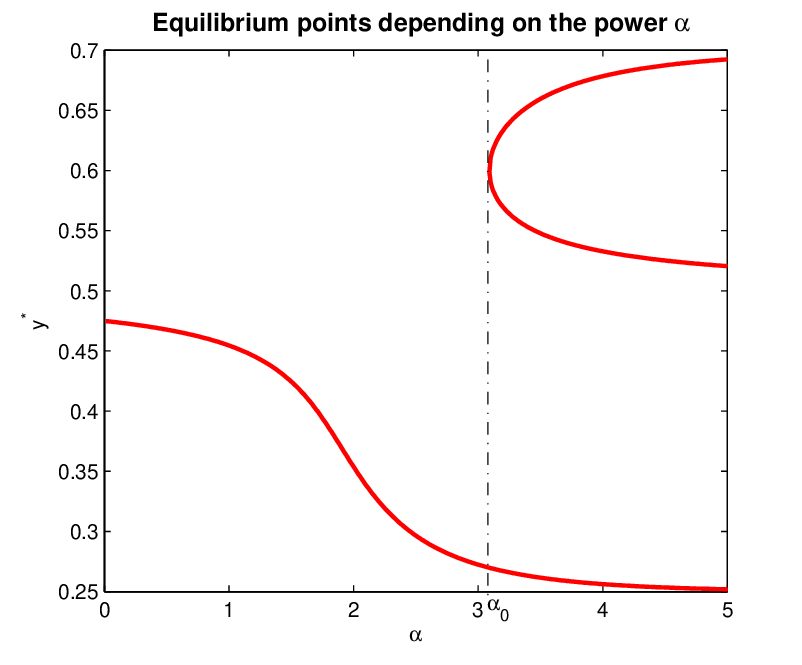}
\vskip-0.7cm
\caption{\em Equilibrium points for $f(u)=u^{\alpha}$ depending on $\alpha\!\in (0,5]$ with $p_1=0.7$ and $p_2=0.75$. Bifurcation appears at $\alpha_0\approx3.09$.}
\end{figure}

\subsection{Stability of equilibrium points of \texorpdfstring{$ODE_{h_0}$}{ODE h0}}\label{subsec:StabODEh2d}

Once  elucidated  the  existence and the  number of  equilibrium points  $u^*$ of $ODE_{h_0}$, we have to determine their nature {\em i.e.} the sign of $h'_0(u^*)$.

\begin{figure}[!ht]\label{fig:h1}
\centering
\vskip-0.2cm
\includegraphics[width=10cm]{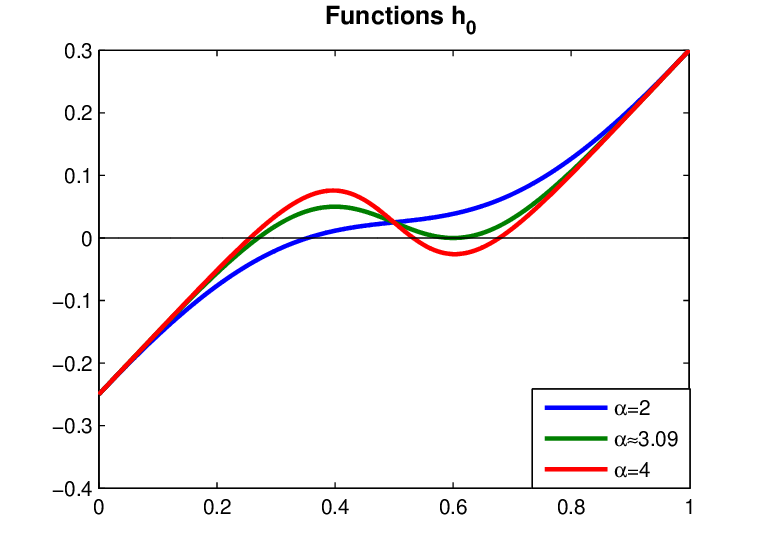}
\vskip-0.7cm
\caption{\em Examples of function $h_0$: $p_1=0.7$ and $p_2=0.75$ and  $f(u)=u^{\alpha}$. One zero for $\alpha=2$, two zeros for $\alpha\approx3.09$, three zeros for $\alpha=4$.}
\end{figure}

Let $u^*\!\in \{h_0=0\}\subset I^*$.  We deduce from~\eqref{eq:h0}  that the condition $h_0(u^*)=0$ reads
\[
(u-p_1)f(u)+(u-1-p_2)f(1-u)=0.
\]
Plugging this equality  into the expression~\eqref{eq:h'0} for $h'_0(u^*)$ shows that, for such equilibrium points $u^*$, 
\begin{equation}\label{eq:h'_0}
h_0'(u^*)=1-\frac{f'(u^*)(p_1-u^*)+f'(1-u^*)(u^*-1+p_2)}{f(u^*)+f(1-u^*)}.
\end{equation}
The function $h_0$ has at most three zeros owing to Proposition~\ref{PropEq}. Therefore, each such zero $u^*$ of $h_0$ has three possible ``status''  (see $e.g.$ Figure~\ref{fig:h1}):
\begin{itemize}
\item If there exists $\varepsilon_0>0$ such that $h_{0|(u^*-\varepsilon_0,u^*)}<0$  and $h_{0|(u^*,u^*+\varepsilon_0)}>0$, then $u^*$ is a   {\em stable} equilibrium point. This is the case $e.g.$ if $h'_0(u^*)>0$ or, more generally, if the first non zero derivative at $u^*$ has an odd order and is positive.
\item If there exists $\varepsilon_0>0$ such that $h_{0|(u^*-\varepsilon_0,u^*)}>0$ and $h_{0|(u^*,u^*+\varepsilon_0)}<0$, then $u^*$ is an {\em unstable} equilibrium point (or a repeller).  This is the case $e.g.$ if $h'_0(u^*)<0$ or if the first non-zero derivative at $u^*$ has an odd order and is negative.
\item If $h_0$ has a constant sign over an interval $(u^*-\varepsilon_0,u^*+\varepsilon_0)$, for some $\varepsilon_0>0$, then $u^*$ is a {\em semi-stable equilibrium point} (one-sided stable and one-sided unstable). This is the case if $h''_0(u^*)=0$ and the first non-zero derivative at $u^*$  occurs at an even order. The simplest case is when $h''_0(u^*)$ exists and is non-zero. 
\end{itemize}
Then it is an elementary exercise to derive the following proposition, if one keeps in mind that $h_0(0)=-(1-p_2)<0$ and $h_0(1)=1-p_1>0$. 

\begin{prop}\label{ThmAttract} Let $p_1,p_2\!\in (0,1)$. 

\smallskip
\noindent 	$(i)$  If $h_0$ has a unique equilibrium point, then it is stable.
	
\smallskip
\noindent $(ii)$  If $h_0$ has two equilibrium points, then one  is stable and one is semi-stable.
	
\smallskip
\noindent 	$(iii)$ If $h_0$ has three equilibrium points, then either
\begin{itemize}
\item the lowest and the highest ones are stable and the one in the middle is unstable.

or 

\item  the one in the middle is stable and the  other two are   semi-stable.
\end{itemize}
\end{prop}

 
\paragraph{Remark.} The second setting in item $(iii)$ does not correspond to a ``generic'' situation but rather to a limiting case of a situation with $5$ equilibrium points and $3$  extrema. Thus, it never occurs with the family $f_{\alpha}(y)= y^{\alpha}$, $\alpha>0$,  of skewing functions. That is why we will not consider this configuration  in what follows. 

\paragraph{Example of phase transition~(II).} We still consider the family  $f(u)=u^{\alpha}$, $\alpha>1$.  The ``bifurcation'' has been established in part~I. Now we want to elucidate what happens at $\alpha_0$. The critical two equilibrium case occurs at $\alpha_0\approx3.09$. By continuity of $(\alpha,u)\mapsto u^{\alpha}$ and monotony in $\alpha$ and $u$, it is clear that at least one~--~in fact exactly one~--~of the two equilibrium point $u^*$ satisfies $h_0(u^*)=0$ and $h_0$ does not change its sign around $u^*$ so that $h'_0(u^*)=0$. To determine  its status we have to look at the sign of the second derivative $h''_0(u^*)$. In fact elementary computations show that , for any $u\!\in[0,1]$ such that $h_0(u)=h'_0(u)=0$,
$$
h_0''(u)=(1-p_1-p_2)\alpha u^{\alpha-2}(1-u)^{\alpha-2}\frac{(\alpha-1)(1-2u)(u^{\alpha}+(1-u)^{\alpha})-2\alpha u(1-u)(u^{\alpha-1}-(1-u)^{\alpha-1})}{(u^{\alpha}+(1-u)^{\alpha})^3}>0.
$$
Hence, $u^*$ is semi-stable (stable on the right and unstable on the left).


\begin{figure}[!ht]
\centering
\vskip-0.2cm
\includegraphics[width=10cm]{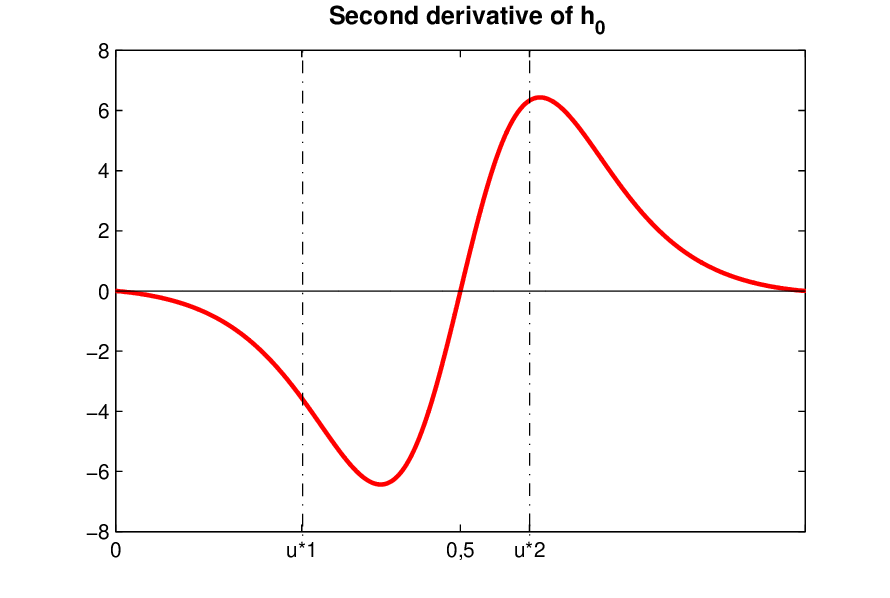}
\vskip-0.7cm
\caption{Second derivative of $h_0$ for $f(u)=u^{3.09}$ with $p_1=0.7$ and $p_2=0.75$.}
\end{figure}

\medskip
It remains to show that the algorithm does not converge towards the repulsive equilibrium point, denoted by $\hat u$ in what follows. To show that there is an excitation in the repulsive direction, we have to prove that assumption~(\ref{BruitMgle}) holds  (see Theorem~\ref{ThmPiege} in the  Appendix).

\begin{prop}\label{ThmRepuls}
Let $\hat u$ be an unstable equilibrium point for $h_0$  satisfying  $h_0'(\hat u)<0$. Then
$$
\P\big(\widetilde{Y}_n\to(\hat u,1-\hat u)\big)=0.
$$
\end{prop}

\noindent {\bf Proof.} It follows from~\eqref{eq:h'0} that $h'_0$ is locally  Lipschitz continuous. We know that $\{\widetilde{Y}_n\to(\hat u,1-\hat u)\}=\{\widetilde{Y}^1_n\to \hat u\}$, so we may focus on the first component. We rely  on Theorem~\ref{ThmPiege} in the  Appendix and adopt its notation, namely
$$
\Delta M_{n+1}^{(r)}=\Delta M_{n+1}^1.
$$
Using Assumption {\bf (A1)}, we obtain
$$
\Delta M_{n+1}^{(r)}=D_{n+1}^{11}X_{n+1}^1+D_{n+1}^{12}X_{n+1}^2-\frac{H_{n+1}^{11}f(\widetilde{Y}^1_n)+H_{n+1}^{12}f(\widetilde{Y}^2_n)}{\Tr(\widetilde{f}(\widetilde{Y}_n))}.
$$
Therefore
\begin{eqnarray*}
	\E\left[\left\|\Delta M_{n+1}^{(r)}\right\|\big | \F_n\right]&=&\E\big[\left|\Delta M_{n+1}^1\right|\big | \F_n\big] \\
&=&\frac{f(\widetilde{Y}^1_n)}{\Tr(\widetilde{f}(\widetilde{Y}_n))}\E\left[\left|D_{n+1}^{11}-\frac{H_{n+1}^{11}f(\widetilde{Y}^1_n)+H_{n+1}^{12}f(\widetilde{Y}^2_n)}{\Tr(\widetilde{f}(\widetilde{Y}_n))}\right|\Big | \F_n\right]\\
	& &+\frac{f(\widetilde{Y}^2_n)}{\Tr(\widetilde{f}(\widetilde{Y}_n))}\E\left[\left|D_{n+1}^{12}-\frac{H_{n+1}^{11}f(\widetilde{Y}^1_n)+H_{n+1}^{12}f(\widetilde{Y}^2_n)}{\Tr(\widetilde{f}(\widetilde{Y}_n))}\right|\Big | \F_n\right].
\end{eqnarray*}
By Jensen's inequality applied to both conditional expectations in the right hand side, we obtain
\begin{eqnarray*}
	\E\left[\left\|\Delta M_{n+1}^{(r)}\right\| \Big | \F_n\right]&\geq&\frac{f(\widetilde{Y}^1_n)}{\Tr(\widetilde{f}(\widetilde{Y}_n))}\left|H_{n+1}^{11}-\frac{H_{n+1}^{11}f(\widetilde{Y}^1_n)+H_{n+1}^{12}f(\widetilde{Y}^2_n)}{\Tr(\widetilde{f}(\widetilde{Y}_n))}\right|\\
	& &+\frac{f(\widetilde{Y}^2_n)}{\Tr(\widetilde{f}(\widetilde{Y}_n))}\left|H_{n+1}^{12}-\frac{H_{n+1}^{11}f(\widetilde{Y}^1_n)+H_{n+1}^{12}f(\widetilde{Y}^2_n)}{\Tr(\widetilde{f}(\widetilde{Y}_n))}\right|.
\end{eqnarray*}
Owing to {\bf (A5)}$_v$, $H_{n+1}^{ij}\overset{a.s.}{\underset{n\to+\infty}{\longrightarrow}}H^{ij}$, where $H^{ii}=p_i$ and $H^{ij}=1-p_j$, $1\leq i,j\leq 2$. Furthermore, on $\hat{\cal Y}_{\infty}=\left\{\omega:\widetilde{Y}^1_n(\omega)\to\hat u\right\}$, 
$$
\frac{H_{n+1}^{11}f(\widetilde{Y}^1_n)+H_{n+1}^{12}f(\widetilde{Y}^2_n)}{\Tr(\widetilde{f}(\widetilde{Y}_n))}\overset{a.s.}{\underset{n\to+\infty}{\longrightarrow}}\hat u.
$$
Consequently,
$$
\frac{f(\widetilde{Y}^1_n)}{\Tr(\widetilde{f}(\widetilde{Y}_n))}\left|H_{n+1}^{11}-\frac{H_{n+1}^{11}f(\widetilde{Y}^1_n)+H_{n+1}^{12}f(\widetilde{Y}^2_n)}{\Tr(\widetilde{f}(\widetilde{Y}_n))}\right|\overset{a.s.}{\underset{n\to+\infty}{\longrightarrow}}\frac{f(\hat u)}{f(\hat u)+f(1-\hat u))}\left|p_1-\hat u\right|>0
$$
and	
$$\frac{f(\widetilde{Y}^2_n)}{\Tr(\widetilde{f}(\widetilde{Y}_n))}\left|H_{n+1}^{12}-\frac{H_{n+1}^{11}f(\widetilde{Y}^1_n)+H_{n+1}^{12}f(\widetilde{Y}^2_n)}{\Tr(\widetilde{f}(\widetilde{Y}_n))}\right|\overset{a.s.}{\underset{n\to+\infty}{\longrightarrow}}\frac{f(\hat u)}{f(\hat u)+f(1-\hat u))}\left|1-p_2-\hat u\right|>0
$$
since  $\hat u\!\in I^*$. Thus (\ref{BruitMgle}) is satisfied. Then, by using (\ref{A3}) and by applying Theorem~\ref{ThmPiege} in the  Appendix,
$\P(\hat{\cal Y}_{\infty})=0$. \hfill$\cqfd$

\bigskip
\noindent{\bf Example (bi-stochastic generating matrix $H$).} If $p_1=p_2=p\!\in (0,1)$ then ($H$ is bi-stochastic) $(\frac 12;\frac 12)$ is an equilibrium and $h'_0(1/2) = 1-\frac{f'(1/2)}{2f(1/2)}(2p-1)$.

\smallskip
\noindent -- If $f$ is convex and $2p-1\leq0$, then $h'_0\left(\frac12\right)>0$ and the equilibrium is unique (since $h'_0$ is increasing). 

\smallskip
\noindent -- If $f$ is convex and $2p-1>0$, then we have three possibilities:
\begin{itemize}
	\item If $p>\frac 12+\frac{f(1/2)}{f'(1/2)}$, then $h'_0(1/2)<0$, so that $\P(\hat{\cal Y}_{\infty})=0$.
	\item If $p=\frac 12+\frac{f(1/2)}{f'(1/2)}$, then $h'_0(1/2)=0$. But $h''_0(1/2)=0$, so we need to investigate higher order.
	\item If $p<\frac 12+\frac{f(1/2)}{f'(1/2)}$, then $h'_0(1/2)>0$, so $u^*$ is stable.
\end{itemize}
If $f(u)=u^{\alpha}$, the first case also reads $\alpha >\frac{1}{2p-1}$. For the second case, higher order analysis leads to $h^{(3)}_0(1/2)>0$, so that the equilibrium is stable.

\subsection{\texorpdfstring{$A.s.$}{A.s.} convergence}

\begin{theo}\label{ThmCvx}
Let $(Y_n)_{n\geq0}$ be the urn composition sequence defined by~(\ref{dynamic})-(\ref{ConstructX}). Under the assumptions {\bf (A1)}, {\bf (A2)}, {\bf (A3)} and the $f$-skewed drawing rule,
	
	\smallskip \noindent	$(a)$ $\displaystyle\frac{Y_n}{\Tr(Y_n)}\overset{a.s.}{\underset{n\rightarrow+\infty}{\longrightarrow}}(u^*,1-u^*)\in{\cal E}_2$ ({\it i.e.} $u^*\!\in\{ h_0=0\}$ and $u^*$ is not unstable for $ODE_{h_0}$).
	
	\smallskip \noindent $(b)$ $\displaystyle \widetilde{N}_n= \frac 1n \sum_{k=1}^n X_k\overset{a.s.}{\underset{n\rightarrow+\infty}{\longrightarrow}}\displaystyle\begin{pmatrix} \frac{f(u^*)}{f(u^*)+f(1-u^*))} \cr  \frac{f(u^*)}{f(u^*)+f(1-u^*))}\end{pmatrix}$.
\end{theo}

\noindent {\bf Proof.}  First, we will prove that $(a)\Rightarrow(b)$, then   $(a)$. 

\smallskip
\noindent $(a)\Rightarrow(b).$ For every $n\ge 1$, we have
$$
\E\left[X_n\,|\,\F_{n-1}\right]=\sum_{i=1}^d\frac{f(\widetilde{Y}_{n-1}^i)}{\Tr(\widetilde{f}(\widetilde{Y}_{n-1}))}e^i=\frac{\widetilde{f}(\widetilde{Y}_{n-1})}{\Tr(\widetilde{f}(\widetilde{Y}_{n-1}))}
$$ 
and, by construction $\left\|X_n\right\|^2=1$, so that $\E\left[\left\|X_n\right\|^2\,|\,\F_{n-1}\right]=1$. Hence, the martingale
$$
\widetilde{M}_n=\sum_{k=1}^n\frac{X_k-\E\left[X_k\,|\,\F_{k-1}\right]}{k}\overset{a.s.\&\, L^2}{\underset{n\rightarrow+\infty}{\longrightarrow}}\widetilde{M}_{\infty}\in L^2.
$$
Finally, it follows from the Kronecker Lemma that
$$
\frac{1}{n}\sum_{k=1}^nX_k-\frac{1}{n}\sum_{k=1}^n\frac{\widetilde{f}(\widetilde{Y}_{k-1})}{\Tr(\widetilde{f}(\widetilde{Y}_{k-1}))}\overset{a.s.}{\underset{n\rightarrow+\infty}{\longrightarrow}}0.
$$
This proves the announced implication owing  to the C\'esaro Lemma since $\frac{\widetilde f}{\Tr(\widetilde f)}$ is continuous at every point of ${\cal S}_2$.

\smallskip
\noindent  $(a)$ The algorithm is bounded by construction since it is $[0,1]$-valued. Assumption {\bf (A2)} implies that $\sup_{n\geq 0}\E\left[\left\|\Delta M_{n+1}\right\|^2\,|\,\F_n\right]<+\infty$ $a.s.$ and Assumption {\bf (A3)} implies that
$r_n\overset{a.s.}{\underset{n\rightarrow+\infty}{\longrightarrow}}0$. The set $\{h=0\}$ is finite hence $\widetilde{Y}_n$ $a.s.$ converges toward a zero of $h$ (see Theorem~\ref{ThmODE}$(c)$ in the appendix). Moreover, it follows from Proposition~\ref{ThmAttract} devoted to attractiveness that this zero cannot be a repulsive  (a point at which $h_0'$ is negative).


%

\section{Weak rate of convergence}

\noindent To establish a $CLT$ for the sequence $(\widetilde{Y}_n)_{n\geq0}$ on  a convergence event $\big\{\widetilde Y_n \to y^*\big\}$, we need to make  the following additional assumptions:

\medskip
\noindent {\bf (A4)} The addition rules $D_n$ $a.s.$ satisfy on the event $\big\{\widetilde Y_n \to y^*\big\}$
\begin{equation*}\label{A4}
	\forall k\!\in \{1,\ldots,d\}, \quad\left\{
	\begin{array}{lll}
	 	(i)&\sup_{n\geq1}\E\left[\|D^{\cdot k}_n\|^{2+\delta}\,|\,\F_{n-1}\right]\leq \kappa< +\infty & \mbox{for a  $\delta>0$,} \\
		(ii)&\E\left[D^{\cdot k}_n(D^{\cdot k}_n)^t\,|\,\F_{n-1}\right]\underset{n\rightarrow+\infty}{\longrightarrow} C^k_{y^*}, & \\
	\end{array}\right.
\end{equation*}
where $C_{y^*}^k=(C_{_{y^*},ij}^k)_{1\leq i,j\leq d}$, $k=1,\ldots,d$, are $d\times d$ positive definite matrices. 

Note that {\bf (A4)} $\Rightarrow${\bf (A2)} since $\E\left[\|D^{\cdot k}_n\|^2\,|\,\F_{n-1}\right]\leq\left(\E\left[\|D^{\cdot k}_n\|^{2+\delta}\,|\,\F_{n-1}\right]\right)^{\frac{2}{2+\delta}}$.

\medskip
Also note that, if the matrices $D_n $ are {\em deterministic}  (hence co-stochastic) then  $C^k =\ds  \lim_n D_n^{.k}\otimes D_n^{.k}$ so that $ \mbox{\bf 1}^tC^k  \mbox{\bf 1}=1$ for every $ k\!\in \{1,\ldots,d\}$.

\medskip
Let $v=(v_n)_{n\ge 1}$ be  a sequence of positive real numbers. 

\medskip
\noindent {\bf (A5)}$_v$ The matrices $H_n$ and  $H$ satisfy on the event $\big\{\widetilde Y_n \to y^*\big\}$
\begin{equation}\label{A5}
	n\,v_n\,\E\left[|\!|\!|H_n-H|\!|\!|^2\right]\underset{n\rightarrow+\infty}{\longrightarrow}0.
\end{equation}

To establish the weak rate of convergence we cannot restrict ourselves to simplex as we essentially did for the convergence result since we must take into account the rate of convergence of the algorithm $\widetilde Y_n $ toward the simplex which is itself non trivial in general. Actually, we saw in~\eqref{ASTrace} that
\[
\Tr\big(\widetilde Y_n \big)-1= \frac{\Tr(M_n)}{n+\Tr(Y_0)}.
\]

Elementary computations show, under Assumptions~{\bf (A1)}, {\bf (A3)}, {\bf (A4)}-$(i)$ and if $\widetilde Y_n \to y^*$, that
\[
\E\left[\Tr(M_n)^2\,|\, \F_{n-1}\right]\to \sigma^2(y^*) =   \mbox{\bf 1}^t \frac{\sum_{k=1}^df(y^k_*)C_{y^*}^k}{\Tr(\widetilde f(y^*))} \mbox{\bf 1}-1.
\]
Note that, if the matrices $D_n $ are deterministic then $\sigma^2(y^*) =0$ which is expected since in that case $\widetilde Y_n $ is ${\cal S}_d$-valued. Owing to Condition~{\bf (A4)}
Lindeberg's {\it CLT} for arrays of martingale applies  to $\big(\frac{\Tr(M_\ell)}{\sqrt{n}}\big)_{1\le \ell\le n}$ (see {\it e.g.} Corollary~3.1, p.58 in~\cite{HallHeyde}) if one keeps in mind that or that {\bf (A4)}-$(ii)$ implies the usual condition as a straightforward application of Markov inequality). Consequently
\[
\sqrt{n}\Big(\Tr\big(\widetilde Y_n \big)-1\Big)= \frac{n}{n+\Tr(Y_0)}\frac{\Tr(M_n)}{\sqrt{n}} \stackrel{{\cal L}_{stably}}{\longrightarrow} {\cal N}\big(0;\sigma^*(y^*)\big)\;\mbox{ on the event $\big\{\widetilde Y_n \to y^*\big\}$}.
\]

\subsection{Strictly concave case with irreducible bi-stochastic limiting generating matrix}

Assume $H$ is bi-stochastic and irreducible and the skew-function is strictly concave. We know from Proposition~\ref{prop:convy(d)} that $\widetilde Y_n \to y(d)= \frac{\mbox{ 1}}{d}$ $a.s.$ On the other hand we know from (the proof of) Lemma~\ref{lem:jacphi} that $J_{\varphi}(y(d))= \lambda_1I_{d|\mbox{\bf 1}^{\perp}}$ so that
\[
J_h\big(y(d)\big)= I_d -\lambda_1H_{|\mbox{\bf 1}^{\perp}} \quad\mbox{still with}\quad \lambda_1= \frac{f'(1/d)}{d\cdot f(1/d)}<1.
\]

As $H$ is (co-)stochastic,  we know that $1$ is its eigenvalue with the highest real part so that the eigenvalue of $J_h(y(d))$ with the lowest real part is $1-\lambda_1>0$.
\begin{prop}\label{prop:concavebistochTCL} Assume that the skewed drawing rule $f$ is strictly concave and $H$ is irreducible and bi-stochastic. Assume {\bf (A1)}, {\bf (A3)}, {\bf (A4)} and {\bf (A5)}$_v$ hold. 

\smallskip
\noindent $(a)$ If $f'(1/d) <\frac d2 f(1/d)$ and  {\bf (A5)}$_v$ holds with $v_n=1$, then 
\[
\sqrt{n}\big(\widetilde Y_n-y(d)\big)\stackrel{{\cal L}_{stably}}{\longrightarrow} {\cal N}\big(0;\Sigma^*\big) \quad \mbox{ with }\quad \Sigma^*= \int_0^{+\infty}e^{-u(J_h(y^*)-\frac{I_2}{2})^t}\Gamma^* e^{-u(J_h(y^*)-\frac{I_2}{2})}du
\]
with
\begin{equation}\label{Gamma}
\Gamma^*= \frac{\sum_{k=1}^d f(y^{*,k})C^k_{y^*}}{\Tr\big(\widetilde f(y^*)\big)}-y^*(y^*)^t .
\end{equation}
\noindent $(b)$ If $f'(1/d) = \frac d2 f(1/d)$ and  {\bf (A5)}$_v$ holds with $v_n=\log n$, then 
\[
\sqrt{\frac{n}{\log n}}\big(\widetilde Y_n-y(d)\big)\stackrel{{\cal L}_{stably}}{\longrightarrow} {\cal N}\big(0;\Sigma^*\big) 
\]
with 
\begin{equation}\label{eq:Sigma^*}
\Sigma^*=\lim_n \frac{1}{n}\int_0^n e^{-u\left(J_h(y^*)-\frac{I_d}{2}\right)^t} \Gamma^* e^{-u\left(J_h(y^*)-\frac{I_d}{2}\right)}du.
\end{equation}

\smallskip
\noindent $(c)$ If $f'(1/d) > \frac d2 f(1/d)$ and  {\bf (A5)}$_v$ holds with $v_n=n ^{1-2\lambda_1+\eta}$ for some $\eta>0$, then 
\[
n^{1-\lambda_1}\big(\widetilde Y_n-y(d)\big)\mbox{ converges toward a finite random variable.}
\]
 \end{prop}

\noindent {\bf Proof.} First note that $\widetilde Y_n \to y(d)$ $a.s.$  under the assumptions. We will check the three assumptions of Theorem~\ref{ThmCLT} ($CLT$ for $SA$ algorithms) recalled in the  Appendix.
The parameter $\Lambda$ which rules the regime of the rate is given here by $\Lambda = 1-\lambda_1$ which justifies the above three cases.
%
Secondly Assumption~{\bf (A4)} ensures that Condition~(\ref{HypDM}) is satisfied since
$$
\sup_{n\geq1}\E\left[\left\|\Delta M_n\right\|^{2+\delta}\,|\,\F_{n-1}\right]<+\infty \quad a.s.\quad \mbox{and} \quad
\E\left[\Delta M_n\Delta M_n^t\left.\right|\F_{n-1}\right]\overset{a.s.}{\underset{n\rightarrow+\infty}{\longrightarrow}}\Gamma^*.
$$
To be more precise on the convergence on the right-hand side
\begin{eqnarray*}
\E\left[\Delta M_{n+1}\Delta M_{n+1}^t\,|\,\F_n\right]\!\!&\!\!=\!\!&\!\!\sum_{k=1}^d\P(X_{n+1}=e^k\,|\,\F_n) \! \left(\E\left[D_{n+1}^{\cdot k}(D_{n+1}^{\cdot k})^t\,|\,\F_n\right] \right.\\
 \!\!&\!\!  \!\!&\!\!\left.-\E\left[D_{n+1}X_{n+1}\,|\,\F_n\right]\E\left[D_{n+1}X_{n+1}\,|\,\F_n\right]^t\right) \\
	 \!\!&\!\!=\!\!&\!\!\sum_{k=1}^d\frac{f(\widetilde{Y}_n^q)}{\Tr(\widetilde{f}(\widetilde{Y}_n))}\E\left(D_{n+1}^{\cdot k}(D_{n+1}^{\cdot k})^t\,|\,\F_n\right)-\left(H_{n+1}\frac{\widetilde{f}(\widetilde{Y}_n)}{\Tr(\widetilde{f}(\widetilde{Y}_n))}\right)\left(H_{n+1}\frac{\widetilde{f}(\widetilde{Y}_n)}{\Tr(\widetilde{f}(\widetilde{Y}_n))}\right)^t \\
	 \!\!&\!\!\overset{a.s.}{\underset{n\rightarrow+\infty}{\longrightarrow}}\!\!&\!\!\Gamma^*=\frac{\sum_{k=1}^df(y^{*k})C^k_{y^*}}{\Tr(\widetilde{f}(y^*))}-y^*(y^*)^t.
\end{eqnarray*}
\noindent Finally, using {\bf (A5)}$_v$ with the appropriate sequence $(v_n)_{n\ge1}$, one proves in the three cases that the remainder sequence $(r_n)_{n\geq1}$ defined by~\eqref{resteConvex} satisfies~\eqref{HypReste} since one checks that $\frac{\widetilde{f}(\widetilde{Y}_n)}{\Tr(\widetilde{f}(\widetilde{Y}_n))}$ is bounded.~$\,\Box$
 
\subsection{Back to the ``convex'' bi-dimensional case}

As for the weak rate, we need to deal with the original $2$-dimensional algorithm on its whole.
%
\begin{theo}\label{Thm2}
Assume {\bf (A1)}, {\bf (A3)}, {\bf (A4)} hold. Every equilibrium point $y^*\!\in {\cal S}_2$ is of the form $y^*=(u^*,1-u^*)$ where  $u^*$ is solution to~\eqref{EConvex} and lies in $I^*= [p_1\wedge (1-p_2), p_1\vee (1-p_2)]$, and 
$$
\mbox{Sp}\big(J_h(y^*)\big)=\big\{1,1-\rho^*\big\},
$$
where
$$
\rho^*= \rho(y^*)=\frac{f'(u^{*})(p_1-u^{*})+f'(1-u^{*})(u^{*}-1+p_2)}{f(u^{*})+f(1-u^{*})}.
$$

\noindent $(a)$ If $p_1+p_2\leq1$    and {\bf (A5)}$_v$ holds with $v_n=1$, $n\geq1$, then $y^*$ is unique, stable,   $Y_n\to y^*$ $a.s.$ owing to Theorem~\ref{ThmCvx} $(a)$ and 
	$$
	\sqrt{n}\left(\widetilde{Y}_n-y^*\right)\overset{\L_{stably}}{\underset{n\rightarrow+\infty}{\longrightarrow}}{\cal N}\left(0,\Sigma^*\right)\quad \mbox{ with }\quad \Sigma^*=\int_0^{+\infty}e^{-u(J_h(y^*)-\frac{I_2}{2})^t}\Gamma^* e^{-u(J_h(y^*)-\frac{I_2}{2})}du
	$$
\begin{equation*}
\mbox{and  }\quad\Gamma^*=\frac{f(u^{*})C^1+f(1-u^{*})C^2}{\Tr(\widetilde{f}(y^*))}-y^*(y^*)^t.
\end{equation*}

\noindent $(b)$  If $p_1+p_2>1$,  we have three possible rates of convergence on an event $\big\{\widetilde Y_n \to y^*\big\}$, where $y^*\!\in {\cal E}_2$ is not unstable, depending on $\rho^*= \lambda(y^*)$:
	\begin{enumerate}
		\item[(i)] If $0<\rho^*<\frac12$ and    {\bf (A5)}$_v$ holds with  $v_n=1$, $n\geq1$, then 
			$$
			\sqrt{n}\left(\widetilde{Y}_n-y^*\right)\overset{\L_{stably}}{\underset{n\rightarrow+\infty}{\longrightarrow}}{\cal N}\left(0,\Sigma^*\right) \quad\mbox{on}\quad \big\{\widetilde Y_n \to y^*\big\},
			$$
			where $\Sigma^*$ is formally defined like in item~$(a)$. 
		\item[(ii)] If $\rho^*=\frac12$  and {\bf (A5)}$_v$ holds with $v_n =\log n$, $n\ge 1$, then 
			$$
			\sqrt{\frac{n}{\log n}}\left(\widetilde{Y}_n-y^*\right)\overset{\L}{\underset{n\rightarrow+\infty}{\longrightarrow}}{\cal N}\left(0,\Sigma^*\right)\quad\mbox{ where $\Sigma^*$ is given by~\eqref{eq:Sigma^*}.}
			$$
			
		\item[(iii)] If $\frac12<\rho^*<1$  and {\bf (A5)}$_v$ holds  with $v_n=n^{1-2\rho^*+\eta}$, $\eta>0$, then $n^{\rho^*}\big(\widetilde{Y}_n-y^*\big)\to \Upsilon$ $a.s.$ converges as $n\to+\infty$ where $\Upsilon$ is a  finite random variable.
	\end{enumerate}
\end{theo}

\paragraph{Remark.} 
\noindent $\bullet$ The condition $0<\rho^*<\frac12$ is satisfied as soon as 
$$\left\{\begin{array}{rl}
	\displaystyle\frac{(f'(1-p_1)+f'(1-p_2))(p_1+p_2-1)}{f(1-p_1)+f(1-p_2)}<\frac12 & \mbox{if $f$ is concave}, \\
	\\
	\displaystyle\frac{(f'(p_1)+f'(p_2))(p_1+p_2-1)}{f(1-p_1)+f(1-p_2)}<\frac12 & \mbox{if $f$ is convex}, \\
\end{array}	
\right.$$	
by using the monotonicity of $f$ and $f'$ and that $u^*\in(1-p_2,p_1)$.

\medskip
\noindent$\bullet$ If $f(y)=y$, then $y^*$ is unique, is given by 
\[
y^*= \Big(\frac{1-p_2}{2-p_1-p_2},\frac{1-p_1}{2-p_1-p_2} \Big) \quad \mbox{ and }\quad \rho^*= p_1+p_2-1.
\]
Thus, if $p_1+p_2<\frac32$, then the recursive procedure~(\ref{ASNLConvex}) satisfies a regular $CLT$; if $p_1+p_2=\frac32$, (\ref{ASNLConvex}) satisfies Theorem~\ref{Thm2}$(b)$-$(ii)$; and if $p_1+p_2>\frac32$, (\ref{ASNLConvex}) admits an $a.s.$-rate of convergence. 

\medskip\noindent $\bullet$ In~\cite{ChaMaiPou,ChaPouSah}, the properties of the random variable $\Upsilon$ are deeply investigated in the more standard framework of P\'olya's urn with deterministic addition rule matrix. It is shown to be solution to a smoothing equation obtained by a smart decomposition of the urn into canonical components.  Thus, it is proved that its distribution is characterized by its moments. It is clear that such results are out of reach of standard $SA$ techniques although it would be challenging to check whether similar results about $\Upsilon$ in the randomized  and nonlinear framework are true. 

\bigskip
\noindent {\bf Proof.} We again rely on Theorem~\ref{ThmCLT} from the Appendix. Elementary though tedious computations yield the following formula  for the Jacobian $J_h(y)$ of $h$ at $y\!\in \R_+^d\setminus\{0\}$. We obtain
$$
J_h(y)=\begin{pmatrix} 1+\frac{f'(y^1)}{f(y^1)+f(y^2)}\left(\frac{p_1f(y^1)+(1-p_2)f(y^2)}{f(y^1)+f(y^2)}-p_1\right) & \frac{f'(y^2)}{f(y^1)+f(y^2)}\left(\frac{p_1f(y^1)+(1-p_2)f(y^2)}{f(y^1)+f(y^2)}-(1-p_2)\right) \cr \cr \frac{f'(y^1)}{f(y^1)+f(y^2)}\left(\frac{(1-p_1)f(y^1)+p_2f(y^2)}{f(y^1)+f(y^2)}-(1-p_1)\right) & 1+\frac{f'(y^2)}{f(y^1)+f(y^2)}\left(\frac{(1-p_1)f(y^1)+p_2f(y^2)}{f(y^1)+f(y^2)}-p_2\right)\end{pmatrix}.
$$
As all equilibrium points $y^*$ lie in the simplex ${\cal S}_2$, we have  $y^{*2}=1-y^{*1}$. Combined with the constraint $h(y^*)=0$ (see~\eqref{EConvex}), we finally obtain the following formula only true at equilibrium points:
$$
J_h(y^*)=\begin{pmatrix} 1+\frac{f'(y^{*1})}{f(y^{*1})+f(1-y^{*1})}\left(y^{*1}-p_1\right) & \frac{f'(1-y^{*1})}{f(y^{*1})+f(1-y^{*1})}\left(y^{*1}-(1-p_2)\right) \cr \cr \frac{f'(y^{*1})}{f(y^{*1})+f(1-y^{*1})}\left(p_1-y^{*1}\right) & 1+\frac{f'(1-y^{*1})}{f(y^{*1})+f(1-y^{*1})}\left(1-p_2-y^{*1}\right)\end{pmatrix}.
$$
Then, one easily checks  that  the spectrum of $J_h(y^*)$ is real given by 
$$
\mbox{Sp}\left(J_h(y^*)\right)=\left\{1,1-\rho^*
\right\}.
$$

\noindent $(a)$ When $p_1\le 1-p_2$ , we know from Proposition~\ref{PropEq}$(a)$  that the equilibrium point $y^*$ is unique and $u^*$ lies in $I^*$. Hence $\rho^*\le 0$ so that $1$ is the lowest eigenvalue of $J_h(y^*)$. Consequently, we are in the ``regular'' case of the $CLT$ for $SA$ (Theorem~\ref{ThmCLT}$(a)$ in the Appendix) since $1>\frac 12$.
%
%
Then, following the lines of the proof of Proposition~\ref{prop:concavebistochTCL}, we check that Assumption~{\bf (A4)} ensures that Condition~(\ref{HypDM}) is satisfied 
with $\Gamma^*= \displaystyle \frac{f(y^{*1})C^1+f(1-y^{*1})C^2}{\Tr(\widetilde{f}(y^*))}-y^*(y^*)^t$. Finally, using {\bf (A5)}$_v$, the remainder sequence $(r_n)_{n\geq1}$ defined by~\eqref{resteConvex} satisfies~(\ref{HypReste}) since $\frac{\widetilde{f}(\widetilde{Y}_n)}{\Tr(\widetilde{f}(\widetilde{Y}_n))}$ is bounded. 

\smallskip
\noindent $(b)$ If $p_1+p_2>1$, then $\rho^*>0$ which explains the three cases. The rest of the proof is the same as above, given the convergence event $\{\widetilde Y_n\to y^*\}$.  \hfill$\cqfd$

\section{P\'olya urn with concave reinforced drawing rule: a bandit approach}\label{quatre}

By {\em P\'olya urn}, we mean in this section that the matrices $D_n$ involved in the drawing rule all satisfy $D_n=I_d$ , $n\geq1$. Moreover we assume that the drawing rue is still skewed  following~(\ref{LoiXFConvex})  where the function $f$ is concave/convex and that the initial urn composition vector $Y_0\!\in \R_+^d\setminus\{0\}$.  Note that when $f(u)=u$, then the urn dynamics is that of a regular P\'olya urn with $d$ colors. 

In such a framework, $H=H_n=I_d$, $n\geq1$, therefore $H$ is no more irreducible and we cannot use the results proved in Sections~\ref{deux} and~\ref{trois}. 

We still normalize $Y_n$ by setting $\widetilde{Y}_n:=\frac{Y_n}{n+\Tr(Y_0)}$, $n\geq0$. The sequence $(\widetilde{Y}_n)_{n\geq0}$ satisfies the following recursive stochastic algorithm (obvious consequence of~(\ref{ASC2}))
\begin{equation}\label{ASBandit}	
\widetilde{Y}_{n+1}=\widetilde{Y}_n-\frac{1}{n+1+\Tr(Y_0)}\left(\widetilde{Y}_n-\frac{\widetilde{f}(\widetilde{Y}_n)}{\Tr\left(\widetilde{f}(\widetilde{Y}_n)\right)}\right)+\frac{1}{n+1+\Tr(Y_0)}\Delta M_{n+1}, \quad n\geq1,
\end{equation}	                   
where 
\begin{equation}\label{MgleBandit}
\Delta M_{n+1}:=X_{n+1}-\E\left[X_{n+1}\,|\,\F_{n}\right]
\end{equation}
is a true  $(\F_{n})_{n\geq0}$ martingale increment. Let us remark that, in this setting, 
$$
\frac{\Tr(Y_n)}{n+\Tr(Y_0)}=\frac{n+\Tr(Y_0)}{n+\Tr(Y_0)}=1,\quad n\ge 0,
$$ 
so that the sequence 
$$
\widetilde{Y}_n\!\in {\cal S}_d,\quad \mbox{for every $n\ge 0$},
$$
 since it has non-negative components.

 
\noindent The special case of a linear drawing rule  $f(y)=y$ is entirely elucidated by the celebrated Athreya theorem  recalled below for completeness.
\begin{theo}[Athreya's Theorem, see~\cite{AthKar2}]\label{ThmfId}
Let $(Y_n)_{n\geq0}$ be the urn composition sequence defined by~(\ref{dynamic}) and (\ref{ConstructX}) with $D_n=I_d$, $n\geq1$, and a  linear drawing rule ($i.e.$~(\ref{LoiXFConvex}) with $f(u)=u$). Then, if $Y_0\!\in \R^d\setminus\{0\}$ is deterministic,   there exists a random vector $\widetilde{Y}_{\infty}$ having values in the simplex ${\cal S}_d$ such that
$$\widetilde{Y}_n=\frac{Y_n}{\Tr(Y_n)}\overset{a.s.}{\underset{n\rightarrow+\infty}{\longrightarrow}}\widetilde{Y}_{\infty} \quad a.s.$$
Furthermore,
\begin{itemize}
	\item[(i)] $\widetilde{Y}_{\infty}$ has a Dirichlet distribution with parameter $Y_0$.

	\item[(ii)] In particular, if $d=2$, $\widetilde{Y}_{\infty}^1$ has a Beta distribution with parameters $Y_0^1$ and $Y_0^2$ (in particular, $\widetilde{Y}_{\infty}^1$ has a uniform distribution on $[0,1]$ if $Y_0^1=Y_0^2=1$).

\end{itemize}
\end{theo}
Now, we investigate the case $f\neq\operatorname{Id}_{\R_+}$ by borrowing tools to adaptive bandit models analysis (see~\cite{LamPagTar, LamPag2, LamPag3}). First note that, as $\widetilde Y_n $ lives in the simplex ${\cal S}_d$,  {\em the  function $f$  only needs to be defined on $[0,1]$}. Moreover, we will no longer ask $f$ to be  convex or concave on $(0,1)$ but require finiteness of the derivatives at $0$ and $1$.

\medskip
\begin{center}$f$   is continuous, non-decreasing, $f(0)\!=\!0$, $f(1)\!=\!1$, $f\!>\!0$ on $(0,1]$, with finite  right and left derivatives at $0$ and $1$.
\end{center}

\medskip 
A typical example could be $f(u)=4\left(u-\frac12\right)^3+\frac12$, $u\in[0,1]$.
We will see that  our study  requires new tools, especially a method to avoid noiseless repulsive equilibrium points (sometimes called {\em noiseless traps} in the $SA$ literature).

\begin{theo}\label{ThmBandit} $(a)$ Let $I\subsetneq\{1,\ldots,d\}$ be non-empty. If $f$ satisfies $f'_r(0)>|I|f\big(\frac{1}{|I|}\big)$, 
then, for every deterministic initial value such that $Y_0^j>0$ for some  $j\!\notin I$, 
$$
\P\left(\widetilde{Y}_{\infty}=\widetilde e_{I} \right)=0.
$$
\noindent \smallskip $(b)$ If $d=2$, the above conclusion still holds if $f'_r(0)=1$ and $f'_l(1)+\frac{f''_r(0)}{2}>1$.

\noindent \smallskip $(c)$ If $f$ is strictly concave then ${\cal E}_d=\left\{\widetilde e_{_I}, I\subset\{1,\ldots,d\},I\neq\varnothing\right\}$ by Proposition~\ref{prop:fixpt}$(c)$. Then, for every starting value $Y_0\!\in (0, +\infty)^d$,   
$$
\widetilde{Y}_n \overset{a.s.}{\longrightarrow}\widetilde e_{\{1,\ldots,d\}}=y(d)\; \mbox{ as }\; n\to+\infty.
$$
\end{theo}

\paragraph{Remarks.} 
%
$\bullet$ In claim  $(b)$, if $f(0)=0$, $f(1)=1$, $f_r'(0)=1$ and $f$ is convex or concave, then $f=\operatorname{Id}$, so this case is interesting only out of the concave/convex framework. 

\smallskip 
\noindent $\bullet$ If $Y_0^i=0$ for some $i\!\in \{1,\ldots,d\}$, then, as proved below, $Y_n^i=0$ for every $n\ge 0$. So, as soon as $Y_0\!\in \R_+^d\setminus\{0\}$, one may apply the above result $(c)$ to the urn restricted to $I'=\{i\!\in I,\; Y_0^i>0\}$ to prove that $\widetilde Y_n \to \widetilde e_{I'}$ as $n\to +\infty$. 

\bigskip
\noindent {\bf Proof.} $(a)$-$(b)$ It follows from~\eqref{dynamic} and the fact that $D_n\equiv I_d$ that, if $Y^i_0=0$, then, for every instant $n\ge 0$, $Y^i_n =0$. So, up to a reduction of the dimension $d$, we may always assume that all components $Y^i_0>0$. 
%
%
%

As a consequence we may assume that,  for every  $n\ge 0$, $\min_i\widetilde Y_n^i>0$. Our aim is to prove that $\P\left(\widetilde{Y}_{\infty}^j=0\right)=0$, for every $j\notin I$. To this end, we will show that $\{\widetilde Y^j_{\infty} = 0\}\subset \{\widetilde L_{\infty}=0\}$ where $\widetilde L_{\infty}$ is the terminal value of a non-negative martingale. Then we will apply an ``oracle'' inequality to this martingale. Without loss of generality,  we may assume that, up to a permutation, $1\notin I$ and $j=1$ in what follows. 

\medskip\noindent {\sc Step~1}:  First, we define the function $\widetilde{h}$ by 
$$
\widetilde{h}(y)=1-\frac{f(y^1)}{y^1\Tr(\tilde f(y))}\mbox{\bf 1}_{\{y^1\neq 0\}},\; y\!\in{\cal S}_d,
$$ 
which satisfies $\widetilde{h}(y)<1$, for every $y\!\in{\cal S}_d\setminus \{y:y^1=0\}$.
%

Starting from the dynamics of $\widetilde{Y}_n^1$ given by~(\ref{ASBandit}), we have, for every $n\geq0$,
\begin{eqnarray*}	
\widetilde{Y}_{n+1}^1&=&\widetilde{Y}_n^1-\frac{1}{n+1+\Tr(Y_0)}\left(\widetilde{Y}_n^1-\frac{f(\widetilde{Y}_n^1)}{\Tr\big(\widetilde f(\widetilde Y_n)\big)}\right)+\frac{1}{n+1+\Tr(Y_0)}\Delta M^1_{n+1} \\	
				&=&\widetilde{Y}_n^1\left(1-\frac{1}{n+1+\Tr(Y_0)} \widetilde h\big(\widetilde{Y}_n^1\big)\right)+\frac{1}{n+1+\Tr(Y_0)}\Delta M^1_{n+1}.
\end{eqnarray*}
 We derive that the (non-negative) sequence 
\begin{equation}\label{eq:Ln}
\widetilde{L}_n:=\frac{\widetilde{Y}_n^1}{\prod_{k=1}^n\left(1-\frac{1}{k+\Tr(Y_0)}\widetilde{h}(\widetilde{Y}_{k-1}^1)\right)}, \quad n\geq0,
\end{equation}
is a non-negative martingale satisfying the recursive equation $\widetilde L_0=\widetilde Y^1_0$ and 
$$
\widetilde{L}_{n+1}=\widetilde{L}_n+\frac{1}{n+1+\Tr(Y_0)}\frac{\Delta M^1_{n+1}}{\prod_{k=1}^{n+1}\left(1-\frac{1}{k+\Tr(Y_0)}\widetilde{h}(\widetilde{Y}_{k-1}^1)\right)}, \quad n\geq0.
$$

\noindent$\bullet$ If $ f'_r(0)>|I|f\big(\frac{1}{|I|}\big)$, then $\widetilde{h}(y) \longrightarrow \kappa:=1-\frac{f'_r(0)}{|I|f\big(\frac{1}{|I|}\big)}<0$ as $y\to \widetilde e_{_I}$. Therefore, on the event $\left\{ \widetilde Y_n \to \widetilde e_{_I}\right\}$, $ \widetilde Y^1_n \to 0$ so that $\widetilde{h}\left(\widetilde Y^1_{n-1}\right)\overset{a.s.}{\sim}\kappa<0$, which in turn implies $\prod_{k=1}^n\left(1-\frac{1}{k+\Tr(Y_0)}\widetilde{h}(\widetilde{Y}_{k-1}^1)\right)\overset{a.s.}{\longrightarrow}+\infty$. It follows from its definition in~\eqref{eq:Ln}  that $\widetilde{L}_n\overset{a.s.}{\longrightarrow}0$ on $\left\{\widetilde Y_n \to \widetilde e_{_I}\right\}$ since $0\le \widetilde Y^1_n \le \frac{\Tr(Y_n)}{n+\Tr(Y_0)}\overset{a.s.}{\longrightarrow} 1$. 
Consequently, 
$$
\left\{\widetilde Y_n \to \widetilde e_{_I}\right\} \subset\left\{\widetilde{L}_n\to 0\right\}.
$$

\noindent$\bullet$ (Case $d=2$) If $f'_r(0)=1$ and $f'_l(1)+\frac{f''_r(0)}{2}>1$, then $\widetilde{h}(y)<0$ for $y$ in the neighborhood of $0$. So, with in mind $\widetilde e_{\{2\}}=e^2$,  we still have  $\left\{\widetilde Y_n\to e^2\right\}= \left\{\widetilde Y^1_n\to0\right\}\subset\left\{\widetilde{L}_n\to0\right\}$.


\medskip
\noindent{\sc Step~2}: The end of the proof is based on the following (short) ``oracle'' lemma (see~\cite{LamPagTar}) reproduced here for the reader's convenience.
\begin{lem}[Oracle inequality]\label{LemMgle}
Let $(M_n)_{n\geq0}$ be a non-negative martingale. Then
$$
\forall n\geq0, \quad\P\big(M_{\infty}=0\,\big |\, \F_n\big)\leq\frac{\E\left[\Delta\left\langle M\right\rangle_{n+1}^{\infty}\Big | \F_n\right]}{M_n^2}.
$$
\end{lem}

\paragraph{Proof of Lemma~\ref{LemMgle}.} It is sufficient to observe that, for every $n\geq0$,
\begin{eqnarray*}
\P\big(M_{\infty}=0\,\big |\, \F_n\big)&=&\frac{\E\left[\mathds{1}_{\{M_{\infty}=0\}}M_n^2\Big | \F_n\right]}{M_n^2}\leq\frac{\E\left[\left(M_{\infty}-M_n\right)^2\Big | \F_n\right]}{M_n^2} =\frac{\E\left[\Delta\left\langle M\right\rangle_{n+1}^{\infty}\Big | \F_n\right]}{M_n^2}. \cqfd
\end{eqnarray*}

\bigskip
\noindent First we note that 
\[
\E\left[(\Delta\widetilde{L}_{n+1})^2\Big | \F_n\right]=\left(\frac{1}{n+1+\Tr(Y_0)}\right)^2\frac{\E\left[(\Delta M^1_{n+1})^2\Big | \F_n\right]}{\left(\prod_{k=1}^{n+1}\left(1-\frac{\widetilde{h}(\widetilde{Y}_{k-1}^1)}{k+\Tr(Y_0)}\right)\right)^2}
\]
and 	
\[
\E\left[\big(\Delta M^1_{n+1}\big)^2\Big | \F_n\right]= \frac{f(\widetilde Y^1_n)\big(\Tr\big(\widetilde f(\widetilde Y_n)\big)-f(\widetilde Y^1_n)\big)}{\Tr\big(\widetilde f(\widetilde Y_n)\big)^2}\le \frac{d-1}{f(1/d)^2}
\]
since $\Tr(f(y))\ge f\big(\max_i y_i\big)\ge f(1/d)$ and $\widetilde Y_n \!\in {\cal S}_d$. As a consequence	
\[	
\E\left[(\Delta\widetilde{L}_{n+1})^2\,\big |\, \F_n\right]= \frac{1}{(n+1+\Tr(Y_0))^2\left(\prod_{k=1}^{n+1}\left(1-\frac{\widetilde{h}(\widetilde{Y}_{k-1}^1)}{k+\Tr(Y_0)}\right)\right)^2}\frac{f(\widetilde{Y}_n^1)\big(\Tr\big(\widetilde f(\widetilde Y_n)\big)-  f(\widetilde{Y}_n^1)\big)}{\left(\Tr\big(\widetilde f(\widetilde Y_n)\big)\right)^2}.
\]
Then,  applying Lemma~\ref{LemMgle} to the non-negative martingale $(\widetilde{L}_n)_{n\ge 1}$ yields
\begin{eqnarray*}
 \P\big(\widetilde{L}_{\infty}=0\,\big |\, \F_n\big) &\leq&\frac{\E\left[\left.\Delta\langle\widetilde{L}\rangle_{n+1}^{\infty}\right|\F_n\right]}{\widetilde{L}_n^2} \\	                                   &=&\frac{1}{\widetilde{L}_n^2}\E\left[\sum_{k=n+1}^{\infty}\frac{F(\widetilde{Y}_{k-1}^1)}{(k+\Tr(Y_0))^2\left(\prod_{\ell=1}^k\left(1-\frac{\widetilde{h}(\widetilde{Y}_{\ell-1}^1)}{\ell+\Tr(Y_0)}\right)\right)^2}\Big|\F_n\right],
\end{eqnarray*}
where the function $F$,  defined by  $F(y):=\frac{f(y^1)\big(\Tr(\widetilde f(y))-f(y^1)\big)}{\left(\Tr(\widetilde f(y))\right)^2}$, $y\!\in (0,1]\times [0,1]^{d-1}$, is clearly non-negative and bounded by $\kappa_d:=(d-1)/f(1/d)^2$. 
Consequently,
\begin{align*}
 \P\big(\widetilde{L}_{\infty}=0\,&\big |\, \F_n\big)   \\                                
&\leq\frac{\kappa_d}{\widetilde{L}_n^2}\sum_{k=n+1}^{\infty}\frac{1}{(k+\Tr(Y_0))^2}\E\left[\underbrace{\frac{\widetilde{Y}_{k-1}^1}{\prod_{\ell=1}^{k-1}\left(1-\frac{\widetilde{h}(\widetilde{Y}_{\ell-1}^1)}{\ell+\Tr(Y_0)}\right)}}_{=\widetilde{L}_{k-1}}\frac{\left(1-\frac{1}{k+\Tr(Y_0)}\widetilde{h}(\widetilde{Y}_{k-1}^1)\right)^{-1}}{\prod_{\ell=1}^k\left(1-\frac{\widetilde{h}(\widetilde{Y}_{\ell-1}^1)}{\ell+\Tr(Y_0)}\right)}\,\Big|\,\F_n\right].
\end{align*}
Let $\widetilde{h}_{_+}:=\max(\widetilde{h},0)$, so that $\widetilde{h}\leq\widetilde{h}_{_+}\leq\big\|\widetilde{h}_{_+}\big\|_{\infty}$. First note that   $\widetilde{h}_{_+}(y)\le 1-\frac{f(y^1)}{d\cdot y^1}<1$, for every $y\in{\cal S}_d\setminus\{y^1=0\}$ since $\Tr(\tilde f(y))\le d$ and we assumed $f(y^1)>0$ on $(0,1]$. Moreover,  $\limsup_{y^1\to 0}\widetilde{h}(y) \le 1-f'_r(0)/d<1$ since $f'_r(0)>0$ by assumption. As a consequence $\|\widetilde{h}_{_+}\big\|_{\infty}<1$.

In $2$-dimensions, under the additional assumption in the critical case ($f'_r(0)=1$), the extension of the function $h$ over $[0,1]$ is negative so that $\widetilde{h}_{_+}\equiv 0$.  

\smallskip
Finally, we obtain 
$$
\left(1-\frac{\widetilde{h}(\widetilde{Y}_{k-1}^1)}{k+\Tr(Y_0)}\right)^{-1}\leq\left(1-\frac{\|\widetilde{h}_{_+}\big\|_{\infty}}{k+\Tr(Y_0)}\right)^{-1}, \quad k\geq1.
$$
Then, as $\E\left[\widetilde{L}_{k-1}\Big | \F_n\right]=\widetilde{L}_n$ for every $k\ge n+1$, since $\widetilde{M}$ is a $(\P,\F_n)$-martingale,
\begin{align}
\P(\widetilde{L}_{\infty}=0&\big | \F_n)  \nonumber \\  	                                   
&\leq\frac{\kappa_d}{\widetilde{L}_n^2}\sum_{k=n+1}^{\infty}\frac{1}{(k+\Tr(Y_0))^2}\frac{\widetilde{L}_n}{\left(1-\frac{\|\widetilde{h}_{_+}\|_{_\infty}}{k+\Tr(Y_0)}\right)\prod_{\ell=1}^n\left(1-\frac{\widetilde{h}(\widetilde{Y}_{\ell-1}^1)}{\ell+\Tr(Y_0)}\right)\prod_{\ell=n+1}^k\left(1-\frac{\|\widetilde{h}_{_+}\|_{_\infty}}{\ell+\Tr(Y_0)} \right)} \nonumber \\                       
&=\frac{\kappa_d}{\widetilde{L}_n\prod_{\ell=1}^n\left(1-\frac{\widetilde{h}(\widetilde{Y}_{\ell-1}^1)}{\ell+\Tr(Y_0)}\right)}\sum_{k=n+1}^{\infty} \frac{1}{\left(1-\frac{\|\widetilde{h}_{_+}\|_{_\infty}}{k+\Tr(Y_0)}\right)(k+\Tr(Y_0))^2\prod_{\ell=n+1}^k\left(1-\frac{\|\widetilde{h}_{_+}\|_{_\infty}}{\ell+\Tr(Y_0)} \right)}\nonumber \\
&\le \frac{\kappa_d}{C_{0}\widetilde{Y}^1_n}\sum_{k=n+1}^{\infty}\frac{1}{(k+\Tr(Y_0))^2}\exp\left(-\sum_{\ell=n+1}^{k}\log \left(1-\frac{\|\widetilde{h}_{_+}\|_{_\infty}}{\ell+\Tr(Y_0)}\right)\right)\label{eq:sumlog}	                                  
\end{align}
where $C_0= \frac{\Tr(Y_0)}{1+\Tr(Y_0)}\!\in (0,1)$ since $
1-\frac{\|\widetilde{h}_{_+}\|_{_\infty}}{k+\Tr(Y_0)}>1-\frac{1}{k+\Tr(Y_0)}\ge C_0$.

Note that $\log(1-u)\ge -\frac{u}{1-u_0}$, $u\!\in (0,u_0)$. Applying this inequality with $u_0= \frac{1}{n+1+\Tr(Y_0)}$ yields, for every $\ell\ge n+1$,
\[
\log \left(1-\frac{\|\widetilde{h}_{_+}\|_{_\infty}}{\ell+\Tr(Y_0)}\right)\ge -\left(1+\frac{1}{n+\Tr(Y_0)}\right)\frac{\|\widetilde{h}_{_+}\|_{_\infty}}{\ell+\Tr(Y_0)}.
\]
Hence
\begin{align*}
\sum_{\ell=n+1}^k\log \left(1-\frac{\|\widetilde{h}_{_+}\|_{_\infty}}{\ell+\Tr(Y_0)}\right)&\ge -\|\widetilde{h}_{_+}\|_{_\infty}\left(1+\frac{1}{n+\Tr(Y_0)}\right)\sum_{\ell= n+1}^k\frac{1}{\ell+\Tr(Y_0)}\\
&\ge  -\|\widetilde{h}_{_+}\|_{_\infty}\left(1+\frac{1}{n+\Tr(Y_0)}\right)\int_n^k\frac{du}{u+\Tr(Y_0)}\\
&=  -\|\widetilde{h}_{_+}\|_{_\infty}\left(1+\frac{1}{n+\Tr(Y_0)}\right)\log\left(\frac{k+\Tr(Y_0)}{n+\Tr(Y_0)}\right).
\end{align*}
Plugging this inequality into~\eqref{eq:sumlog} and using that $\widetilde Y_n^1= \frac{Y^1_{n}}{n+\Tr(Y_0)}$, 
\begin{align*}
\P\big(\widetilde{L}_{\infty}=0&\, |\, \F_n\big)  \leq\frac{\kappa_d(n+\Tr(Y_0))}{C_{0}Y^1_n}\sum_{k=n+1}^{+\infty}\frac{e^{(1+\frac{1}{n+\Tr(Y_0)})\|\widetilde{h}_{_+}\|_{_\infty}\log\big(\frac{k+\Tr(Y_0)}{n+\Tr(Y_0)}\big)}}{(k+\Tr(Y_0))^2}.
\end{align*}
As $\|\widetilde{h}_{_+}\|_{_\infty}<1$, there exists $n_0$ be such that for every $n\ge n_0$, $1-(1+\frac{1}{n+\Tr(Y_0)})\|\widetilde{h}_{_+}\|_{_\infty}>0$. Hence  

 \begin{eqnarray*}
\P\big(\widetilde{L}_{\infty}=0\, |\, \F_n\big) & \le &   \frac{\kappa_d(n+\Tr(Y_0))}{C_{0}Y^1_n}  \sum_{k=n+1}^{\infty}\frac{(n+\Tr(Y_0))^{-(1+\frac{1}{n+\Tr(Y_0)})\|\widetilde{h}_{_+}\|_{_\infty}}}{(k+\Tr(Y_0))^{2-(1+\frac{1}{n+\Tr(Y_0)})\|\widetilde{h}_{_+}\|_{_\infty}}} \\
	                                   &=&\frac{\kappa_d(n+\Tr(Y_0))^{1-(1+\frac{1}{n+\Tr(Y_0)})\|\widetilde{h}_{_+}\|_{_\infty}}}{C_0Y^1_n}\sum_{k=n+1}^{\infty}\frac{1}{(k+\Tr(Y_0))^{2-(1+\frac{1}{n+\Tr(Y_0)})\|\widetilde{h}_{_+}\|_{_\infty}}} \\
	                                   &\le &\frac{\kappa_d(n+\Tr(Y_0))^{1-(1+\frac{1}{n+\Tr(Y_0)})\|\widetilde{h}_{_+}\|_{_\infty}}}{C_0Y^1_n} \int_n^{+\infty}\frac{du}{(u+\Tr(Y_0))^{2-(1+\frac{1}{n+\Tr(Y_0)})\|\widetilde{h}_{_+}\|_{_\infty}}}\\
	                                   &\leq&\frac{\kappa_d}{C_0Y^1_n} \frac{1}{1-(1+\frac{1}{n+\Tr(Y_0)})\|\widetilde{h}_{_+}\|_{_\infty}}.
\end{eqnarray*}
Now, it remains to  prove that $Y^1_n\overset{a.s.}{\longrightarrow}+\infty$. Let $Y^1_{\infty}:=\lim_n Y^1_n$ (the components of $Y_n$ are non-decreasing). One checks that
%
$$\left\{Y^1_{\infty}<+\infty\right\}=\bigcup_{n\geq0}\bigcap_{k>n}\left\{U_k>\frac{f(\frac{Y^1_n}{k-1+\Tr(Y_0)})}{f(\frac{Y^1_n}{k-1+\Tr(Y_0)})+f(1-\frac{Y^1_n}{k-1+\Tr(Y_0)})}\right\},$$
then
$$\forall n\in\N, \quad \P\left(Y^1_{\infty}<+\infty\left.\right|Y^1_n=y\right)=\prod_{k>n}\left(1-\frac{f(y/(k-1+\Tr(Y_0)))}{f(y/(k-1+\Tr(Y_0)))+f(1-y/(k-1+\Tr(Y_0)))}\right)=0,$$
since $\displaystyle \sum_{k\ge1}\frac{f(y/k+\Tr(Y_0))}{f(y/k+\Tr(Y_0))+f(1-y/k+\Tr(Y_0))}=+\infty$ because $f'_r(0)>0$. Therefore $Y^1_{\infty}=+\infty $ $a.s.$ 

On the other hand,  the closed martingale $\P\big(\widetilde{L}_{\infty}=0|\, \F_n\big)\to \mbox{\bf 1}_{\{\widetilde{L}_{\infty}=0\}} $ $a.s.$ and in $L^1$ so that 
$$
\mbox{\bf 1}_{\{\widetilde{L}_{\infty}=0\}} = 0\quad a.s.\quad \mbox{ i.e.}\quad\P(\widetilde L_n\to 0)=  \P(\widetilde{L}_{\infty}=0)=0
$$
which in turn implies that $\P\left(\widetilde{Y}_{\infty}^1=0\right)$ since it was proved in {\sc Step~1} that $\{\widetilde{Y}_{n}^1\to 0\}\subset \{\widetilde L_n\to0\}$.

\smallskip
\noindent $(c)$   As $Y_0^i>0$, $i\!\in \{1,\ldots,d\}$, we derive from what precedes that $\P\big(\widetilde Y_n \to \partial{\cal S}_d\big)=0$. As a consequence, following Theorem~\ref{ThmODE}, $\P(d\omega)$-$a.s.$, the compact connected flow invariant set $\Theta^{\infty}(\omega)$ of limiting values of $(\widetilde Y_n(\omega))_{n\ge 0}$ is a minimal connected attractor    of   $ODE_h$ in $\overset{\circ}{{\cal S}_d}$.  We know from Proposition~\ref{pro:stabinstab}$(b)$ that $y(d)=\frac 1d \mbox{\bf 1}$ is a uniformly attracting point $ODE_h$ and from Proposition~\ref{prop:flot}$(b)$  that the flow of $ODE_h$ $(y(y_0,t))_{t\ge 0, y_0 \in\overset{\circ}{{\cal S}_d}}$  converges toward $\frac 1d \mbox{\bf 1}$, so it converges uniformly with respect to $y_0\!\in \Theta^{\infty}(\omega)$. Consequently, one concludes by  Theorem~\ref{ThmODE}  that  $\Theta^{\infty}(\omega)= \{y(d)\}$ (otherwise it would have an internal attractor). Hence $\widetilde{Y}_n\overset{a.s.}{\longrightarrow} y(d)$.\hfill$\cqfd$
%


\section{Applications}\label{cinq}

\subsection{A drawing rule based on  a function with regular variation}\label{sec:regvar}

Let define the law of the drawings as follows
\begin{equation}\label{LoiXFRegular}
\forall 1\leq i\leq d, \quad \P(X_{n+1}=e^i\,|\, \F_n)=\frac{f(Y_n^i)}{\sum_{j=1}^df(Y_n^j)}, \quad n\ge 0,
\end{equation}
where $f$ has regular variation with index $\alpha>0$ in the sense that   for every $t>0$, $\displaystyle\frac{f(tx)}{f(x)}\underset{x\to\infty}{\longrightarrow}t^{\alpha}$ and $f$ is bounded on each interval $(0,M]$. Then, by applying Theorem 1.5.2 p.22 in~\cite{BinGolTeu}, $\displaystyle\frac{f(tx)}{f(x)}\underset{x\to\infty}{\longrightarrow}t^{\alpha}$ uniformly in $t$ on each $(0,b]$, $0<b<\infty$.

We can reformulate the dynamics (\ref{dynamic})-(\ref{ConstructX}) into a recursive stochastic algorithm like in the Section~\ref{sec:sa}, and we obtain the following recursive procedure satisfied by the sequence $(\widetilde{Y}_n)_{n\geq0}$:
\begin{equation}\label{ASNLReg}
	\widetilde{Y}_{n+1}=\widetilde{Y}_n-\gamma_{n+1}
	\left(\widetilde{Y}_n-H\frac{\widetilde{Y}_n^{\alpha}}{\Tr(\widetilde{Y}_n^{\alpha})}\right)+\gamma_{n+1}
\left(\Delta M_{n+1}+\widehat{r}_{n+1}\right),
\end{equation}	                   
with the step $\gamma_n =\frac{1}{n+\Tr(Y_0)}$, $\widetilde{Y}_n^{\alpha}=\left((\widetilde{Y}_n^i)^{\alpha}\right)_{1\leq i\leq d}$   and an $\F_n$-measurable  remainder term given by
\begin{equation}\label{resteReg}
	\widehat{r}_{n+1}:=H_{n+1}\frac{\widetilde{f}(Y_n)}{\Tr(\widetilde{f}(Y_n))}-H\frac{\widetilde{Y}_n^{\alpha}}{\Tr(\widetilde{Y}_n^{\alpha})}.
\end{equation}
Notice that, in the skewing case, the remainder term was $r_{n+1}=(H_{n+1}-H)\frac{\widetilde{f}(\widetilde{Y}_n)}{\Tr(\widetilde{f}(\widetilde{Y}_n))}$, therefore assumption {\bf (A3)} implied directly that $r_n\overset{a.s.}{\underset{n\to+\infty}{\longrightarrow}}0$. Here we have to use the uniform convergence of the regular variation to prove the required assumption on $\widehat{r}_{n+1}$. 

By the same arguments used  in Section~\ref{sec:sa}, $\Tr(Y_n)$ satisfies~(\ref{ASTrace}). Moreover, for the quantity $\widetilde{N}_n:=\frac 1n\sum_{k=1}^nX_k$, we also devise a stochastic recursive procedure in the same way as before, namely
$$\widetilde{N}_{n+1}=\widetilde{N}_n-\frac{1}{n+1}\left(\widetilde{N}_n-\frac{\widetilde{Y}_n^{\alpha}}{\Tr(\widetilde{Y}_n^{\alpha})}\right)+\frac{1}{n+1}\left(\Delta\widetilde{M}_{n+1}+\widetilde{r}_{n+1}\right),$$
where $\widetilde{r}_{n+1}=\displaystyle\frac{\widetilde{f}(Y_n)}{\Tr(\widetilde{f}(Y_n))}-\frac{\widetilde{Y}_n^{\alpha}}{\Tr(\widetilde{Y}_n^{\alpha})}$, thus $\widetilde{r}_{n+1}\in\F_n$. 

\begin{theo} 
Let  $d=2$. Assume that {\bf (A1)}, {\bf (A2)} and {\bf (A3)} hold. 
	\item  $(a)$ If $0<\alpha\leq1$, then $h$ has a unique zero $y^*\in I^*$ and 
		$$\frac{\Tr(Y_n)}{n+\Tr(Y_0)}\overset{a.s.}{\underset{n\rightarrow+\infty}{\longrightarrow}}1, \quad \frac{Y_n}{\Tr(Y_n)}\overset{a.s.}{\underset{n\rightarrow+\infty}{\longrightarrow}}y^* \quad\mbox{and}\quad\widetilde{N}_n\overset{a.s.}{\underset{n\rightarrow+\infty}{\longrightarrow}}\frac{(y^*)^{\alpha}}{\Tr((y^*)^{\alpha})}.$$
	\item $(b)$  If $\alpha>1$, then $h$ has a unique zero $y^*\in I^*$ or $ODE_h$ has two attracting equilibrium points in $I^*$ (as we have established in Section~\ref{sec:bound}). Thus, the stochastic recursive procedure $a.s.$ converges to one of the possible limit values.
\end{theo}

\noindent{\bf Proof.} By the same arguments like in Section~\ref{sec:bound}, $\Tr(Y_n)$ satisfies~(\ref{ASTrace}), therefore Proposition~\ref{ThmCvxTrace} holds. Consequently, $\widetilde{Y}_n$ lies in a compact of $\R_+$, thus
$$\max_{1\leq i\leq d}\left|\frac{f(Y_n^i)}{f(n+\Tr(Y_0))}-\left(\frac{Y_n}{n+\Tr(Y_0)}\right)^{\alpha}\right|\underset{n\to+\infty}{\longrightarrow}0.$$
Set $a^i_n=\frac{f(Y_n^i)}{f(n+\Tr(Y_0))}$ and $b^i_n=(\widetilde{Y}_n^i)^{\alpha}$, $i\in\{1,\ldots,d\}$. Then, for every $i\in\{1,\ldots,d\}$,
$$\frac{a^i_n}{\Tr(a_n)}-\frac{b^i_n}{\Tr(b_n)}=\frac{a^i_n-b^i_n}{\Tr(b_n)}+\frac{a^i_n}{\Tr(a_n)}\left(1-\frac{\Tr(a_n)}{\Tr(b_n)}\right).$$
But
$$\Tr(b_n)=\sum_{i=1}^d(\widetilde{Y}_n^i)^{\alpha}\geq\left\{
\begin{array}{rcl}
\left(\sum_{i=1}^d\widetilde{Y}_n^i\right)^{\alpha}=\Tr(\widetilde{Y}_n)^{\alpha} & \mbox{if} & \alpha\in[0,1] \\
d^{1-\alpha}\Tr(\widetilde{Y}_n)^{\alpha} & \mbox{if} & \alpha>1 \\
\end{array}\right.,$$
therefore
$$\Tr(b_n)\geq\frac{\Tr(\widetilde{Y}_n)^{\alpha}}{d^{(\alpha-1)_+}}\underset{a.s.}{\sim}\frac{(n+\Tr(\widetilde{Y}_0))^{\alpha}}{d^{(\alpha-1)_+}}.$$
Consequently, for every $i\in\{1,\ldots,d\}$,
$$\frac{a^i_n}{\Tr(a_n)}-\frac{b^i_n}{\Tr(b_n)}\leq\frac{\max_{1\leq i\leq d}|a^i_n-b^i_n|+\sum_{j=1}^d|a^j_n-b^j_n|}{\Tr(b_n)},$$
$i.e.$
$$\max_{1\leq i\leq d}\left|\frac{a^i_n}{\Tr(a_n)}-\frac{b^i_n}{\Tr(b_n)}\right|\leq\frac{d+1}{\Tr(b_n)}\max_{1\leq i\leq d}|a^i_n-b^i_n|\overset{a.s.}{\underset{n\to+\infty}{\longrightarrow}}0.$$
Thus
$$|\widehat{r}_{n+1}|\leq|\!|\!|H|\!|\!|\max_{1\leq i\leq d}\left|\frac{a^i_n}{\Tr(a_n)}-\frac{b^i_n}{\Tr(b_n)}\right|+|\!|\!|H_{n+1}-H|\!|\!|\overset{a.s.}{\underset{n\to+\infty}{\longrightarrow}}0,$$
and in the same way $\widetilde{r}_{n+1}\overset{a.s.}{\underset{n\to+\infty}{\longrightarrow}}0$. Consequently claim~$(a)$ follows from Proposition~\ref{PropEq}$(a)$ and Theorem~\ref{ThmCvx}.                          

We have to check the assumption on the remainder term to apply result on traps for $SA$. We have that
\begin{equation}\label{VitesseReste}
	\max_{1\leq i\leq d}\left|\frac{a^i_n}{\Tr(a_n)}-\frac{b^i_n}{\Tr(b_n)}\right|\lessim\frac{(d+1)d^{(\alpha-1)_+}}{(n+\Tr(Y_0))^{\alpha}}\max_{1\leq i\leq d}|a^i_n-b^i_n|=o(n^{-\alpha}).
\end{equation}	
So, $\alpha>1$, under assumption {\bf (A3)} on the generating matrices, 
$$\sum_{n\geq0}\left\|r_{n+1}\right\|^2<+\infty.$$
The end of the proof follows from Proposition~\ref{PropEq}$(b)$\&$(c)$ and Theorem~\ref{ThmCvx}.\hfill$\cqfd$

\bigskip
\noindent To establish a $CLT$ for the sequence $(\widetilde{Y}_n)_{n\geq0}$ we need that the remainder term $(r_n)_{n\geq1}$ satisfies~(\ref{HypReste}). Then we will assume that the addition rule matrices $(D_n)_{n\geq1}$ satisfy {\bf (A1)}-$(ii)$ to ensure that $(\widetilde{Y}_n)_{n\geq0}$ lies in the simplex (which implies that the rate in~(\ref{VitesseReste}) is no more $a.s.$) and we assume also that $\alpha>1/2$.

\begin{theo}
Assume that the index of regular variation $\alpha>1/2$, that the addition rule matrices $(D_n)_{n\geq1}$ satisfy {\bf (A1)}-$(ii)$, {\bf (A3)}, {\bf (A4)}  and that {\bf (A5)}$_v$ holds. We have
$$
\mbox{Sp}\left(J_h(y^*)\right)=\left\{1,1-\rho^*\right\},
$$
where
$$
\rho^*=\frac{\alpha(y^{*1})^{\alpha-1}(p_1-y^{*1})+\alpha(1-y^{*1})^{\alpha-1}(y^{*1}-1+p_2)}{(y^{*1})^{\alpha}+f(1-y^{*1})^{\alpha}}.
$$
$(a)$ If $p_1+p_2-1\leq0$ and $v_n=1$, $n\geq1$, then
	$$\sqrt{n}\left(\widetilde{Y}_n-y^*\right)\overset{\L_{stably}}{\underset{n\rightarrow+\infty}{\longrightarrow}}{\cal N}\left(0, \Sigma\right)\quad \mbox{ with }\quad \Sigma=\int_0^{+\infty}e^{u(J_h(y^*)-\frac{I}{2})}\Gamma e^{u(J_h(y^*)-\frac{I}{2})^t}du
	$$
\begin{equation}\label{Gamma2}
\mbox{and  }\quad\Gamma=\frac{(y^{*1})^{\alpha}C^1+(1-y^{*1})^{\alpha}C^2}{\Tr((y^*)^{\alpha})}-y^*(y^*)^t=a.s.\mbox{-}\lim_{n\rightarrow+\infty}\E\left[\Delta M_n\Delta M_n^t\,|\,\F_{n-1}\right].
\end{equation}
$(b)$ If $p_1+p_2-1>0$, we have three possible rate of convergence depending on the second eigenvalue:
		\item[(i)] If $0<\rho^*<\frac12$ and $v_n=1$, $n\geq1$, then 
			$$\sqrt{n}\left(\widetilde{Y}_n-y^*\right)\overset{\L_{stably}}{\underset{n\rightarrow+\infty}{\longrightarrow}}{\cal N}\left(0, \Sigma\right).$$
		\item[(ii)] If $\rho^*=\frac12$ and $v_n =\log n$, $n\ge 1$, then 
			$$\sqrt{\frac{n}{\log n}}\left(\widetilde{Y}_n-y^*\right)\overset{\L}{\underset{n\rightarrow+\infty}{\longrightarrow}}{\cal N}\left(0,\Sigma\right)\; \mbox{ where }\; \Sigma = \lim_n \frac{1}{\log n} \int_0^{\log n} e^{u(J_h(y^*)-\frac{I}{2})}\Gamma e^{u(J_h(y^*)-\frac{I}{2})^t}du.
			$$
		\item[(iii)] If $\frac12<\rho^*<1$ and $v_n=n^{1-2\rho^*+\eta}$, $\eta>0$, then $n^{\rho^*}\left(\widetilde{Y}_n-y^*\right)$ $a.s.$ converges as $n\to+\infty$ towards a positive finite random variable $\Upsilon$.
\end{theo}
This result follows from Theorem~\ref{ThmCLT} and Theorem~\ref{Thm2}.

\subsection{An application to Finance: Adaptive asset allocation}

Such urn based recursive procedures can be applied to adaptive portfolio allocation by an asset manager or a trader, or to optimal split across liquidity pools. Indeed the first setting has already been done in~\cite{LamPagTar} and successfully implemented with multi-armed bandit procedure. We develop in this section the adaptive portfolio allocation, but the optimal split across liquidity pools can be implemented in the same way, by considering that the different colors represent the different liquidity pools, and the trader want to optimally split a large volume of a single asset among the different possible destinations. 

\smallskip
Imagine an asset manager who deals with a portfolio of $d$ tradable assets. To optimize the yield of her portfolio, she can modify the proportions invested in each asset. She starts with the initial allocation vector $Y_0$. At stage $n$, she chooses a tradable asset according to the distribution~(\ref{ConstructX}) or~(\ref{ConstructX2}) of $X_n$, then evaluates its performance over one time step and modifies the portfolio composition accordingly (most likely virtually) and proceeds. Thus the normalized urn composition $\widetilde{Y}_n$ represents the allocation vector among the assets and the addition rule matrices $D_n$ model the successive reallocations depending on the past performances of the different assets. The evaluation of the asset performances can be carried out recursively with an estimator like with multi-arm clinical trials (see~\cite{BaiHuShe,LarPag}). In practice, it can be used to design the addition rule matrices $D_n$. For example, we may consider $(T_n^i)_{n\geq1}$, $i\in\{1,\ldots,d\}$, a {\em success indicator}, namely $d$ independent sequences of i.i.d. $\{0,1\}$-valued Bernoulli trials with respective parameter $p_i$ (with convention $T^i_n=1$ if the return of the $i^{th}$ asset in the $n^{th}$ reallocation is positive and $T^i_n=0$ otherwise).

Let $N^i_n:=\sum_{k=1}^nX^i_k$ be the number of times the $i^{th}$ asset is selected among the first $n$ stages with $N_0^i=0$, $i\in\{1,\ldots,d\}$, and let $S_n$ be the $d$ dimensional vector defined by
$$
S^i_n=S^i_{n-1}+T^i_nX^i_n, \quad n\geq1, \quad S_0^i=1, \quad i\in\{1,\ldots,d\},
$$
denoting the number of successes of the $i^{th}$ asset among these $N^i_n$ reallocations. Define $\Pi_n$ an estimator of the vector of success probabilities, namely $\Pi_n^i=\frac{S_n^i}{N_n^i}$, $i\in\{1,\ldots,d\}$. We can prove that $\Pi_n\overset{a.s.}{\underset{n\to+\infty}{\longrightarrow}}p:=(p^1,\ldots,p^d)^t$ (see~\cite{BaiHuShe,LarPag}). Then we build the following addition rule matrices
\begin{equation}\label{AddMatrixMACT}
D_{n+1}=\begin{pmatrix}T_{n+1}^1 & \frac{\Pi_n^1(1-T_{n+1}^2)}{\sum_{j\neq2}\Pi_n^j} & \cdots & \frac{\Pi_n^1(1-T_{n+1}^d)}{\sum_{j\neq d}\Pi_n^j} \cr  &  &  &  \cr \frac{\Pi_n^2(1-T_{n+1}^1)}{\sum_{j\neq1}\Pi_n^j} & T_{n+1}^2 & \cdots & \frac{\Pi_n^2(1-T_{n+1}^d)}{\sum_{j\neq d}\Pi_n^j} \cr \vdots & \vdots & \ddots & \vdots \cr \frac{\Pi_n^d(1-T_{n+1}^1)}{\sum_{j\neq1}^d\Pi_n^j} & \frac{\Pi_n^d(1-T_{n+1}^2)}{\sum_{j\neq2}^d\Pi_n^j} & \cdots & T_{n+1}^d \end{pmatrix},
\end{equation}
{\it i.e.} at stage $n+1$, if the return of the $j^{th}$ asset is positive, then one ball of type $j$ is added in the urn. Otherwise, $\frac{\Pi_n^i}{\sum_{k\neq j}\Pi_n^k}$ (virtual) balls of type $i$, $i\neq j$, are added. This addition rule matrix clearly satisfies {\bf (A1)}-$(i)$ and {\bf (A2)}. Then, one easily checks that the generating matrices $H_n= \E\,[D_{n+1}\,|\, \F_n]$ satisfy {\bf (A1)}-$(ii)$ and, as soon as $Y_0\in\R_+^d\setminus\{0\}$, $H_n\overset{a.s.}{\longrightarrow}H$ (see~\cite{BaiHuShe,LarPag}), where 
$$
H_{n+1}=\begin{pmatrix}p_1 & \frac{\Pi_n^1(1-p_2)}{\sum_{j\neq2}\Pi_n^j} & \cdots & \frac{\Pi_n^1(1-p_d)}{\sum_{j\neq d}\Pi_n^j} \cr  &  &  &  \cr \frac{\Pi_n^2(1-p_1)}{\sum_{j\neq1}\Pi_n^j} & p_2 & \cdots & \frac{\Pi_n^2(1-p_d)}{\sum_{j\neq d}\Pi_n^j} \cr \vdots & \vdots & \ddots & \vdots \cr \frac{\Pi_n^d(1-p_1)}{\sum_{j\neq1}\Pi_n^j} & \frac{\Pi_n^d(1-p_2)}{\sum_{j\neq2}\Pi_n^j} & \cdots & p_d \end{pmatrix},
\quad
H=\begin{pmatrix}p^1 & \frac{p^1(1-p^2)}{\sum_{j\neq2}p^j}  & \cdots & \frac{p^1(1-p^d)}{\sum_{j\neq d}p^j} \cr  &  &  & \cr \frac{p^2(1-p^1)}{\sum_{j\neq1}p^j} & p^2 & \cdots &\frac{p^2(1-p^d)}{\sum_{j\neq d}p^j} \cr \vdots & \vdots & \ddots & \vdots \cr \frac{p^d(1-p^1)}{\sum_{j\neq1}p^j} & \frac{p^d(1-p^2)}{\sum_{j\neq2}p^j} & \cdots & p^d\end{pmatrix}.
$$
Let remark that $H$ is $\R$-diagonalizable since it is symmetric with respect to its invariant measure. Therefore, the number of each asset in the portfolio $Y_n$ follows the dynamics~(\ref{dynamic}) and the distribution of the portfolio in each asset follows the dynamics~(\ref{ASNLConvex}) or~(\ref{ASNLReg}) depending on the drawing rule. 

Here the components of the limiting generating matrix $H$ can be interpreted as constraints on the composition of the portfolio. Indeed, in presence of two assets (or colors), we prove that the first component of the allocation vector $y^{*1}$ lies in $I^*$ (see Proposition~\ref{PropEq}), therefore the portfolio will contain at least $p_1\vee(1-p_2)$\% and no more than $p_1\wedge(1-p_2)\%$ of the first asset. Such rules may be prescribed by the regulation, the bank policy or the bank customer, and our approach is a natural way to have them satisfied (at least asymptotically). 

The idea of reinforcing the drawing rule (instead of considering the uniform drawing) like in~(\ref{ConstructX}) or~(\ref{ConstructX2}) can be interpreted as a way to take into account the risk aversion of the trader or the customer. Indeed, if $f$ is concave the equilibrium point will be in the middle of the simplex (see Theorem~\ref{ThmAttract} and Theorem~\ref{ThmRepuls}), so the trader prefers to have diversification in her portfolio. On the contrary, if $f$ is convex, the equilibrium points will lie on the boundary of the set of constraints induced by the limiting generating matrix $H$, so she prefers to take advantage of the most money-making asset (like in a ``winner take all'' or a ``0-1'' strategy). 

\bigskip
\noindent {\bf Numerical experiments.} We present some numerical experiments for the drawing rule defined by~(\ref{ConstructX}), firstly with a concave function $f:y\mapsto\sqrt{y}$ and secondly with a convex function $f:y\mapsto y^4$. Therefore we have a unique equilibrium point in the first setting and two attracting targets in the second framework. We consider an asset manager who deal with a portfolio of 2 tradable assets.  We model the addition rule matrices like in the multi-arm clinical trials, namely $D_n$ is defined by~(\ref{AddMatrixMACT}). We use the same success probabilities, namely $p_1=0.7$ and $p_2=0.75$, and the initial urn composition is chosen randomly in the simplex ${\cal S}_2$.

\clearpage
\medskip
\noindent$\rhd$ {\em Convergence of the portfolio allocation with concave drawing rule.}\\

\begin{figure}[!ht]
\centering
\includegraphics[width=15cm]{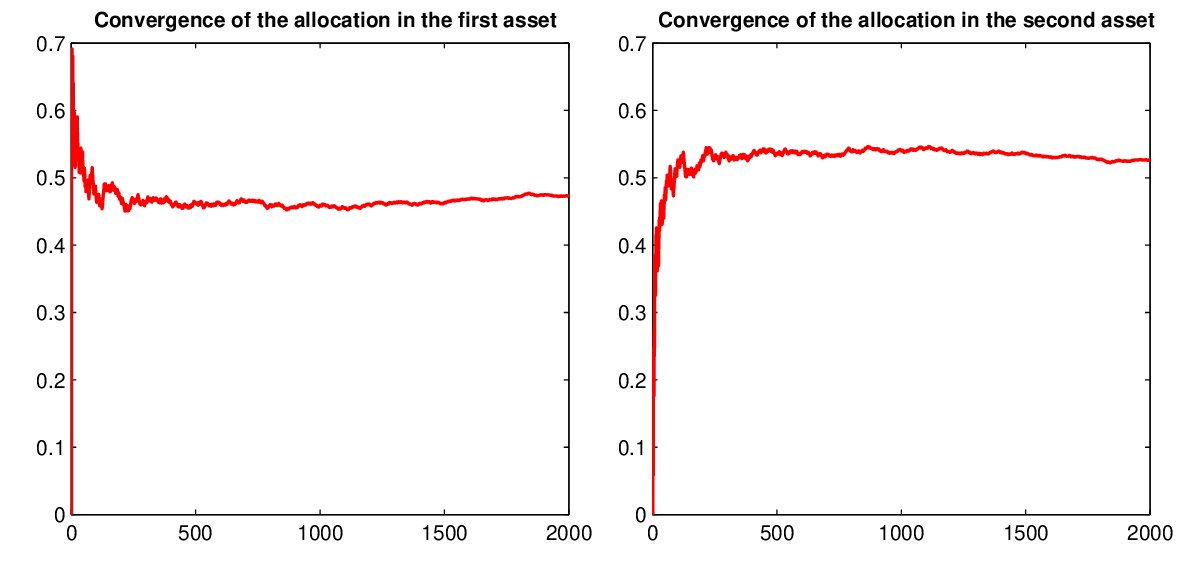}
\caption{Convergence of $\widetilde{Y}_n$ toward $y^*$ for $f(y)=\sqrt{y}$ with $p_1=0.7$ and $p_2=0.75$.}
\end{figure}
We have that $y^{*1}\in(0.25,0.7)$ and $y^{*1}$ and $y^{*2}$ are close to $\frac12$, so the portfolio is diversified because in this case the investor is risk adverse.

\bigskip

\noindent$\rhd$ {\em Convergence of the portfolio allocation with convex drawing rule.}
\begin{figure}[!ht]
\centering
\includegraphics[width=15cm]{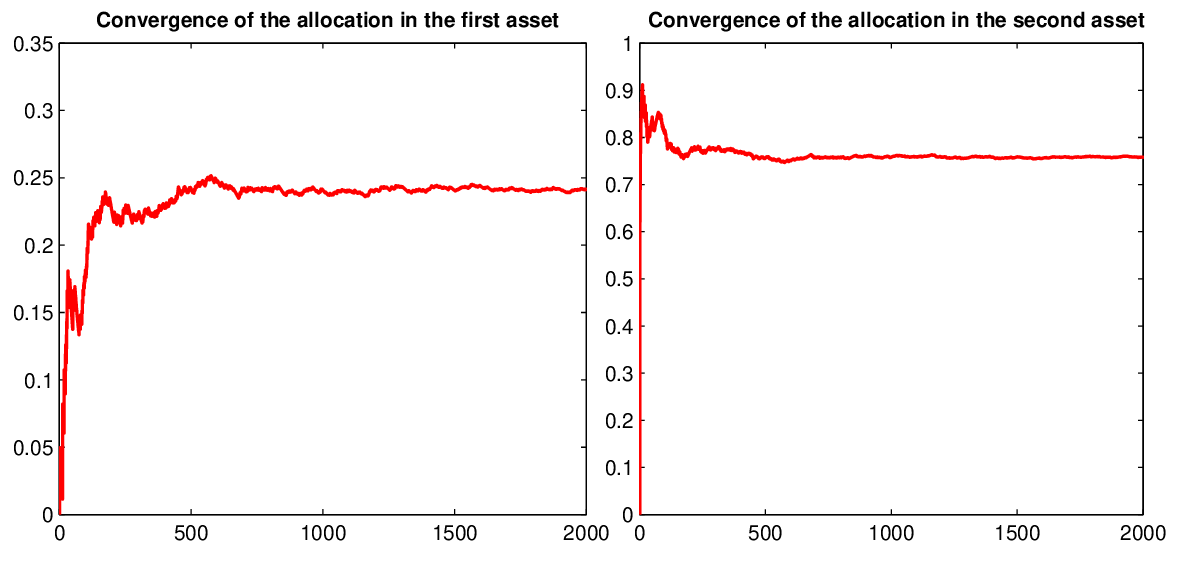}
\caption{Convergence of $\widetilde{Y}_n$ toward $y^*$ for $f(y)=y^4$ with $p_1=0.7$ and $p_2=0.75$.}
\end{figure}


\begin{figure}[!t]
\centering
\includegraphics[width=15cm]{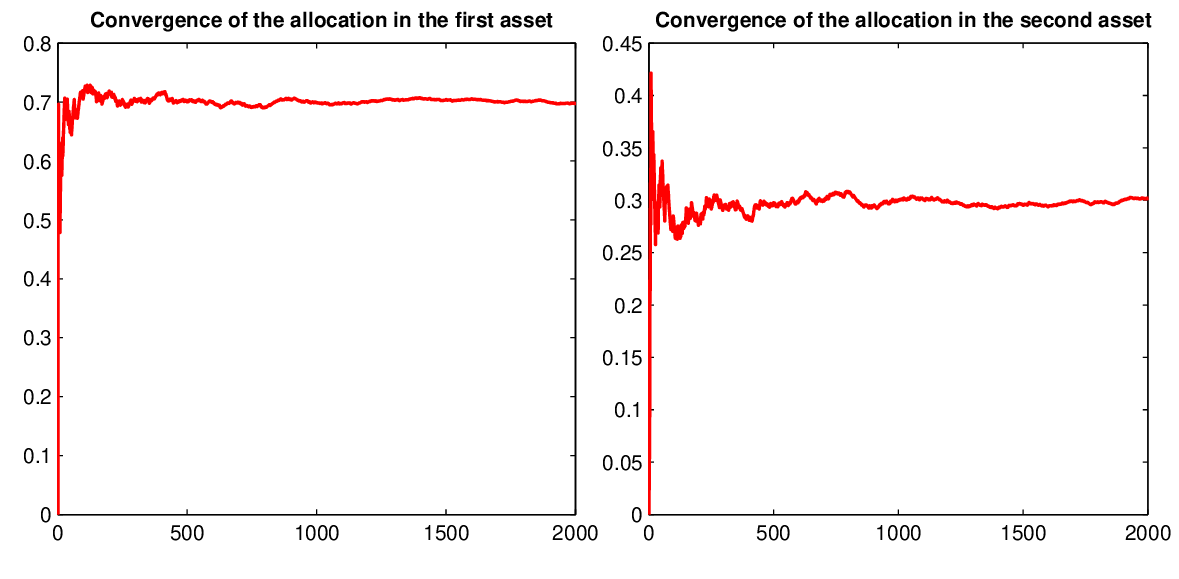}
\caption{Convergence of $\widetilde{Y}_n$ toward $y^*$ for $f(y)=y^4$ with $p_1=0.7$ and $p_2=0.75$.}
\end{figure}

In the convex framework, we have two possible strategies and they are close to the boundaries defined by regulation. Moreover the distribution of the portfolio between the two assets is more asymmetric, because the trader chooses to invest two times more in one asset than in the other.

\paragraph{Acknowledgement.} We thank Fr\'ed\'eric Abergel (MAS Laboratory, ECP)suggesting us to investigate the nonlinear case in randomized urn models and for helpful discussions in view of application to financial frameworks.

\small
\bibliographystyle{plain}
\bibliography{biblio}

\appendix
\begin{center}
\huge{{\bf  Appendix}}
\end{center}

\section{Basic tools from  Stochastic Approximation}
Consider the following recursive procedure defined on a filtered probability space $(\Omega,{\cal A},(\F_n)_{n\geq0},\P)$ having values in a convex set $C\subset \R^d$,
\begin{equation}\label{SAP}
\forall\, n\ge 0,\quad \theta_{n+1}=\theta_n-\gamma_{n+1}h(\theta_n)+\gamma_{n+1}\left(\Delta M_{n+1}+r_{n+1}\right),
\end{equation}
where $(\gamma_n)_{n\ge 1}$ is a $(0, \bar \gamma]$-valued step sequence for some $\bar \gamma>0$, $h:C\rightarrow\R^d$  is a
continuous  function with linear growth (the {\em mean field} of the algorithm) such that  
\begin{equation}\label{eq:Stabh}(I_d-\gamma h)(C)\subset C \;\mbox{ for every }\quad \gamma\!\in (0, \bar \gamma],
\end{equation} 
and $\theta_{0}$ is an $\F_{0}$-measurable finite random vector and, for every $n\ge 1$,  $\Delta M_{n}$ is an $\F_n$-martingale increment and $r_{n}$ is an $\F_n$-adapted remainder term. 

\medskip Note that the assumptions of the theorems recalled below are possibly not minimal, but adapted to the problems we want to solve.

\paragraph{$\rhd$ $A.s.$ Convergence.} Let us introduce a few additional notions on differential systems. We consider the differential system $ODE_h\equiv \dot x=-h(x)$ associated to the  (continuous)   {\em mean field} $h:C\to \R^d$. We assume that this system has a $C$-valued {\em flow}~(\footnote{this is satisfied   under the above stability condition~\eqref{eq:Stabh}.})  $\Phi(t,\xi)_{t\in \R_+, \xi \in C}$: For every $\xi\!\in C$, $(\Phi(t,\xi))_{t\ge 0}$ is the unique solution to $ODE_h$ defined on the whole positive real line. This flow exists as soon as $h$ is locally Lipschitz with linear growth.


Let $K$ be a compact connected, flow invariant subset of $C$, {\it i.e.} such that $\Phi(t,K)\subset K$ for every $t\!\in \R_+$. 

A non-empty subset $A\subset K$ is an {\em  internal attractor} of $K$ for $ODE_h$ if 
\begin{itemize}
\item[$(i)$] $A\varsubsetneq K$,
\item[$(ii)$] $\exists\, \varepsilon_0>0 $ such that $\displaystyle \sup_{x \in K, {\rm dist}(x,A)\le \varepsilon_0}{\rm dist}\big(\Phi(t,x),A\big) \to 0$ as $t\to+\infty$.
\end{itemize}

A compact connected flow invariant set $K$ is a {\em minimal attractor} for $ODE_h$ if it contains no internal attractor. This terminology coming from dynamical systems may be misleading: Thus any equilibrium point of $ODE_h$ (zero of $h$) is a minimal attractor by this definition, regardless of its stability (see~Claim~$(b)$ in Theorem~\ref{ThmODE} below, see also Definition~\ref{def:stabeq}).

\medskip
\noindent {\bf Remark.} When the flow does not exist, the  above definition should be understood as follows: One replaces the flow $\Phi(x,\cdot)$ by the family of all solutions of $ODE_h$ starting from $x$ at time $0$ (whose existence follow from Peano's Theorem). For more details on this natural extension, we refer to~\cite{ForPag2} (see Appendix ``the ODE method without flow").   Up to this extension, the theorem below remains true even when uniqueness of solutions of $ODE_h$ fails.

\begin{theo}\label{ThmODE}($A.s.$ convergence with $ODE$ method, see $e.g.$~\cite{BMP, Duf2, KusYin, ForPag, Ben}). Assume that $h:C\to \R^d$ satisfies~\eqref{eq:Stabh} and that $ODE_h$ has a  $C$-valued flow ($e.g.$ because $h$ is a locally Lipschitz function with linear growth). Assume furthermore that 
$$
r_n\overset{a.s.}{\underset{n\rightarrow+\infty}{\longrightarrow}}0 \quad \mbox{and} \quad \sup_{n\geq 0}\E\left[\left\|\Delta M_{n+1}\right\|^2\,|\,\F_n\right]<+\infty \quad a.s.,
$$
and that $(\gamma_n)_{n\geq1}$ is a positive sequence satisfying ($\gamma_n \!\in (0, \bar \gamma]$, $n\ge 1$)  and 
$$
\sum_{n\geq1}\gamma_n=+\infty \quad \mbox{and} \quad \sum_{n\geq1}\gamma_n^2<+\infty.
$$
On the event $A_{\infty}=\big\{\omega: (h(\theta_n(\omega)))_{n\ge 0} \;\mbox{ is  bounded}\big\}$, $\P(d\omega)$-$a.s.$, the set $\Theta^{\infty}(\omega)$ of the limiting values of $(\theta_n(\omega)_{n\ge 0})$ as $n\rightarrow+\infty$ is  a compact connected flow invariant minimal attractor for $ODE_h$ (see Proposition 5.3 in Section 5.1 in~\cite{Ben}).

\medskip
\noindent Furthermore:

\medskip
\noindent $(a)$ {\em Equilibrium point(s) as limiting value(s)}. If ${\rm dist}\big(\Phi(\theta_0,t),\{h=0\}\big)\to 0$ as $t\to +\infty$, for every $\theta_0\!\in \R^d$,  then $\Theta^{\infty}(\omega)\cap \{h=0\} \neq \varnothing$. 

\medskip
\noindent $(b)$ {\em Single stable equilibrium point}.  If $\{h=0\}= \{\theta^*\}$ and $\Phi(\theta_0,t)\to \theta^*$ as $t\to +\infty$ locally uniformly in $\theta_0$, then $\Theta^{\infty}(\omega)=\{\theta^*\}$ {\it i.e.} $\displaystyle \theta_n \stackrel{a.s.}{\longrightarrow} \theta^*$ as $n\to+\infty$. 

\medskip
\noindent $(c)$ {\em $1$-dimensional setting}. If $d=1$ and $\{h=0\}$ is locally finite, then $\Theta^{\infty}(\omega)=\{\theta_{\infty}\}\subset\{h=0\}$ {\it i.e.} $\displaystyle \theta_n \stackrel{a.s.}{\longrightarrow} \theta_{\infty}\!\in \{h=0\}$.
 
%
%

\end{theo}
A stochastic algorithm may $a.s.$ converge under the existence of multiple equilibrium points, typically stochastic gradient or pseudo-descents, but we do not need such results to solve the urn problems under consideration in this paper. We refer to~\cite{Duf2, KusYin,ForPag,Ben}, among others. Note also that examples of situation~$(a)$ where the algorithm $a.s.$ {\em does not converge} are developed  in~\cite{ForPag},~\cite{ForPag2} or~\cite{Ben} (necessarily with $d\ge 2$ owing to Claim~$(c)$).

\paragraph{$\rhd$ Traps (Unstable equilibrium point).} This second theorem deals with ``traps" {\it i.e.} repulsive zeros of the mean function $h$. It shows that, provided such a trap is noisy enough, such a trap cannot be a limiting point of the algorithm. 

\begin{theo}\label{ThmPiege}($A.s.$ non-convergence toward a noisy  trap, see $e.g.$~\cite{BraDuf,Duf}). Assume that $z^*\in\R^d$ is a trap for the stochastic algorithm~(\ref{SAP}), {\it i.e.} \\
	$(i)$ $h(z^*)=0$,\\
	$(ii)$ there exists a neighborhood $V(z^*)$ of $z^*$ in which $h$ is differentiable with a Lipschitz differential,\\
	$(iii)$ the eigenvalue of $J_h(z^*)$ with the lowest real part, denoted by $\lambda_{{\rm min}}$, satisfies $\Re e(\lambda_{{\rm min}})<0$.\\
Assume furthermore that $a.s.$ on $\Gamma(z^*)=\{\theta_n\underset{n\rightarrow+\infty}{\longrightarrow}z^*\}$,
$$
\sum_{n\geq1}\left\|r_n\right\|^2<+\infty \quad \mbox{and} \quad \lim\sup_{n}\E\left[\left\|\Delta M_{n+1}\right\|^2\,|\,\F_n\right]<+\infty.
$$
Let $K_+$ the subspace of $\R^d$ spanned by the eigenvectors whose associated eigenvalues have a non-negative real part and $K_-$ the subset of $\R^d$ spanned by the eigenvectors whose associated eigenvalues have a negative real part (then $\R^d=K_+\oplus K_-$). By setting $\Delta M_{n+1}^{(r)}$ the projection of $\Delta M_{n+1}$ on $K_-$ alongside $K_+$, assume that $a.s.$ on $\Gamma(z^*)$,
\begin{equation}\label{BruitMgle}
\lim\inf_n\E\left[\left\|\Delta M_{n+1}^{(r)}\right\|\Big | \F_n\right]>0.
\end{equation}
Moreover, if the positive sequence $(\gamma_n)_{n\geq1}$ satisfies
$$
\sum_{n\geq1}\gamma_n=+\infty \quad \mbox{and} \quad \sum_{n\geq1}\gamma_n^2<+\infty,
$$
then   $\P(\Gamma(z^*))=0$.
\end{theo}

\paragraph{$\rhd$ Rate(s) of convergence.} We will say that $h$ is $\epsilon$-differentiable ($\epsilon>0$) at $\theta^*$ if
$$
h(\theta)=h(\theta^*)+J_h(\theta^*)(\theta-\theta^*)+o(\left\|\theta-\theta^*\right\|^{1+\epsilon}) \quad\mbox{as}\quad \theta\to\theta^*.
$$

\begin{theo}\label{ThmCLT}(Rate of convergence see~\cite{Duf2} Theorem 3.III.14 p.131~(for the $CLT$ see also $e.g.$~\cite{BMP,KusYin})). Let $\theta^*$ be an equilibrium point of $\{h = 0\}$ and $\{\theta_n \to \theta^*\}$ the convergence event associated to $\theta^*$ (supposed to have a positive probability). Set  the gain parameter sequence $(\g_n)_{n\ge 1}$ as follows
\begin{equation}\label{HypoPas}
\forall n\geq1, \quad \gamma_n=\frac{1}{n}.
\end{equation}
Assume that the function $h$ is differentiable at $\theta^*$ and all the eigenvalues of $J_h(\theta^*)$ have positive real parts. Assume that, for a real number  $\delta>0$,
\begin{equation}\label{HypDM}
	\sup_{n\geq 0}\E\left[\left\|\Delta M_{n+1}\right\|^{2+\delta}\,|\,\F_n\right]<+\infty \, a.s., \quad \E\left[\Delta M_{n+1}\Delta M_{n+1}^t\,|\,\F_n\right]\overset{a.s.}{\underset{n\rightarrow+\infty}{\longrightarrow}}\Gamma^*\quad \mbox{on$\quad\{\theta_n \to \theta^*\}$,}
\end{equation}
where $\Gamma^*\!\in {\cal S}^+(d, \R)$  (deterministic symmetric   positive matrix) and for  an $\varepsilon>0$  and a positive sequence $(v_n)_{n\ge 1}$ (specified below), 
\begin{equation}\label{HypReste}
n\,v_n\E\left[ \left\|r_{n+1}\right\|^2\mathds{1}_{\{\left\|\theta_n-\theta^*\right\|\leq\varepsilon\}}\right]\underset{n\rightarrow+\infty}{\longrightarrow}0.
\end{equation}

Let $\lambda_{{\rm min}}$ denote the eigenvalue of $J_h(\theta^*)$ with the lowest real part and set $\Lambda:=\Re e(\lambda_{{\rm min}})$.

\smallskip
\noindent $(a)$ If $\Lambda>\frac{1}{2}$ and $v_n = 1$, $n\ge 1$, then, the  weak convergence rate  is ruled on the convergence event $\{\theta_n\overset{a.s.}{\longrightarrow}\theta^*\}$ by the following Central Limit Theorem
$$
\sqrt{n}\left(\theta_n-\theta^*\right)\overset{{\cal L}_{stably}}{\underset{n\rightarrow+\infty}{\longrightarrow}}{\cal N}\left(0, \Sigma^*\right) \quad
\mbox{with} \quad \Sigma^*:=\displaystyle\int_0^{+\infty}e^{-u\left(J_h(\theta^*)^t-\frac{I_d}{2}\right)}\Gamma^* e^{-u\left(J_h(\theta^*)-\frac{I_d}{2}\right)}du.
$$
\smallskip \noindent $(b)$ If $\Lambda=\frac{1}{2}$,  $v_n = \log n$, $n\ge 1$, and $h$ is $\epsilon$-differentiable at $\theta^*$, then
\[
\sqrt{\frac{n}{\log n}}\left(\theta_n-\theta^*\right)\overset{{\cal L}_{stably}}{\underset{n\rightarrow+\infty}{\longrightarrow}}{\cal N}(0,\Sigma^*)\quad \mbox{ on $\quad\{\theta_n \to \theta^*\}$},
\]
where $\displaystyle \Sigma^* = \lim_n \frac{1}{n}\int_0^{n} e^{-u\left(J_h(\theta^*)^t-\frac{I_d}{2}\right)}\Gamma e^{-u\left(J_h(\theta^*)-\frac{I_d}{2}\right)}du$.

\smallskip
\noindent $(c)$ If $\Lambda \!\in \big(0,\frac{1}{2}\big)$,  $v_n = n^{2\Lambda-1+\eta}$, $n\ge 1$, for some $\eta>0$, and $h$ is $\epsilon$-differentiable at $\theta^*$, for some $\epsilon >0$, then $n^{\Lambda}\left(\theta_n-\theta^*\right)$ is $a.s.$ bounded  on $\{\theta_n \to \theta^*\}$ as $n\to+\infty$. 

If, moreover, $\Lambda= \lambda_{\min}$ ($\lambda_{\min}$ is real), then
$n^{\Lambda}\left(\theta_n-\theta^*\right)$ $a.s.$ converges as $n\to+\infty$ toward a finite random variable.
\end{theo}

\end{document}